\documentclass[a4paper, 10pt ]{article}

\usepackage[a4paper,left=1cm,right=1cm,top=2.5cm,bottom=2.5cm]{geometry}
\usepackage{hyperref}
\hypersetup{colorlinks,linkcolor={blue},citecolor={blue},urlcolor={black}}
\usepackage{subfig}
\usepackage{pgfplots}
\usepackage{pgfplotstable}
\usepackage{mathtools}
\usepackage{multicol}
\usepackage{comment}
\usepackage{booktabs}
\pgfplotsset{compat=1.5}
\usepackage{amssymb}
\usepackage{url}
\usepackage{bm}

\usepackage{pdfsync}
\usepackage{float}
\usepackage{tabularx}
\usepackage{enumerate}
\usepackage{array}
\usepackage{xspace}
\usepackage{tikz}
\usepackage{tikz-cd}
\usepackage{tikzsymbols}
\usetikzlibrary{calc,trees,positioning,arrows,chains,shapes.geometric,%
    decorations.pathreplacing,decorations.pathmorphing,shapes,%
    matrix,shapes.symbols, decorations.markings, patterns,fit}

\usepackage{siunitx}

\usepackage{authblk}

\usepackage{todonotes}

\usepackage[draft,inline,marginclue]{fixme}
\usepackage{mathrsfs}
\usepackage{color, colortbl}
\usepackage{multirow}

\usepackage{amsthm}
\usepackage{amsmath}
\usepackage{stmaryrd}
\usepackage[ruled,vlined,linesnumbered]{algorithm2e}
\usepackage{graphicx}
\usepackage{booktabs}
\usepackage{cleveref}
\usepackage{arydshln}
\usepackage[inkscapeformat=png]{svg}

\FXRegisterAuthor{dt}{adt}{\color{red}DT}
\FXRegisterAuthor{fr}{afr}{\color{blue}FR}

\newtheorem{theorem}{Theorem}
\newtheorem{lemma}{Lemma}

\newtheorem{definition}{Definition}

\theoremstyle{remark}
\newtheorem{rmk}{Remark}

\newcommand{\x}{\ensuremath{\mathbf{x}}}
\newcommand{\n}{\ensuremath{\mathbf{n}}}
\newcommand{\ub}{\ensuremath{\mathbf{u}}}
\newcommand{\vb}{\ensuremath{\mathbf{v}}}

\newcommand{\Th}{\mathcal{T}_h}
\newcommand{\NDG}{N_h}
\newcommand{\NRB}{r}

\newcommand{\Npar}{N_{\text{par}}}
\newcommand{\RB}{\text{RB}}
\newcommand{\Naff}{N_{\text{aff}}}
\newcommand{\psiRB}{\psi^{\RB}}
\newcommand{\zRB}{z_{\RB}}

\newcommand{\R}{\ensuremath{\mathbb{R}}}
\newcommand{\Eb}{\ensuremath{\mathbf{E}}}
\newcommand{\Hb}{\ensuremath{\mathbf{H}}}

\newcommand{\scpr}[2]{\left\langle {#1},{#2}\right\rangle}

\DeclareMathOperator*{\argmin}{\arg\!\min}
\DeclareMathOperator*{\argmax}{\arg\!\max}

\newcolumntype{C}[1]{>{\centering\arraybackslash}m{#1}}

\definecolor{Gray}{gray}{0.9}

\begin{document}

\title{Friedrichs' systems discretized with the Discontinuous Galerkin method: domain decomposable model order reduction
and Graph Neural Networks approximating vanishing viscosity solutions}

\author[]{Francesco~Romor\footnote{francesco.romor@sissa.it}}
\author[]{Davide Torlo\footnote{davide.torlo@sissa.it}}
\author[]{Gianluigi~Rozza\footnote{gianluigi.rozza@sissa.it}}

\affil{Mathematics Area, mathLab, SISSA, via Bonomea 265, I-34136 Trieste, Italy}

\maketitle

\begin{abstract}
Friedrichs' systems (FS) are symmetric positive linear systems of first-order partial differential equations (PDEs), which provide 
a unified framework for describing various elliptic, parabolic and hyperbolic semi-linear PDEs such as the linearized Euler equations of gas dynamics, the equations of compressible linear elasticity and the Dirac-Klein-Gordon system. FS were studied to approximate PDEs of mixed elliptic and hyperbolic type in the same domain. 
For this and other reasons, the versatility of the discontinuous Galerkin method (DGM) represents the best approximation space for FS. We implement a distributed memory solver for stationary FS in \texttt{deal.II}. 
Our focus is model order reduction. Since FS model hyperbolic PDEs, they often suffer from a slow Kolmogorov n-width decay. We develop two approaches to tackle this problem. The first is domain decomposable reduced-order models (DD-ROMs). We will show that the DGM offers a natural formulation of DD-ROMs, in particular regarding interface penalties, compared to the continuous finite element method. We also develop new repartitioning strategies to obtain more efficient local approximations of the solution manifold. The second approach involves graph neural networks used to infer the limit of a succession of projection-based linear ROMs corresponding to lower viscosity constants: the heuristic behind is to develop a multi-fidelity super-resolution paradigm to mimic the mathematical convergence to vanishing viscosity solutions while exploiting to the most interpretable and certified projection-based ROMs.
\end{abstract}

%\tableofcontents
%\listoffixmes

\section{Introduction}
\label{sec:intro}
%!TEX root = ../main_arxiv.tex

Friedrichs' systems (FS) are a class of symmetric positive linear systems of first-order partial derivative equations (PDEs). 
They were introduced by Friedrichs \cite{friedrichs1958symmetric} as a tool to study hyperbolic and elliptic phenomena in different parts of the domain within a unifying framework.
The main ideas that allow recasting many models into the FS frameworks are the introduction of extra variables to lower the order of the higher derivatives and the linearization of nonlinear problems.
FS are characterized by linear and positive operators and (non-uniquely defined) boundary operators that allow them to impose classical boundary conditions (BCs).
Various works proved uniqueness, existence and well--posedness of the FS in their strong, weak and ultraweak formulation and the necessary conditions to properly define the boundary operators \cite{friedrichs1958symmetric,rauch1985symmetric,rauch1994boundary,ern2006discontinuous,ern2007intrinsic,antonic2009equivalent,di2011mathematical}.

In the last decades, different numerical discretizations of the FS have been proposed to approximate the analytical solutions.
The strategies vary among finite volume \cite{sonar1998dual} and discontinuous Galerkin (DG) formulations \cite{hoang2021domain,jensen2004discontinuous,ern2006discontinuous,Ern2006b,ern2007intrinsic,ern2008discontinuous,bui2013unified,chen2023unified}.
Along with the DG discretization, also error estimation analysis that provide, according to the type of edge penalization, optimal or sub--optimal estimates, have been carried out \cite{ern2006discontinuous,Ern2006b,di2011mathematical}. We focus on the DG method since it is more versatile to approximate both elliptic and hyperbolic PDEs and it fits naturally in the framework of domain decomposable ROMs (DD-ROMs).

In the context of parametric PDEs, for multi--query tasks or real--time simulations, fast and reliable simulations of the same problem for different parameters are often needed. This is especially true when the full-order models (FOMs) are based on expensive and high-order DG discretizations.
Reduced order models (ROMs) decrease the computational costs looking for the solutions of unseen parametric instances on low-dimensional discretization spaces, called reduced basis spaces. This is possible because the new solutions to be predicted are expected to be highly correlated with the database of training DG solutions used to build the reduced spaces.
ROMs have been proven to be a powerful tool for many applications \cite{prud2002reliable,maday2002reduced,hesthaven2016certified,rozza2022advanced}. 
In particular for linear problems, classical Galerkin and Petrov-Galerkin projection methods are very easy to set up and extremely convenient in terms of computational costs. FS are perfectly suited for such algorithms due to their linearity. This is a preliminary step needed to reduce parametric nonlinear PDEs whose linearization results in FS. In the most simple formulation, we will apply singular value decomposition (SVD) to compress a database of snapshots and provide a reliable reduced order model (ROM), with standard \textit{a posteriori} error estimators.

In the context of model order reduction, FS are particularly beneficial as theoretical frameworks for many reasons. They represent a new form of structure-preserving ROMs: the positive symmetric properties of FS are in fact easily inherited by the reduced numerical formulations. This advocates for the employment of FS for reduced order modelling whenever a PDE can be reformulated in the FS framework. This is the case for the Euler equations of gas dynamics, when they are written in terms of entropy variables~\cite{sonar1998dual, parish2022impact}. The same rationale is behind structure-preserving symplectic or Hamiltonian ROMs~\cite{hesthaven2021structure} and port-Hamiltonian ROMs~\cite{van2014port, beattie2022structure}.
Moreover, since FS are often studied in their ultraweak formulation, they are good candidates for optimally stable error estimates~\cite{bui2013constructively} at the full-order level~\cite{bui2013unified}, also in a hybridized DG implementation in~\cite{chen2023unified}, and at the reduced order level, similarly to what has been achieved in the works~\cite{beurer2022ultraweak, hain2022ultra, henning2022ultraweak}. Finally, from the point of view of software design, the possibility to implement in a unique maintainable and generic manner the realization of ROMs for PDEs ascribed to the class of FS is a convenient feature to search for.

Though being linear, FS are hyperbolic systems and often show an advection dominated character, which is not easily approximable through a simple proper orthogonal decomposition (POD). This leads to a slow Kolmogorov $n$-width (KnW) decay that results in very inefficient approximations of the reduced models.
Several approaches have been studied to overcome this difficulty \cite{taddei19registration,iollo14advectionModes,peherstorfer18adaptiveBases,rim18displacementInterpolation,carlberg11,Cagniart2019,amsallem08interpolation,lee2020model,crisovan2019model,torlo20ALEMOR,iollo2022mapping}.

%% DD 
A strategy that has been developed to reduce PDEs solved numerically with domain decomposition approaches, like fluid-structure interaction systems, are domain decomposable ROMS (DD-ROMS). The initial formulations~\cite{maday2002reduced, maday2004reduced, huynh2013static, eftang2012adaptive, huynh2013static} involved continuous finite elements discretizations for which new ways to couple the solutions restricted to different subdomains needed to be devised, especially to enforce continuity at the interfaces. 
We show that the DGM imposes naturally flux interface penalties from the full-order discretizations and it is, thus, amenable for straightforward implementations of DD-ROMs. 
From the point of view of solution manifold approximability and so KnW decay, DD-ROMs are based on local linear approximants that are employed to reach a higher accuracy for unseen solutions. This is useful when the computational domain is divided in subdomains that are independently affected by the parametric instances. The typical case in which this may happen are parametric models for which discontinuous values of the parameters over fixed subdomains cause non correlated responses on their respective subdomains. Similar cases will be studied in~\cref{subsubsec:max,subsubsection:cle_test}. Another example is represented by parametric fluid-structure interaction systems in which the parameters cause complex interdependencies between the structure and fluid components in favor of partitioned linear solution manifold approximations (SVD is performed separately for the fluid, for the structure and for the interface) rather than monolithic ones. 
In our implementation of DD-ROMs, we exploit the partitions obtained from the distributed memory solver in \texttt{deal.II}. 
Since these domain decompositions typically satisfy constraints related to the computational efficiency, we devise some strategies to repartition the domain responding to solution manifold approximability concerns instead. 
Another work that implements this is~\cite{{xiao2019domainb}} where the Reynolds stress tensor is employed, among others, as indicator for partitioning the computational domain. 
Similarly, we develop new indicators.

%% VV 
Another way to approach the problem of a slow KnW decay is exploiting the mathematical proofs of existence of vanishing viscosity solutions \cite{oleinik1963discontinuous,kruvzkov1970first,diperna1982convergence,goodman1992viscous}.
In fact, solutions of hyperbolic problems can be obtained as a limit process of solutions associated to viscosity terms approaching zero.
The crucial point is that ROMs associated to larger viscosity values may not suffer from a slow Kolmogorov $n$-width decay.
Hence, we can set up classical projection based ROMs for the high viscosity solutions, and use graph neural networks (GNNs)~\cite{tencer2021tailored} only to infer the vanishing viscosity solution in a very efficient manner.
This procedure can be applied also to more general hyperbolic problems, not necessarily FS.
The key features of this new methodology are the following: the employment of computationally heavy graph neural networks is reduced to a minimum and, at the same time, interpretable certified projection ROMs are exploited as much as possible in their regime of accurate linear approximability. In fact, GNNs, used generally to perform non-intrusive MOR, have high training computational costs and they are employed mainly for small academic benchmarks in terms of number of degrees of freedom, up to now. We avoid these high computational efforts with our multi-fidelity formulation: the GNNs are employed only to infer the vanishing viscosity solutions from the previous higher viscosity level, not to approximate and perform dimension reduction of the entire solution manifold. The overhead is the collection of additional full-order snapshots corresponding to high viscosity values, but this can be performed on coarser meshes as it will be done in~\cref{sec:vv}. 
Moreover, the support of our GNNs is the DG discretization space, so, we can enrich the typical machine learning framework of GNNs with data structure and operators from numerical analysis. 
We validate the use of data augmentation with numerical filters (discretized Laplacian, gradients), as proposed in~\cite{tencer2021tailored}.

In brief, we summarize our contributions with the present work:
\begin{itemize}
    \item structure-preserving model order reduction for Friedrichs' systems. We synthetically describe the realization of ROMs for FS and the definition of standard \textit{a posteriori} error estimators. Hints towards the implementation of optimally stable ROMs are highlighted.
    \item domain decomposable reduced-order models for full-order models discretized with the discontinuous Galerkin method. We introduce DD-ROMs for DG discretizations and introduce novel indicators to repartition the computational domain with the aim of obtaining more efficient local solution manifold approximants.
    \item surrogate modelling of vanishing viscosity solutions with graph neural networks. We propose a new framework for the MOR of parametric hyperbolic PDEs with a slow Kolmogorov n-width decay.
\end{itemize}

The topics addressed in this work are presented as follows. In Section~\ref{sec:fs}, we introduce the definition of FS and well--posedness results and we will provide several examples of models that fall into this framework: the Maxwell equations in stationary regime, the equations of linear compressible elasticity and the advection diffusion reaction equations. 
Then, we provide a DG discretization of the FS following \cite{di2011mathematical} with related error estimates in Section~\ref{sec:fom}. 
In Section~\ref{sec:mor}, we introduce the projection-based MOR technique and some error bounds that can be effectively used. 
In Section~\ref{sec:results}, we will discuss a new implementation of domain decomposable ROMs for FOMs discretized with the DGM and we will test the approach on three parametric models.
In Section~\ref{sec:vv}, we introduce the concept of vanishing viscosity solutions and how graph neural network are exploited to overcome the problem of a slow Kolmogorov $n$-width decay. We will provide some numerical tests to show the effectiveness of the proposed approach.
Finally, in Section~\ref{sec:conclusions} we summarize our results and we suggest further directions of research.

\section{Friedrichs’ systems}
\label{sec:fs}
In this section, we will provide a summary of FS theory: their definition, existence, uniqueness and well-posedness results, their weak and ultraweak forms and many PDEs which can be rewritten into FS. 
The following discussion collects many results from \cite{friedrichs1958symmetric,rauch1985symmetric,hughes1986new,rauch1994boundary, houston1999posteriori,jensen2004discontinuous,sonar1998dual,jensen2004discontinuous,Ern2006b,ern2007intrinsic,ern2008discontinuous,antonic2009equivalent,bui2013unified,di2011mathematical}, but we will follow the notation in \cite{di2011mathematical}.
Let us represent with $d$ the ambient space dimension and with $m\geq 1$ the number of equations of the FS. We consider a connected Lipschitz domain $\Omega \subset \R^d$, with boundary $\partial\Omega$ and outward unit normal $\mathbf{n}:\partial\Omega\rightarrow\mathbb{R}^d$.

A FS is defined through $(d+1)$ matrix-valued fields $A^0, A^1\dots, A^d\in [L^{\infty}(\Omega) ]^{m \times m}$ and the following differential operators $\mathcal{X}, A, \tilde{A}:\Omega\rightarrow\mathbb{R}^{m \times m}$. We suppose that $\mathcal{X}\in[L^{\infty}(\Omega) ]^{m \times m}$ and define
\begin{equation}
	\mathcal{X}=\sum_{k=1}^d \partial_k A^k\ ,\qquad A= A^0+\sum_{i=1}^d A^i \partial_i \ ,\qquad \tilde{A}= \left({A^0}\right)^t - \mathcal{X} - \sum_{i=1}^d A^i \partial_i \ ,
\end{equation} assuming that
\begin{subequations}
	\begin{align}
		&A^k = (A^k)^T \, \text{a.e. in }\Omega , \qquad \text{for } k=1, \dots, d,&&\text{(symmetry property)}\\
		&A^0+ ({A}^0)^T - \mathcal{X} \,\, \text{is u.p.d.} \quad \text{a.e. in }\Omega,&&\text{(positivity property)}\label{eq:positivity_classical}
	\end{align}
\end{subequations}
thus, the name \textbf{symmetric positive operators} or \textbf{Friedrichs operators}, which is used to refer to ($A$, $\tilde{A}$). We recall that the operator in~\eqref{eq:positivity_classical} is uniformly positive definite (u.p.d) if and only if
\begin{equation}
	\exists\mu>0: A^0+ ({A}^0)^T - \mathcal{X}>2\mu_0\mathbb{I}\quad \text{a.e. in }\Omega.
\end{equation}
If this property is not satisfied, it can be sometimes recovered as shown in Appendix~\ref{appendix:accretive}. A weaker condition can be required for two-field systems~\cite{ern2008discontinuous}.

The boundary conditions are expressed through two boundary operators $\mathcal{D}:\partial \Omega \to \R^{m\times m}$ with
\begin{equation}
	\mathcal{D}=\sum_{k=1}^{d} n_k A^k,\qquad \text{a.e. in }\partial\Omega
\end{equation}
and $\mathcal{M}:\partial \Omega \to \R^{m\times m}$ satisfying the following \textbf{admissible boundary conditions}
\begin{subequations}
	\label{eq: admissible strong form}
\begin{align}
	&\mathcal{M} \quad \text{is nonnegative}\quad\text{a.e. on}\quad\partial\Omega,&&\text{(monotonicity property)}\\
	&\ker( \mathcal{D}-\mathcal{M}) + \ker(\mathcal{D}+\mathcal{M}) =\R^m \quad\text{a.e. on}\quad\partial\Omega.&&\text{(strict adjointness property)}
\end{align}
\end{subequations}

\begin{rmk}[Strict adjointness]
	The term strict adjointness property comes from Jensen \cite[Theorem 31]{jensen2004discontinuous}. The strict adjointness property is needed for the solution of the ultra-weak formulation of the FS to uniquely satisfy the boundary conditions: in a slightly different framework from the one presented here, see \cite[Theorem 29]{jensen2004discontinuous} and \cite[proof of Lemma 2.4]{bui2013unified}.
\end{rmk}

\begin{theorem}[Friedrichs' system strong solution~\cite{friedrichs1958symmetric}]
    \label{def: fs_continuous}
    Let $f\in [L^2(\Omega)]^m$, the strong solution $z\in [C^1(\overline{\Omega})]^m$ to Friedrichs' system
    \begin{equation}\label{eq:FS}
    	\begin{cases}
    		Az=f,&\text{in }\Omega,\\
    		(\mathcal{D}-\mathcal{M})z=0, & \text{on }\partial\Omega.
    	\end{cases}
    \end{equation}% find $z\in V_0$ s.t. $A^0z + \sum_{i=1}^{d}A^i\partial_i z = f$ in $L$.
		is unique. Moreover, there exists a solution of the ultra-weak formulation
		\begin{equation}
			\label{eq:strong uw}
			(z, \tilde{A}y)_{L^2} = (f, y)_{L^2},\qquad \forall y\in [C^1(\overline{\Omega})]^m\ s.t.\ (\mathcal{D}+\mathcal{M}^t)y=0.
		\end{equation}
\end{theorem}
Let $L=[L^2(\Omega)]^m$. We define the weak formulation on the graph space $V=\lbrace z \in L: Az \in L\rbrace$, which amounts to differentiability
in the characteristics directions: $A\in\mathcal{L}(V, L)$ and $\tilde{A}\in\mathcal{L}(V', L)$. The boundary operator $\mathcal{D}$ is translated into the abstract operator $D\in\mathcal{L}(V, V')$:
\begin{equation}
	\label{eq:def_D}
	\left\langle Dz, y\right\rangle _{V, V'} = (Az, y)_{L} - (z, \tilde{A}y)_{L},\quad\forall z, y\in V.
\end{equation}
When $z$ is smooth, it can be seen as the integration by parts formula~\cite{jensen2004discontinuous, bui2013unified}:
\begin{equation}
	\left\langle Dz, y\right\rangle_{V, V'} = \left\langle \mathcal{D}z, y\right\rangle_{H^{\frac{1}{2}}(\partial\Omega), H^{-\frac{1}{2}}(\partial\Omega)},\quad \forall z\in H^1(\Omega),\ y\in H^1(\Omega).
\end{equation}
A sufficient condition for well-posedness of the weak formulation is provided by the cone formalism~\cite{antonic2009equivalent, di2011mathematical} that poses the existence of two linear subdomains $(V_0, V^*_0)$ of $V$:
\begin{subequations}
	\label{eq: cone formalism}
	\begin{align}
		&V_0 \ \text{maximal in}\  C^+, \quad V^*_0\ \text{maximal in}\  C^-\\
		& V_0 = D(V_0^*)^{\perp},\quad V_0^* = D(V_0)^{\perp},
	\end{align}
\end{subequations}
such that $A:V_0\rightarrow L$ and $\tilde{A}:V^*_0\rightarrow L$ are isomorphism, where $C^{\pm} = \{w\in V| \pm \scpr{Dw}{w}_{V, V'}\geq 0\}$.

Provided that $V_0+V_0^*\subset V$ is closed~\cite{antonic2009equivalent}, the conditions in~\eqref{eq: cone formalism} are equivalent to the existence of the boundary operator $M\in\mathcal{L}(V, V')$ that satisfies \textbf{admissible boundary conditions} analogue to the ones in~\eqref{eq: admissible strong form}:
\begin{subequations}
	\label{eq: admissible weak form}
\begin{align}
	&M \quad \text{is monotone},&&\text{(monotonicity property)}\\
	&\ker( D-M) + \ker(D+M) =V,&&\text{(strict adjointness property)}
\end{align}
\end{subequations}
identifying $V_0 = \ker(D-M)$ and $V^*_0 = \ker(D+M^*)$. 

\begin{theorem}[Friedrichs' System weak form~\cite{ern2006discontinuous,Ern2006b,di2011mathematical,bui2013unified}]
    \label{def: fs_weak_continuous}
    Let us assume that the boundary operator $M\in\mathcal{L}(V, V')$ satisfies the monotonicity and strict adjointness properties~\eqref{eq: admissible weak form}. Let us define for $z, z^*\in V$ the bilinear forms
		\begin{subequations}
			\begin{align}
				\label{eq:bilinear}
        a(z, y)
        &= \left( Az, y\right)_L+\tfrac{1}{2}\langle(D-M)z, y\rangle_{V', V}, \quad\forall y\in V,\\
				a^*(z^*, y)
        &= \left(\tilde{A}z^*, y\right)_L+\tfrac{1}{2}\langle(D+M^*)z^*, y\rangle_{V', V}, \quad\forall y\in V.
			\end{align}
		\end{subequations}
    Then, Friedrichs' operators $A:V_0\rightarrow L$ and $\tilde{A}:V^*_0\rightarrow L$ are isomorphisms: for all $f\in L$ and $g\in V$  there exists unique $z, z^*\in V$ s.t.
		\begin{subequations}
			\begin{align}
        a(z, y) &= (f, y)_L + \langle(D-M)g, y\rangle_{V', V} \quad\forall y\in V,\label{eq:weak_continuous_FS}\\
				a^*(z^*, y) &= (f, y)_L + \langle(D+M^*)g, y\rangle_{V', V} \quad\forall y\in V,
			\end{align}
		\end{subequations}
		that is 
		\begin{equation}
			\begin{cases}
				Az = f, & \text{in }L,\\
				(M-D)(z-g) = 0 &\text{in }V',
			\end{cases}
			\qquad
			\begin{cases}
				\tilde{A}z^* = f, & \text{in }L,\\
				(M^*+D)(z^*-g) = 0 &\text{in }V'.
			\end{cases}
		\end{equation}
\end{theorem}

\subsection{A unifying framework}\label{sec:examples}
The theory of Friedrichs' systems provides a unified framework to study different classes of PDEs~\cite{jensen2004discontinuous}: first-order uniformly hyperbolic, second-order uniformly hyperbolic, elliptic and parabolic partial differential equations. Originally, Friedrichs' aim was to study equations of mixed type (hyperbolic, parabolic, elliptic) inside the same domain such as the Tricomi equation~\cite{friedrichs1958symmetric} (or more generally the Frankl equation~\cite{jensen2004discontinuous}) inspired by models from compressible gas dynamics for which the domain is subdivided in a hyperbolic supersonic and an elliptic subsonic part.

Some examples of FS can be found in the literature:
\begin{subequations}
	{\allowdisplaybreaks
	\begin{align}
		(x_2\partial^2_1\bullet + \partial^2_1\bullet)u &= 0, &&(\text{Tricomi~\cite{friedrichs1958symmetric}})\\
		\begin{bmatrix}
			-\partial_1\bullet & \partial_2\bullet\\
			-\partial_2\bullet & \partial_1\bullet
		\end{bmatrix}\begin{pmatrix}
			u_1 \\ u_2
		\end{pmatrix} &=0, &&(\text{Cauchy-Riemann~\cite{friedrichs1958symmetric}})\\
		(A(x_2)\partial^2_1 + \partial^2_1)u &= 0, &&(\text{Frankl~\cite{jensen2004discontinuous}})\\
		\begin{bmatrix}
			\mathbb{I}_3 & - \lambda^{-1}\mathbb{I}_3(\nabla\cdot\bullet) -\frac{\left(\nabla\bullet  + (\nabla \bullet)^t\right)}{2}\\
			-\frac{1}{2}\nabla\cdot\left(\bullet+\bullet^t\right) &  \alpha\mathbb{I}_3
		\end{bmatrix}
		\begin{pmatrix}
			\sigma \\ \mathbf{u}
		\end{pmatrix} &= 0, &&(\text{Compressible linear elasticity~\cite{ern2006discontinuous}})\label{eq:cle}\\
		\begin{bmatrix}
			\mu\mathbb{I}_3 & \nabla \times\bullet\\
			-\nabla \times\bullet &  \sigma\mathbb{I}_3
		\end{bmatrix}
		\begin{pmatrix}
			\mathbf{H} \\ \mathbf{E}
		\end{pmatrix} &= 0, &&(\text{Maxwell eq. in stationary regime~\cite{ern2006discontinuous}})\label{eq:msr}\\
		(-\nabla \cdot (\kappa \nabla \bullet) +\boldsymbol{\beta} \cdot \nabla \bullet +\mu \bullet)\mathbf{u}&= 0, &&(\text{Diffusion advection reaction~\cite{ern2006discontinuous}})\label{eq:dar}\\
		(A_0\partial_t\bullet+\Sigma_{i=1}^{3} \tilde{A}_i\partial_i\bullet)\mathbf{V}&= 0, &&(\text{Linearized symmetric Euler~\cite{hughes1986new,sonar1998dual}})\label{eq:leu}\\
		(a\gamma^0 \partial_t\bullet + \gamma^1\partial_1\bullet+ \gamma^2\partial_2\bullet+ \gamma^3\partial_3\bullet+ B)\boldsymbol{\psi}&= 0, &&(\text{Dirac system~\cite{antonic2017complex}})\label{eq:complexD}\\
		\begin{bmatrix}
			-ai\gamma^0 \partial_t\bullet -i\gamma^1\partial_1\bullet-i\gamma^2\partial_2\bullet-i\gamma^3\partial_3\bullet+ M & \mathbf{1}_4\\
			\mathbf{0}_4 &  \partial_t^2\bullet -\Delta\bullet+m^2\mathbb{I}_1
		\end{bmatrix}
		\begin{pmatrix}
			\boldsymbol{\psi} \\ \phi
		\end{pmatrix} &= \mathbf{f}, &&(\text{Dirac-Klein-Gordon system~\cite{antonic2017complex}})\label{eq:complexDKG}\\
		\begin{bmatrix}
			-\frac{i}{2\pi}(a \gamma^0 \partial_t\bullet +\gamma^1\partial_1\bullet+\gamma^2\partial_2\bullet+\gamma^3\partial_3\bullet+ B) & \mathbf{1}_4\\
			\mathbf{0}_4 &  -\partial_t^2\bullet +\Delta\bullet
		\end{bmatrix}
		\begin{pmatrix}
			\boldsymbol{\psi} \\ \mathcal{A}
		\end{pmatrix} &= \mathbf{f}, &&(\text{Maxwell-Dirac system~\cite{antonic2017complex}})\label{eq:complexMD}\\
		\begin{bmatrix}
			-i\omega\mu\mathbb{I}_3 & \nabla \times\bullet\\
			-\nabla \times\bullet &  (-i\omega\epsilon+\sigma)\mathbb{I}_3
		\end{bmatrix}
		\begin{pmatrix}
			\mathbf{H} \\ \mathbf{E}
		\end{pmatrix} &= 0, &&(\text{Time-harmonic Maxwell~\cite{antonic2017complex}})\label{eq:complexTM}\\
		\begin{bmatrix}
			\nu\mathbb{I}_3 & \nabla \times\bullet\\
			\mu\boldsymbol{\beta}\times\bullet-\nabla \times\bullet &  \sigma\mathbb{I}_3
		\end{bmatrix}
		\begin{pmatrix}
			\mathbf{H} \\ \mathbf{E}
		\end{pmatrix} &= \mathbf{f}, &&(\text{Magneto-hydrodynamics~\cite{Ern2006b}})\\
		\begin{bmatrix}
			\nu^{-1}\mathbb{I}_3 & \mathbf{1}_3 &-\frac{\left(\nabla\bullet  + (\nabla \bullet)^t\right)}{2}\\
			\text{tr}(\bullet) &  d\mathbb{I}_1 & 0\\
			-\frac{1}{2}\nabla\cdot\left(\bullet+\bullet^t\right) & 0 & \boldsymbol{\beta} \cdot \nabla \bullet
		\end{bmatrix}
		\begin{pmatrix}
			\sigma \\ p \\ \mathbf{u}
		\end{pmatrix} &= \mathbf{f}, &&(\text{Incompressible linearized Navier-Stokes~\cite{ern2008discontinuous}})\label{eq:ifs}
	\end{align}
}
	for the employed notation we refer to the respective reported references. A non-stationary version of \eqref{eq:msr} and \eqref{eq:dar} from~\cite{di2011mathematical} is omitted. The FS framework here presented easily extends to complex-valued systems as in \eqref{eq:complexD}, \eqref{eq:complexDKG}, \eqref{eq:complexMD} and \eqref{eq:complexTM} from~\cite{antonic2017complex}. We will consider only semi-linear PDEs but FS can be encountered as intermediate steps when solving quasi-linear PDEs: for example solving the compressible Euler equations of gas dynamics in entropy variables with the Newton method brings to the FS~\eqref{eq:leu}, as studied in~\cite{sonar1998dual}.

	One of the critical points of FS is the definition of the boundary conditions. Friedrichs's idea~\cite{friedrichs1958symmetric} was to impose boundary conditions through a matrix-valued boundary operator. Ern, Guermond and Caplain~\cite{ern2007intrinsic} revised the FS theory, without employing the trace of functions in the graph space as developed in~\cite{rauch1985symmetric, rauch1994boundary, houston1999posteriori, jensen2004discontinuous}, but in terms of operators acting in abstract Hilbert
	spaces, as presented here.
	
	The most common homogeneous boundary conditions (homogeneous Dirichlet, Neumann, Robin) for \eqref{eq:cle}, \eqref{eq:msr} and \eqref{eq:dar} can be found in the literature~\cite{ern2006discontinuous,Ern2006b,ern2007intrinsic}. 
	For a choice of boundary conditions, and thus for a choice of spaces $(V_0, V_0^*)$, there can be more than one definition of the boundary operator $M$ \cite[Remark 5.3]{ern2007intrinsic}. A constructive methodology for defining the boundary operator $M\in\mathcal{L}(V, V')$ from specific boundary conditions can be found in~\cite{ern2007intrinsic} and it will be employed for the compressible linear elasticity test case in Section~\ref{subsubsec:comprLinElast_def}. Also inhomogeneous boundary conditions can be imposed through the definition of traces of functions in graph spaces as in~\cite{jensen2004discontinuous} or through a Petrov-Galerkin formulation as in~\cite{bui2013unified}.
\end{subequations}
In the following, we will present in detail three FS on which we will focus in the numerical test section.

\subsubsection{Curl--curl problem: Maxwell equations in stationary regime}
\label{subsubsec:curl-curl_def}
We will consider the Maxwell equations in stationary regime, also known as the curl-curl problem. Let $\mathbf{E} \in \R^d =\R^3$ be the electric field and $\mathbf{H}\in \R^3$ be the magnetic field. The curl--curl problem is defined as
\begin{align*}
	&\begin{cases}
		{\mu \mathbf{H} }  +  {\nabla \times \mathbf{E}} = \mathbf{r},\\
		{\sigma \mathbf{E}} - {\nabla \times \mathbf{H}} = \mathbf{g},
	\end{cases}
\end{align*}
with $\mu, \sigma >0$, the permeability and permittivity constants. The FS is obtained by setting
\begin{align*}
	& {{A}^0 = \begin{bmatrix}
			\mu \mathbb{I}_{d,d} & 0_{d,d}\\
			0_{d,d} & \sigma \mathbb{I}_{d,d}
	\end{bmatrix}}, \; {{A}^k = \begin{bmatrix}
			0_{d,d} & \mathcal{R}^k\\ (\mathcal{R}^k)^T & 0_{d,d}
	\end{bmatrix}}, \; f = \begin{pmatrix}
		\mathbf{r}\\\mathbf{g}
	\end{pmatrix}
\end{align*}
with $\mathcal{R}^k_{ij} = \epsilon_{ikj} $  being the Levi-Civita tensor. The graph space is $V=H(\text{curl}, \Omega)\times H(\text{curl}, \Omega)$. The boundary operator is 
\begin{align}
	\mathcal{D}=\sum_{k=1}^d n_k {A}^k = \begin{bmatrix} 0_{d,d} & \mathcal{T}\\ \mathcal{T}^T &0_{d,d}
	\end{bmatrix} ,\qquad
	 \text{with}\ \mathcal{T}\xi :=\mathbf{n} \times \boldsymbol{\xi},\\
		\langle D(\mathbf{H}, \mathbf{E}), (\mathbf{h}, \mathbf{e})\rangle_{V^{'}, V} = 
	(\mathbf{n}\times\mathbf{E}, \mathbf{e})_{L^2(\partial\Omega)}
	- (\mathbf{n}\times\mathbf{H}, \mathbf{h})_{L^2(\partial\Omega)}.
\end{align}
We impose homogeneous Dirichlet boundary conditions tangential to the electric field $\left(\mathbf{n}\times\mathbf{E}\right)_{|\partial\Omega}=0$ through
		\begin{equation}\label{eq:boundary_max}
		\mathcal{M}=\begin{bmatrix} 0_{d,d} & -\mathcal{T}\\ \mathcal{T}^T &0_{d,d}
	\end{bmatrix},\qquad
	\langle M(\mathbf{H}, \mathbf{E}), (\mathbf{h}, \mathbf{e})\rangle_{V^{'}, V} = 
	- (\mathbf{n}\times\mathbf{E}, \mathbf{e})_{L^2(\partial\Omega)}
	- (\mathbf{n}\times\mathbf{H}, \mathbf{h})_{L^2(\partial\Omega)}.
\end{equation}

\subsubsection{Compressible linear elasticity}
\label{subsubsec:comprLinElast_def}
We consider the parametric compressible linear elasticity system in $\R^d=\R^3$, where $\boldsymbol{\sigma} \in \R^{d\times d}$ is the stress tensor and $ \mathbf{u} \in \R^d$ is the displacement vector. The system can be written as
\begin{equation}\label{eq:elastic}
	\left(
		\begin{array}{c}
			\boldsymbol{\sigma} - \mu_1(\nabla\cdot \ub)\mathbb{I}_{3,3} -2\mu_2\frac{\left(\nabla \ub + (\nabla \ub)^t\right)}{2}\\
			-\frac{1}{2}\nabla\cdot\left(\boldsymbol{\sigma}+\boldsymbol{\sigma}^t\right) +\mu_3 \ub
	\end{array}
	\right)
	= 
	\left(
	\begin{array}{c}
		0 \\
		\mathbf{r}
	\end{array}
	\right),
	\quad \forall x\in\Omega,
\end{equation}
where $\mathbf{r} \in \R^3$, and $\mu_1,\ \mu_2>0$ are the Lam\'e constants. Rescaling the displacement $\mathbf{u}$ by $2\mu_2$, we obtain
\begin{equation}\label{eq:elastic_FS}
	\left(
	\begin{array}{c}
		{\boldsymbol{\sigma}} - \frac{\mu_1}{2\mu_2+3\mu_1}\text{tr}({\boldsymbol{\sigma}})\mathbb{I}_{3,3} -\frac{\left(\nabla \boldsymbol u + (\nabla \boldsymbol u)^T\right)}{2}\\
		-\frac{1}{2}\nabla\cdot\left({\boldsymbol{\sigma}}+{\boldsymbol{\sigma}}^T\right) +\frac{\mu_3}{2\mu_2} \boldsymbol u
	\end{array}
	\right)
	= 
	\left(
	\begin{array}{c}
		\mathbf{0} \\
		\mathbf{r}
	\end{array}
	\right),
	\quad \forall x\in\Omega.
\end{equation}
In this case, we consider the graph space
\begin{equation}
	V=H_{\boldsymbol{\sigma}}\times [H^1(\Omega)]^d,\quad H_{\boldsymbol{\sigma}}=\{\boldsymbol{\sigma}\in [L^2(\Omega)]^{d\times d}\mid \nabla\cdot(\boldsymbol{\sigma}+\boldsymbol{\sigma}^t)\in [L^2(\Omega)]^d\}.
\end{equation}
If we reorder the coefficients of ${\boldsymbol{\sigma}}$ into a vector, we can define $\mathbf{z}=\begin{pmatrix}
	\boldsymbol{\sigma}\\
	\mathbf{u}
\end{pmatrix}$ and have
\begin{equation}
	A^0 = \begin{bmatrix}
		\mathbb I_{d^2,d^2} -\frac{\mu_1}{2\mu_2+3\mu_1}\mathcal{Z}  & 0_{d^2,d}\\
		0_{d,d^2} & \frac{\mu_3}{2\mu_2} \mathbb I_{d,d}
	\end{bmatrix},\qquad A^k=\begin{bmatrix}
	0_{d^2,d^2} & \mathcal{E}^k\\ (\mathcal{E}^k)^T & 0_{d,d} 
\end{bmatrix}, \qquad {f} = \begin{bmatrix}
	\mathbf{0}_{d^2} \\ \mathbf{r}
\end{bmatrix},
\end{equation}
with $\mathcal{Z}_{[ij],[kl]}=\delta_{ij}\delta_{kl}$ and $\mathcal{E}^k_{[ij],l} = -\frac12 \left(\delta_{ik} \delta_{jl}+ \delta_{il} \delta_{jk}\right)$.
This leads to the definition of the boundary operator 
\begin{subequations}
	\begin{align}
		\mathcal{D}=\sum_{k=1}^d n_k {A}^k = \begin{bmatrix} 0_{d^2,d^2} & \mathcal{N}\\ \mathcal{N}^T &0_{d,d}
		\end{bmatrix} \text{ with }\mathcal{N}\boldsymbol{\xi} := -\frac12(\mathbf{n} \otimes\boldsymbol{\xi} + \boldsymbol{\xi} \otimes \mathbf{n}),\\
		\langle D({\boldsymbol{\sigma}}, \mathbf{u}), ({\boldsymbol{\tau}}, \mathbf{v})\rangle_{V^{'}, V} = 
	-\langle \tfrac{1}{2}(\boldsymbol{\sigma}+\boldsymbol{\sigma}^t)\cdot \n, \mathbf{v}\rangle_{-\frac{1}{2}, \frac{1}{2}}
	- \langle  \tfrac{1}{2}(\boldsymbol{\tau}+\boldsymbol{\tau}^t)\cdot \n, \mathbf{u}\rangle_{-\frac{1}{2}, \frac{1}{2}}.
	\end{align}
\end{subequations}
Mixed boundary conditions $\mathbf{u}_{|\Gamma_D} = 0$ and $(\boldsymbol{\sigma}\cdot\mathbf{n})_{|\Gamma_N}=0$ can be applied through the following boundary operator on the Dirichlet boundary $\Gamma_D$ and on the Neumann boundary $\Gamma_N$:
\begin{equation}\label{eq:boundary_bar}
	\begin{split}
	\langle M({\boldsymbol{\sigma}}, \mathbf{u}), ({\boldsymbol{\tau}}, \mathbf{v})\rangle_{V^{'}, V}= & 
	-\langle \tfrac{1}{2}(\boldsymbol{\sigma}+\boldsymbol{\sigma}^t)\cdot \n, \mathbf{v}\rangle_{-\frac{1}{2}, \frac{1}{2}, \Gamma_D}
	+ \langle  \tfrac{1}{2}(\boldsymbol{\tau}+\boldsymbol{\tau}^t)\cdot \n, \mathbf{u}\rangle_{-\frac{1}{2}, \frac{1}{2}, \Gamma_D}\\
	&+\langle \tfrac{1}{2}(\boldsymbol{\sigma}+\boldsymbol{\sigma}^t)\cdot \n, \mathbf{v}\rangle_{-\frac{1}{2}, \frac{1}{2}, \Gamma_N}
	- \langle  \tfrac{1}{2}(\boldsymbol{\tau}+\boldsymbol{\tau}^t)\cdot \n, \mathbf{u}\rangle_{-\frac{1}{2}, \frac{1}{2}, \Gamma_N},
	\end{split}
\end{equation}
the constructive procedure employed to define the boundary operator $M\in\mathcal{L}(V, V')$ is reported in the Appendix~\ref{appendix:boundary_operator}.

\subsubsection{Grad--div problem: advection--diffusion--reaction equations}
\label{subsubsec:adr_def}

Another example is the advection--diffusion--reaction equation
\begin{equation}\label{eq:ADR}
	-\nabla \cdot (\kappa \nabla u) +\boldsymbol{\beta} \cdot \nabla u +\mu u =\mathbf{r},
\end{equation}
with $\kappa\in [L^{\infty}(\Omega)]^{d\times d}$ and $\boldsymbol{\beta}\in [W^{1, \infty}(\Omega)]^d$, under the hypothesis that $\kappa\in [L^{\infty}(\Omega)]^{d\times d}$ and $\mu-\nabla\cdot\boldsymbol{\beta}\in L^{\infty}(\Omega)$ are uniformly bounded from below to satisfy the positivity property~\eqref{eq:positivity_classical}. 
Let us write the equation in the mixed form with $\boldsymbol{\sigma}=-\kappa \nabla u$ and $\mathbf{z}=\begin{pmatrix}
	\boldsymbol{\sigma}\\ u
\end{pmatrix}$.
Then, \eqref{eq:ADR} can be rewritten as \eqref{eq:FS} with 
\begin{equation}
	 A^0 = \begin{bmatrix}
		\kappa^{-1} & 0_{d,1}\\
		0_{1,d} & \mu
	\end{bmatrix}, \qquad A^k = \begin{bmatrix}
		0_{d,d} & \mathbf{e}_k\\ (\mathbf{e}_k)^T & \beta_k
	\end{bmatrix}, \qquad \mathbf{f} = \begin{pmatrix}
		0\\\mathbf{r}
	\end{pmatrix}.
\end{equation}
Here, $0_{m,\ell} \in \R^{m\times \ell}$ is a matrix of zeros and $\mathbf{e}_k$ is the unitary vector with the $k$-th entry equal to 1. The graph space is $V=H(\text{div}, \Omega)\times H^2(\Omega)$. The boundary operator $D$ becomes
\begin{equation}
	D = \sum_{k=1}^d n_k A^k = \begin{bmatrix}
		0_{d,d} &\mathbf{n} \\ \mathbf{n}^t& \beta \cdot \mathbf{n}
	\end{bmatrix},\qquad\langle D({\boldsymbol{\sigma}}, u), ({\boldsymbol{\tau}}, v)\rangle_{V^{'}, V} = \langle \boldsymbol{\sigma}\cdot \n, v\rangle_{-\frac{1}{2}, \frac{1}{2}}
- \langle  \boldsymbol{\tau}\cdot \n, u\rangle_{-\frac{1}{2}, \frac{1}{2}}.
\end{equation}
Homogeneous Dirichlet boundary conditions $u_{|\partial\Omega}=0$ can be imposed with 
\begin{equation}
	\mathcal{M}=\begin{bmatrix}
		0_{d,d} & -\mathbf{n}\\
		\mathbf{n}^t & 0
	\end{bmatrix},
\end{equation}
while Robin/Neumann boundary conditions of the type $\boldsymbol{\sigma}\cdot \mathbf{n} = \gamma u$ are imposed with 
\begin{equation}
	\mathcal{M}=\begin{bmatrix}
		0_{d,d} & \mathbf{n}\\
		-\mathbf{n}^t & 2\gamma + \boldsymbol{\beta} \cdot \mathbf{n}
	\end{bmatrix}.
\end{equation}
For our test case in Section~\ref{sec:vv}, we will consider as advection field $\boldsymbol{\beta}:\Omega\rightarrow\mathbb{R}^d$ an incompressible velocity field from the solution of the 2d incompressible Navier-Stokes equations as described later. Similarly to the linear compressible elasticity mixed boundary conditions in Section~\ref{subsubsec:comprLinElast_def}, we want to impose $u_{|\Gamma_D} = \mathbf{g}\in [L^{\frac{1}{2}}(\Gamma_D)]^{d}$ and $\left(\boldsymbol{\sigma}\cdot\mathbf{n}\right)_{|\Gamma_N} = 0$. This is possible with
\begin{equation}\label{eq:boundary_adr}
	\langle M(\boldsymbol{\sigma}, u), (\boldsymbol{\tau}, v)\rangle_{V^{'}, V} = 
	\langle \boldsymbol{\sigma}\cdot \n, v\rangle_{-\frac{1}{2}, \frac{1}{2}, \Gamma_D}
	+ \langle  \boldsymbol{\tau}\cdot \n, u\rangle_{-\frac{1}{2}, \frac{1}{2}, \Gamma_D}
	-\langle \boldsymbol{\sigma}\cdot \n, v\rangle_{-\frac{1}{2}, \frac{1}{2}, \Gamma_N}
	- \langle  \boldsymbol{\tau}\cdot \n, u\rangle_{-\frac{1}{2}, \frac{1}{2}, \Gamma_N},
\end{equation}
the proof is similar to the one reported in Appendix~\ref{appendix:boundary_operator}.

\section{Discontinuous Galerkin discretization}
\label{sec:fom}
In the literature, a few discretization approaches for FS are presented, e.g. finite volume method \cite{sonar1998dual} or discontinuous Galerkin (DG) method~\cite{ern2006discontinuous,di2011mathematical,bui2013unified}. More recently, a hp-adaptive hybridizable DG formulation was introduced in~\cite{chen2023unified}.
In this work, we perform a DG discretization following the notation reported in \cite{ern2006discontinuous,di2011mathematical}. 
Consider a shape-regular tessellation $\Th$ of the domain $\Omega$ and take a piecewise polynomial space $V_h$ over $\Th$, defined by $V_h= \lbrace z \in V: z|_T \in \mathbb{P}^k_d(T), \forall T \in \Th\rbrace$, where $k$ is the polynomial degree. 
We assume that there is a partition $P_{\Omega}=\{\Omega_i\}_{1\leq i\leq N_\Omega}$ of $\Omega$ into disjoint polyhedra such that the exact solution $z$ belongs to $V^*=V\cap[H^1(P_\Omega)]^m$. 
We define the discrete bilinear form $\forall y_h\in V_h,\ z\in V^*$
\begin{subequations}\label{eq:DGbilinears}
	\begin{align}\label{eq:DGbilinear}
		a^{cf}_h(z,y_h) =& \sum_{T \in \mathcal{T}_h}  (Az,y_h)_{L^2(T)} +{\tfrac{1}2 \sum_{F \in \mathcal{F}^b_h}\left((\mathcal{M}-\mathcal{D})z,y_h\right)_{L^2(F)}} -{\sum_{F \in \mathcal{F}^i_h} \left(\mathcal{D}_F[\![z]\!],\lbrace\!\lbrace y_h \rbrace \!\rbrace  \right)_{L^2(F)}}\\
		=&\sum_{T \in \mathcal{T}_h}  (z,\tilde{A}y_h)_{L^2(T)} +{\tfrac{1}2 \sum_{F \in \mathcal{F}^b_h}\left((\mathcal{M}+\mathcal{D})z,y_h\right)_{L^2(F)}} +{\sum_{F \in \mathcal{F}^i_h} \left(\mathcal{D}_F\lbrace\!\lbrace z \rbrace \!\rbrace,[\![y_h]\!]  \right)_{L^2(F)}},\label{ref:dg_adj_sys}
	\end{align}
\end{subequations}
	where the first two terms are the piece-wise discontinuous discretization of the bilinear form \eqref{eq:bilinear} and the last term penalizes the jump across neighboring cells and stabilizes the method. Here, $\mathcal{F}^b_h$ is the collection of the faces  of the triangulation $\Th$ belonging to the boundary of $\Omega_h$, while $\mathcal{F}^i_h$ is the collection of internal faces. The jump and the average of a function on a face $F$ shared by two elements $T_1$ and $T_2$ are defined as $[\![u]\!]=u|_{T_1}-u|_{T_2}$ and $\{\!\{ u\}\!\} =\frac12(u|_{T_1}+u|_{T_2})$, respectively. The boundary operator $\mathcal{D}:\partial\Omega\rightarrow\mathbb{R}^{m\times m}$ can be extended also on the internal faces $F \in \mathcal{F}^i_h$ as $\mathcal{D}_F = \sum_{k=1}^d n_k^F A^k$, where $n^F$ is a normal to the face $F$ and it is well-defined.

In order to obtain quasi-optimal error estimates, extra stabilization terms are needed. We additionally impose that $A^i\in [C^{0,\tfrac{1}{2}}(\overline{\Omega}_j)]^{m\times m},\ \forall T\in\mathcal{T}_h,\ i=1,\dots,d,\ \forall \Omega_j\in P_{\Omega}$. A possibility is given by the following stabilization term
\begin{equation}\label{eq:stabilization}
	s_h(z,y_h)=\sum_{F\in \mathcal{F}^{b}_h} (S^b_Fz,y_h)_{L^2(F)} + \sum_{F\in \mathcal{F}^i_h} (S^i_h [\![z]\!],[\![y_h]\!])_{L^2(F)},
\end{equation}
where the operators $S_h^i$ and $S_F^b$ have to satisfy the following constraints for some $\alpha_j>0$ for $j=1,\dots,5$:
\begin{subequations}
\begin{align}
	&S^b_F z = 0\quad\forall F \in \mathcal{F}^b_h, \qquad S^i_F [\![z]\!]=0\quad\forall F \in \mathcal{F}^i_h,\qquad \text{with $z$ the exact solution,}\\
	&S^b_F\text{ and }S^i_F\text{ are symmetric and nonnegative,}\\
	&S^b_F\leq \alpha_1 \mathbb I_{m,m}, \qquad \alpha_2 |D_F| \leq S^i_F \leq \alpha_3 \mathbb I_{m,m},\\
	& |((M-D)y,z)_{L^2(F)}| \leq \alpha_4 ((S^b_F+M)y,y)^{1/2}_{L^2(F)}\lVert z\rVert_{L^2(F)},\\
	& |((M+D)y,z)_{L^2(F)}| \leq \alpha_5 ((S^b_F+M)z,z)^{1/2}_{L^2(F)}\lVert y\rVert_{L^2(F)}.
\end{align}
\end{subequations}
Specific definitions of these operators for our test cases are presented in~\cite{ern2006discontinuous, di2011mathematical}, properly declined for our mixed boundary conditions in the compressible linear elasticity and advection--diffusion--reaction test cases, see Sections~\ref{subsubsec:comprLinElast_def} and~\ref{subsubsec:adr_def}, respectively.
Finally, we can define the bilinear form and the right-hand side
\begin{equation}\label{eq:DGbilinearstab}
	a_h(z,y_h) = a^{cf}_h(z,y_h) +s_h(z,y_h),\quad l_h(y_h) = \sum_{T \in \mathcal{T}_h}  (f,y_h)_{L^2(T)} + {\tfrac{1}2 \sum_{F \in \mathcal{F}^b_h}\left((M-D)g,y_h\right)_{L^2(F)}}
\end{equation}
that lead to the definition of the discrete problem.
\begin{definition}[DG Friedrichs' System]
	Given $f\in L$ and $g\in V_h$, the DG approximation of the FS constitute in finding a $z_h\in V_h$ such that
	\begin{equation}\label{eq:DGFS}
		a_h(z_h,y_h)=l_h(y_h),\qquad \forall y_h \in V_h.
	\end{equation} 
\end{definition}
To prove the accuracy of the discrete problem, it is necessary to have the following conditions:
\begin{itemize}
	\item Consistency, i.e., $a_h(z,y_h)=a(z,y_h)$ for $z\in V^*$;
	\item $L^2$-coercivity, i.e., $a_h(y_h,y_h) \geq \mu_0 \lVert y_h\rVert_{L}^2 + \frac12 \lvert y_h\rvert^2_M$, with $\lvert y_h\rvert^2_M= \int_{\partial \Omega} y^tMy$;
	\item Inf-sup stability
	\begin{equation}\label{eq:inf_sup_tight}
		||| z_h|||\lesssim \sup_{y_h \neq 0} \frac{a_h(z_h,y_h)}{|||y_h|||}
	\end{equation}
with  $|||y|||^2=||y||_{L^2}^2 + |y|_M^2 +|y|_S^2 +\!\! \sum_{T\in \mathcal{T}_h} h_T \left \lVert A^k \partial_k y\right \rVert_{L^2(T)}^2$ and $|y|_S^2=s_h(y,y)$;
	\item Boundedness $a_h(w,y_h)\lesssim |||w|||_* |||y_h|||$ with
	\begin{equation}
		\label{eq:inf-sup_norm}
		\!\!\!|||y|||_*^2=|||y|||^2+ \sum_{T\in \mathcal{T}_h} \left(h^{-1}_T ||y||_{L^2(T)}^2 + ||y||_{L^2(\partial T)}^2\right).
	\end{equation}
\end{itemize} 

\begin{theorem}[Error estimate from~\cite{ern2006discontinuous,di2011mathematical}]
	\label{theo:error_estimate}
	Let $z\in V^*$ be the solution of the weak problem~\eqref{eq:weak_continuous_FS} and $z\in V_h$ be the solution of the discrete DG problem~\eqref{eq:DGFS}. Then, the consistency and inf-sup stability of the discrete system~\eqref{eq:DGFS} imply
	\begin{equation}
		|||z-z_h|||\lesssim \inf_{y_h\in V_h} |||z-y_h|||_*,
	\end{equation}
	in particular, if $z\in [H^{k+1}(\Omega)]^m$ the following convergence rate holds
	\begin{equation}
		|||z-z_h|||\lesssim h^{k+\tfrac{1}{2}}\lVert z\rVert_{[H^{k+1}(\Omega)]^m}.
	\end{equation}
\end{theorem}

\section{Projection-based model order reduction}
\label{sec:mor}
The computation of discrete solutions of parametrized PDEs can require a not negligible computational time. In particular, in multi-query context, when many evaluations for different parameters are required, the computations may become unbearable. In this section, we introduce a reduced order model (ROM) for the FS in case of parameter dependent problems \cite{hesthaven2016certified,rozza2022advanced}, in order to drastically reduce the computational costs.
To do so, we exploit two aspects of the above presented FS: the linearity of the problems and the affine dependence of the operators on the physical parameters.

As we have seen in Section~\ref{sec:examples}, all the problems are depending on some parameters $\boldsymbol{\rho} \in\mathcal{P} \subset \R^{\Npar}$ and the dependence is affine. 
This means that it is possible to find $\Naff$ terms independent on the parameters for each form, such that they can be affinely combined with some parameter dependent functions to obtain the original operator, i.e.,
\begin{equation}
	a_h(z,y_h;\boldsymbol{\rho} ) = \sum_{\ell=1}^{\Naff} \theta_\ell^a(\boldsymbol{\rho} ) a_{\ell,h}(z,y_h), \qquad l_h(y_h) = \sum_{\ell=1}^{\Naff} \theta_\ell^f(\boldsymbol{\rho} ) l_{h, l}(y_h).
\end{equation}

Then, we select a reduced space $V_\NRB\subset V_h$ provided by a compression algorithm, e.g. SVD/POD/PCA \cite{jolliffe2002principal,kunisch2001galerkin,torlo2021model} or Greedy algorithm \cite{prud2002mathematical,prud2002reliable,hesthaven2016certified,crisovan2019model}. We suppose that the reduced dimension $\NRB$ is much smaller than the dimension $\NDG$ of the full order model space $V_h$. 
We look as ansatz for a reduced solution $\zRB\in V_\NRB$ a linear combination of the bases $\lbrace \psiRB_j\rbrace_{j=1}^{\NRB} $ of $ V_\NRB$, i.e.,
\begin{equation}
	\zRB = \sum_{j=1}^{\NRB} \zRB^j \psiRB_j,
\end{equation}
then, performing a standard Galerkin projection, we obtain the following $\RB$ problem.
\begin{definition}[Reduced Basis Problem]
	Find $\zRB\in V_\NRB$, given by the coefficients $\zRB^j$, such that 
	\begin{equation}
		\sum_{i=1}^{\NRB} \zRB^j(\boldsymbol{\rho} ) \sum_{\ell=1}^{\Naff} \theta_\ell^a(\boldsymbol{\rho} ) a_{\ell,h}(\psiRB_j,\psiRB_i) = \sum_{\ell=1}^{\Naff}\theta_\ell^f(\boldsymbol{\rho} )(f_\ell,\psiRB_i), \qquad \text{for all }i=1,\dots, \NRB. 
	\end{equation}
\end{definition}
The obtained problem scales depend on the dimension $\NRB$ and $\Naff$ in its assembly and only on $\NRB$ in its solution, and it is completely independent on $\NDG$.
To obtain computational advantages for the parametric problem, we split the tasks into an expensive \textit{offline phase} and a cheap \textit{online phase}. 
In the \textit{offline phase}, we find the reduced space $V_\NRB$ and we assemble the reduced matrices and right hand sides 
\begin{equation}
	A_\ell:=\lbrace a_{\ell,h}(\psiRB_j,\psiRB_i) \rbrace_{i,j}, \qquad b_\ell:=\lbrace (f_\ell,\psiRB_i)\rbrace_{i}.
\end{equation}
In the \textit{online phase}, we can simply evaluate the coefficients $\theta^a_\ell(\boldsymbol{\rho} )$ and $\theta^f_\ell(\boldsymbol{\rho} )$ and obtain the reduced linear system
\begin{equation}
	A(\boldsymbol{\rho} )z_{\RB} = b(\boldsymbol{\rho} ), \qquad \text{with } A(\boldsymbol{\rho} ):=\sum_\ell\theta^a_\ell(\boldsymbol{\rho} ) A_\ell \text{ and }b(\boldsymbol{\rho} ):= \sum_\ell\theta^f_\ell(\boldsymbol{\rho} ) b_\ell.
\end{equation}
This gives a great speed up in computational times.

\subsection{Reduced basis a posteriori error estimate}
\label{sec:certified}
We derive two error estimators for the energy norm and the $L^2$ norm of the reduced basis error $e_h=z_h-z_{RB}\in V_h$ following the procedure in \cite{hesthaven2016certified}.
\begin{subequations}
	\label{eq:coercivity_nrg_R}
Exploiting the equality in \eqref{eq:DGbilinears}, we obtain the following lower bound 
\begin{align}
	a_h(e_h,e_h) =& a^{cf}_h(e_h,e_h) +s_h(e_h,e_h)\\
    =&\sum_{T \in \mathcal{T}_h}  (Ae_h,e_h)_{L^2(T)} +\tfrac{1}2 \sum_{F \in \mathcal{F}^b}\left((\mathcal{M}-\mathcal{D})e_h,e_h\right)_{L^2(F)} - \sum_{F \in \mathcal{F}^i_h} \left(\mathcal{D}_F[\![e_h]\!],\lbrace\!\lbrace e_h \rbrace \!\rbrace  \right)_{L^2(F)}+\\
    &\sum_{F\in \mathcal{F}^{b}_h} (S^b_Fe_h,e_h)_{L^2(F)} + \sum_{F\in \mathcal{F}^i_h} (S^i_h [\![e_h]\!],[\![e_h]\!])_{L^2(F)}\\ \label{eq:inequality_bound_a}
    =&\sum_{T \in \mathcal{T}_h}  ((A^0 - \tfrac{1}{2} \mathcal{X} )e_h,e_h)_{L^2(T)}+\tfrac{1}2 \sum_{F \in \mathcal{F}^b}\left(\mathcal{M}e_h,e_h\right)_{L^2(F)} +\sum_{F\in \mathcal{F}^{b}_h} (S^b_Fe_h,e_h)_{L^2(F)} + \sum_{F\in \mathcal{F}^i_h} (S^i_h [\![e_h]\!],[\![e_h]\!])_{L^2(F)}\\
    \geq& \mu_0 \lVert e_h\rVert^2_L +\tfrac{1}{2}|e_h|^2_M +|e_h|^2_S,\label{eq:coercivity_bilinear_form}
\end{align}
\end{subequations}
where we have defined $|\cdot|_M^2 = \tfrac{1}2 \sum_{F \in \mathcal{F}^b}\left(\mathcal{M}\cdot,\cdot\right)_{L^2(F)}$ and $|\cdot|_S^2 =s_h(\cdot,\cdot)$. 
We define the $R$-norm
\begin{equation}
	||y_h||_{R}^2 = \mu_0\lVert y_h\rVert^2_L +\tfrac{1}{2}|y_h|^2_M +|y_h|^2_S,\quad \forall y_h\in V_h,
\end{equation}
that may depend on $\rho$ only through $\mu_0$ and is generated by the scalar product
\begin{equation}
    \langle u_h, v_h \rangle_{R} = \mu_0\sum_{T \in \mathcal{T}_h}  (u_h,v_h)_{L^2(T)} + \tfrac{1}{2}\sum_{F \in \mathcal{F}^b}\left(\mathcal{M}^{sym} u_h,v_h\right)_{L^2(F)} +  \sum_{F^b} (S^b u_h,v_h) +\sum_{F^i} (S^i u_h,v_h).
\end{equation}
The boundary operators we will employ in our benchmarks are all skew-symmetric so $\mathcal{M}^{sym}=\frac{\mathcal{M}+\mathcal{M}^t}{2}$ is the null matrix and $|e_h|_M = 0$. Now, we can proceed to provide an \textit{a posteriori} error estimate for the $R$-norm and energy norm. Hence, let us define $r_{RB}(y_h)= l_h(y_h;\rho)- a_h(z_{RB}, y_h;\rho) $ and its $R$ and $L$-Riesz representations as $\hat{r}_{R}$ and $\hat{r}_{L}$ such that
\begin{equation}
    \label{eq:riesz-repr}
    r_{RB}(u_h) = \langle \hat{r}_{R}, u_h \rangle_{R},\quad X_{R} \hat{\mathbf{r}}_{R} = \mathbf{L}_h- A_h \mathbf{z}_{RB},\quad
    r_{RB}(u_h) = \langle \hat{r}_{L}, u_h \rangle_{L},\quad X_{L} \hat{\mathbf{r}}_{L} = \mathbf{L}_h- A_h \mathbf{z}_{RB},
\end{equation} 
where $X_R$ and $X_{L}$ are the $R$-norm and $L$-norm mass matrices, and $L_h$, $A_h$ and $\mathbf{z}_{RB}$ are the representations of $l_h(\cdot;\rho)$, $a_h(\cdot, \cdot;\rho)$ and $z_{RB}$ in the DG basis of $V_h$. The $\hat{r}_{L}$ representation can be computed cheaply when the parametric model is affinely decomposable with respect to the parameters, while $\hat{r}_{R}$ requires the inversion of a possibly parametric dependent matrix $X_{R}$.

Now, consider the energy norm of the error $\lVert e_h\rVert^2_{nrg} = a_h(e_h,e_h)$ and the coercivity constant $\lVert e_h\rVert^2_{nrg} \geq \mu_0\lVert e_h\rVert_L^2$ derived in \eqref{eq:coercivity_bilinear_form}, we have the following \textit{a posteriori} error estimates
\begin{equation}
    \frac{\lVert e_h\rVert_{nrg}}{\lVert z_h\rVert_{nrg}}\leq \frac{||\hat{r}_R||_R}{\lVert z_h\rVert_{nrg}} ,
    \qquad \frac{\lVert e_h\rVert_{R}}{\lVert z_h\rVert_{R}}\leq \frac{||\hat{r}_R||_R}{\lVert z_h\rVert_{R}},\qquad
    \frac{\lVert e_h\rVert_{nrg}}{\lVert z_h\rVert_{nrg}}\leq \frac{||\hat{r}_L||_L}{\sqrt{\mu_0}\lVert z_h\rVert_{nrg}} ,
    \qquad \frac{\lVert e_h\rVert_{L}}{\lVert z_h\rVert_{L}}\leq \frac{||\hat{r}_L||_L}{\mu_0\lVert z_h\rVert_{L}},
\end{equation}
namely, the relative energy error with the corresponding \textit{a posteriori} $R$-norm energy estimate, the relative $R$-norm error with the corresponding \textit{a posteriori} $R$-norm estimate, the relative energy error with the corresponding \textit{a posteriori} $L$-norm energy estimate, and the relative $L$-norm error with the corresponding \textit{a posteriori} $L$-norm estimate.

\subsection{Optimally stable error estimates for the ultraweak Petrov-Galerkin formulation}
\label{subsec:ultraweak}
In this section, we show that Friedrichs' systems are a desirable unifying formulation to consider when performing model order reduction also due to the possibility to achieve an optimally stable formulation.
This can further simplify the error estimator analysis reaching the equality between the error and the residual norm.
This is not the first case in which optimally stable formulations are introduced also at the reduced level, see~\cite{beurer2022ultraweak, hain2022ultra, henning2022ultraweak}. 
In the following, we describe how to achieve this ultraweak formulation and we delineate the path one should follow to use such formulation.
Nevertheless, we will not use this formulation in our numerical tests and we leave the implementation to future studies.

We introduce the following Discontinuous Petrov-Galerkin (DPG) formulation from~\cite{bui2013unified}. 
To do so, we first define $V(\mathcal{T}_h)$ the broken graph space with norm $\lVert\bullet\rVert^2_{V(\mathcal{T}_h)}=\lVert\bullet\rVert^2_L +\sum_{T\in\mathcal{T}_h} \lVert A\bullet\rVert^2_{L^2(T)}$ and  $\Tilde{V}=V/Q(\Omega)$ the quotient of the graph space $V$ with
\begin{equation}
	\begin{split}
	Q(\Omega)=&\left\{z\in V\,\bigg| \, {\sum_{T \in \mathcal{T}_h} \langle D z,y \rangle_{\Tilde{V}(T),V(T)}}+{\tfrac{1}2 \langle (M-D)z,y\rangle_{\Tilde{V}(\Omega),V(\Omega)}}=0,\quad\forall y\in V(\mathcal{T}_h)\right\}\\
	=&\left\{z\in V\,\bigg| \,a(z, y)=0,\quad\forall y\in V(\mathcal{T}_h)\right\}.
	\end{split}
\end{equation}
The DPG formulation reads: find $(z, q)\in L\times \Tilde{V}$ such that, for all $y\in V(\mathcal{T}_h)$,
\begin{equation}
    \label{eq:PGDG}
    \sum_{T \in \mathcal{T}_h}  (z,\tilde{A}y)_{L^2(T)} +{\sum_{T \in \mathcal{T}_h} \langle D q,y \rangle_{\Tilde{V}(T),V(T)}}+{\tfrac{1}2 \langle(M-D)q,y\rangle_{\Tilde{V},V}} = {\sum_{T \in \mathcal{T}_h} \left(f,v \right)_{\Tilde{V}(T),V(T)}} + {\tfrac{1}2 \langle(M-D)g,y\rangle_{\Tilde{V},V}}.
\end{equation}
%where $\Tilde{V}=V/Q(\Omega)$ is the quotient of the graph space $V$ with
%\begin{equation}
%    Q=\left\{z\in V\bigg| {\sum_{T \in \mathcal{T}_h} \langle D z,y \rangle_{\Tilde{V}(T),V(T)}}+{\tfrac{1}2 \langle (M-D)z,y\rangle_{\Tilde{V}(\Omega),V(\Omega)}}=0,\quad\forall y\in V(\mathcal{T}_h)\right\}=\left\{z\in V\bigg| a(z, y)=0,\quad\forall y\in V(\mathcal{T}_h)\right\}
%\end{equation}
%and $V(\mathcal{T}_h)$ is the broken graph space with norm $\lVert\bullet\rVert^2_{V(\mathcal{T}_h)}=\lVert\bullet\rVert^2_L +\sum_{T\in\mathcal{T}_h} \lVert A\bullet\rVert^2_{L^2(T)}$. 
The introduction of the hybrid face variables $q\in\Tilde{V}$ is necessary since $z\in L$ does not satify~\eqref{eq:def_D}. In practice, assuming that the traces of $y\in V(\mathcal{T}_h)$ are well-defined and belong to a space $X(\mathcal{F}_{i,b})$, we can formulate~\eqref{eq:PGDG} as follows: find $(z, q)\in L\times X(\mathcal{F}_{i,b})$ such that, for all $y\in V(\mathcal{T}_h)$,
\begin{equation}
    \label{eq:PGDGb}
    \sum_{T \in \mathcal{T}_h}  (z,\tilde{A}y)_{L^2(T)} +{\sum_{F \in \mathcal{F}_{i}} \left(\mathcal{D} q,[\![ y ]\!] \right)_{X(F)}}+{\tfrac{1}2 \sum_{F \in \mathcal{F}_{b}} \left((\mathcal{M}-\mathcal{D})q,y\right)_{X(F)}} = {\sum_{T \in \mathcal{T}_h} \left(f,v \right)_{\Tilde{V}(T),V(T)}} + {\tfrac{1}2 \sum_{F \in \mathcal{F}_{b}}\left((\mathcal{M}-\mathcal{D})g,y\right)_{X(F)}},
\end{equation}
where $X(\mathcal{F}_{i,b})$ is, for example, 
$[H^{-\frac{1}{2}}(\mathcal{F}_{i,b})]^d\times [H^{\frac{1}2}(\mathcal{F}_{i,b})]^d$ for compressible linear elasticity, 
$H^{-\frac{1}{2}}(\mathcal{F}_{i,b})\times H^{\frac{1}{2}}(\mathcal{F}_{i,b})$ for the scalar advection--diffusion--reaction 
and 
$L^2_\mathcal{T}(\mathcal{F}_{i,b})\times L^2_\mathcal{T}(\mathcal{F}_{i,b})$ for the Maxwell equations in stationary regime, with
$L^2_\mathcal{T}(\mathcal{F}_{i,b})$ being the space of fields in $H(\text{curl}, \mathcal{T}_h)$ whose tangential component belongs to $[L^2(\mathcal{F}_{i,b})]^3$.

The problem~\eqref{eq:PGDG} above is well-posed and consistent~\cite[Lemma 2.4]{bui2013unified} with the previous formulation in~\eqref{eq:weak_continuous_FS}. We consider the optimal norms, 
\begin{equation}
    \lVert(z,q)\rVert^2_{\mathcal{U}} = \sum_{T\in\mathcal{T}_h} \lVert z\rVert^2_{L(T)} + \lVert q\rVert^2_{\Tilde{V}},\qquad \lVert y\rVert^2_{\mathcal{Y}}=\sum_{T\in\mathcal{T}_h}\lVert \Tilde{A}y\rVert^2_{L(T)}+\lVert[\![y]\!]\rVert^2_{\partial\Omega_h},\quad\text{with}\quad\lVert[\![y]\!]\rVert_{\partial\Omega_h} = \sup_{q\in\Tilde{V}}\frac{a(q, y)}{\lVert q\rVert_{\Tilde{V}}},
\end{equation}
or formally, considering~\eqref{eq:PGDGb}
\begin{equation}
    \lVert(z,q)\rVert^2_{\mathcal{U}} = \sum_{T\in\mathcal{T}_h} \lVert z\rVert^2_{L(T)} + \sum_{F\in\mathcal{F}_{i,b}}\lVert q\rVert^2_{X(F)},\qquad \lVert y\rVert^2_{\mathcal{Y}}=\sum_{T\in\mathcal{T}_h}\lVert \Tilde{A}y\rVert^2_{L(T)}+\sum_{F\in\mathcal{F}_{i}}\lVert\mathcal{D}_F [\![y]\!]\rVert^2_{X(F)}+\sum_{F\in\mathcal{F}_{b}}\lVert(\mathcal{M}^t-\mathcal{D}) y\rVert^2_{X(F)}.
\end{equation}
With these optimal norms for the trial and test spaces we have the following result~\cite[Theorem 2.6]{bui2013constructively}.
\begin{theorem}[Optimally stable formulation]
    The bilinear form $b:(L\times\Tilde{V}, \lVert\bullet\rVert_{\mathcal{U}})\rightarrow (V(\mathcal{T}_h), \lVert\bullet\rVert_{\mathcal{Y}})$ defined as
    \begin{equation}
        b(u, y)=\sum_{T \in \mathcal{T}_h}  (z,\tilde{A}y)_{L^2(T)} +{\sum_{T \in \mathcal{T}_h} \langle D q,y \rangle_{\Tilde{V}(T),V(T)}}+{\tfrac{1}2 \langle(M-D)q,y\rangle_{\Tilde{V},V}}
    \end{equation}
    with $u=(z, q)$, is an isometry between $L\times\Tilde{V}$ and $V'(\mathcal{T}_h)$: we have that $\gamma = \beta = \beta^* = 1,$ where
    \begin{equation}
        \gamma := \sup_{u\in \mathcal{U}}\sup_{y\in \mathcal{Y}}\frac{b(u, y)}{\lVert u\rVert_\mathcal{U}\lVert y\rVert_{\mathcal{Y}}},\quad
        \beta := \inf_{u\in \mathcal{U}}\sup_{y\in \mathcal{Y}}\frac{b(u, y)}{\lVert u\rVert_\mathcal{U}\lVert y\rVert_{\mathcal{Y}}},\quad
        \beta^* := \inf_{y\in \mathcal{Y}}\sup_{u\in \mathcal{U}}\frac{b(u, y)}{\lVert u\rVert_\mathcal{U}\lVert y\rVert_{\mathcal{Y}}}
    \end{equation}
    with $\mathcal{U}=(L\times\Tilde{V}, \lVert\bullet\rVert_{\mathcal{U}})$.
\end{theorem}

This property is inherited at the discrete level as long as fixed $U^{N_h}_{h}\subset\mathcal{U}$ a discretization of the trial space with $\text{dim}\ Z_h=N_h$, the discrete test space $Y^{N_h}_h\subset V(\mathcal{T}_h)$ is the set of supremizers
\begin{equation}
    Y^{N_h}_h = \text{span}\left\{y_{u_h}\in V(\mathcal{T}_h)\,\bigg|\, y=\argmax_{y\in V({\mathcal{T}_h})} \frac{b(u_h, y)}{\lVert y\rVert_{\mathcal{Y}}},\quad u_h\in Z^{N_h}_h\right\}
\end{equation}
and $\text{dim}\ U_h=\text{dim}\ Y_h$, see~\cite[Lemma 2.8]{bui2013constructively}. In particular, for every $u_h=(z_h, q_h)\in U^{N_h}_h$, we have the optimal \textit{a posteriori} error estimate
\begin{equation}
    \lVert u-u_h\rVert_{\mathcal{U}} =\sup_{y\in V({\mathcal{T}_h})} \frac{b(u-u_h, y)}{\lVert y\rVert_{\mathcal{Y}}}=\lVert r_h(u_h) \rVert_{(V(\mathcal{T}_h))'}, \quad \langle r_h(u_h), v\rangle_{(V(\mathcal{T}_h))', V(\mathcal{T}_h)} = \langle f, v\rangle_{(V(\mathcal{T}_h))', V(\mathcal{T}_h)} - b(u_h, v).
\end{equation}
The same reasoning can be iterated another time to perform model order reduction with the choice $V_n=\lbrace \psiRB_j\rbrace_{j=1}^{\NRB}\subset U^{N_h}_{h}\subset\mathcal{U}$, and
\begin{equation}
    Y^{RB} = \text{span}\left\{y_{u_{RB}}\in V(\mathcal{T}_h)\,\bigg|\,y=\argmax_{y\in V({\mathcal{T}_h})} \frac{b(u_{RB}, y)}{\lVert y\rVert_{\mathcal{Y}}},\quad u_{RB}\in V_n\right\},
\end{equation}
such that for $u_{RB}\in V_n$,
\begin{equation}
    \lVert u-u_{RB}\rVert_{\mathcal{U}} = \lVert r_{RB}(u_{RB}) \rVert_{(V(\mathcal{T}_h))'}, \quad \langle r_{RB}(u_{RB}), v\rangle_{(V(\mathcal{T}_h))', V(\mathcal{T}_h)} = \langle f, v\rangle_{(V(\mathcal{T}_h))', V(\mathcal{T}_h)} - b(u_{RB}, v).
\end{equation}

The main difficulty is the evaluation of the trial spaces $Y^{N_h}_h$ and $Y^{RB}$ since the bilinear form $b$ may depend on the parameters $\boldsymbol{\rho}$. If the parameters affect only the source terms, the boundary conditions or the initial conditions for time-dependent FS, this problem is avoided. The evaluation of $Y^{N_h}_h$ can be performed locally for each element $T\in\mathcal{T}_h$, differently from $Y^{RB}$. An example of the evaluation of the basis of $Y^{N_h}_h$ is presented in~\cite[Equations 24, 25]{bui2013constructively} for linear scalar hyperbolic equations that can be interpreted as FS.

\section{Domain decomposable Discontinuous Galerkin ROMs}
\label{sec:results}
Extreme-scale parametric models are unfeasible to reduce with standard approaches due to the high computational costs of the offline stage. Parametric multi-physics simulations, such as fluid-structure interaction problems, are reduced inefficiently with a global reduced basis, depending on the complexity of the interactions between the physical models considered and the parametric dependency. In some cases, only a part of a decomposable system is reducible with a ROM, thus a possible solution is to implement a ROM-FOM coupling through an interface. In presence of moving shocks~\cite{baiges2013domain} affected by the parametrization, one may want to isolate these difficult features to approximate and apply different dimension reduction methodologies depending on the subdomain. These are the main reasons to develop domain decomposable or partitioned ROMs (DD-ROMs).

Some approaches from the literature are the reduced basis element methods~\cite{maday2002reduced, maday2004reduced}, the static condensation reduced basis element method~\cite{huynh2013static, eftang2012adaptive}, non-intrusive methods based on local regressions in each subdomain~\cite{xiao2019domain, xiao2019domainb}, overlapping Schwarz methods~\cite{del2023boundary,iollo2023one}, optimization-based MOR approaches \cite{prusak2022optimisation,prusak2023optimisation} and hyper-reduced ROMs~\cite{lucia2002domain}. In this last case, local approximations are useful because the local reduced dimensions are smaller and therefore more accurate local regressions can be designed to perform non-intrusive surrogate modelling. Little has been developed for the DG method, even though its formulation imposes naturally flux and solution interface penalties at the internal boundaries of the subdomains, in perspective of performing model order reduction. In our case, the linear systems associated to the parametric models are algebraically partitioned in disjoint subdomains coupled with the standard penalties from the weak DG formulations, without the need to devise additional operators to perform the coupling as long as the interfaces' cuts fall on the cell boundaries.

Another less explored feature of DD-ROMs is the possibility to repartition the computational domain, while keeping the data structures relative to each subdomain local in memory, with the aim of obtaining more efficient or accurate ROMs. In fact, one additional reason to subdivide the computational domain is to partition the solution manifold into local solution manifolds that have a faster decay of the Kolmogorov n-width. The repartition of the computational domain can be performed with \textit{ad hoc} domain decomposition strategies. To our knowledge, the only case found in the literature is introduced in~\cite{xiao2019domainb}, where the degrees of freedom are split in each subdomain minimising the communication and activity between them and balancing the computational load across them. Relevant is the choice of weights to assign to each degree of freedom: uniform weights, nodal values of Reynolds stresses for the turbulent Navier-Stokes equations or the largest singular value of the discarded local POD modes. In particular, the last option results in a balance of energy in $L^2$ norm retained in each subdomain. We explore a different approach.

It must be remarked that in any case, fixed the value of the reconstruction error of the training dataset in Frobenious norm $\lVert\cdot\rVert_{F}$, there must be at least a local reduced basis dimension greater or equal to the global reduced basis space dimension. In fact, for the Eckhart-Young theorem, if $X\in\mathbb{R}^{d\times n}$ is the snapshots matrix ordered by columns, we have that the projection into the first $k$ modes $\{v_k\}_{i=1}^k$ achieves the best approximation error in the Frobenious norm  in the space of matrices of rank $k$:
\begin{equation}
    P_k = \argmin_{P\in\mathbb{R}^{m\times n} s.t.\ r(P)=k}\lVert X- PX\rVert_{F},\quad P_k = \sum_{i=1}^{k}v_i\otimes v_i,
\end{equation}
where $r(\cdot)$ is the matrix rank. So, in general, it is not possible to achieve a better training approximation error in the Frobenious norm $\lVert\cdot\rVert_{F}$ employing a number of local reduced basis smaller than what would be needed to achieve the same accuracy with a global reduced basis. So, differently from~\cite{xiao2019domainb}, instead of balancing the local reduced basis dimension among subdomains, we repartition the computational domain in regions whose restricted solution manifold is easily approximable by linear subspaces and regions for which more modes are needed. Anyway, for truly decomposable systems we expect that the reconstruction error on the test set is lower when considering local reduced basis instead of global ones, as will be shown for the Maxwell equation in stationary regime test case with discontinuous piecewise constant parameters, see Figure~\ref{fig:a_posteriori_torus_disc}. 

\subsection{Implementation of Domain Decomposable ROMs}
\label{subsec:decoposableROMS}
Let us assume that the full-order model is implemented in parallel with distributed memory parallelism with $K>1$ cores, i.e., each $i$-th core owns locally the data structures relevant only to its assigned subdomain $\Omega_i$ of the whole computational domain $\cup_{i=1}^K \Omega_i = \Omega_h\subset\mathbb{R}^d$, for $i=1,\dots,K$. We will employ the \texttt{deal.II}  library~\cite{dealII93} to discretize the FS with the DG method, assemble the associated linear systems and solve them in parallel~\cite{BBHK11}. In particular, we employ \texttt{p4est}~\cite{BursteddeWilcoxGhattas11} to decompose the computational domain, \texttt{PETSc}~\cite{petsc-web-page} to assemble the linear system and solve it at the full-order level and \texttt{petsc4py}~\cite{DalcinPazKlerCosimo2011} to assemble and solve the reduced order system. At the offline and online stages the computations are performed in a distributed memory setting in which each core assembles its own affine decomposition, so that the evaluation of the reduced basis and of the projected local operators is always performed in parallel.

The weak formulation~\eqref{eq:DGFS} is easily decomposable thanks to the additive properties of the integrals. We recall the definition of the weak formulation $\forall y_h\in V_h,\ z_h\in V^*$
\begin{align}
	\begin{split}
    a^{cf}_h(z_h, y_h) + s_h(z_h, y_h) = &\sum_{T \in \mathcal{T}_h}  (z_h,\tilde{A}y_h)_{L^2(T)} +{\tfrac{1}2 \sum_{F \in \mathcal{F}^b_h}\left((\mathcal{M}+\mathcal{D})z_h,y_h\right)_{L^2(F)}} +{\sum_{F \in \mathcal{F}^i_h} \left(\mathcal{D}_F\lbrace\!\lbrace z_h \rbrace \!\rbrace,[\![y_h]\!]  \right)_{L^2(F)}}+\\
    &\sum_{F\in \mathcal{F}^{b}_h} (S^b_Fz_h,y_h)_{L^2(F)} + \sum_{F\in \mathcal{F}^i_h} (S^i_h [\![z_h]\!],[\![y_h]\!])_{L^2(F)},
	\end{split}
\end{align}
and we decompose it into the $K$ subdomains as
\begin{align}
	\begin{split}
     a^{cf}_h(z_h, y_h) + s_h(z_h, y_h) = &\sum_{i=1}^{K}\left(\sum_{T \in \mathcal{T}_{h, i}}  (z_h,\tilde{A}y_h)_{L^2(T)} +{\tfrac{1}2 \sum_{F \in \mathcal{F}^b_{h, i}}\left((\mathcal{M}+\mathcal{D})z_h,y_h\right)_{L^2(F)}} +{\sum_{F \in \mathcal{F}^i_{h, i}} \left(\mathcal{D}_F\lbrace\!\lbrace z_h \rbrace \!\rbrace,[\![y_h]\!]  \right)_{L^2(F)}}+\right.\\
    & \left.\sum_{F\in \mathcal{F}^{b}_{h, i}} (S^b_Fz_h,y_h)_{L^2(F)} + \sum_{F\in \mathcal{F}^i_{h, i}} (S^i_h [\![z_h]\!],[\![y_h]\!])_{L^2(F)}\right) + \sum_{\substack{i=1\\j=i}}^{K} \left({\sum_{F \in \mathcal{F}^i_{h, i, j}} \left(\mathcal{D}_F\lbrace\!\lbrace z_h \rbrace \!\rbrace,[\![y_h]\!]  \right)_{L^2(F)}}+\right.\\
    &\left.\sum_{F\in \mathcal{F}^i_{h, i, j}} (S^i_h [\![z_h]\!],[\![y_h]\!])_{L^2(F)}\right),
   	\end{split}
\end{align}
\begin{align}
    l_h (y_h) = &\sum_{T \in \mathcal{T}_h}  (f,y_h)_{L^2(T)} + {\tfrac{1}2 \sum_{F \in \mathcal{F}^b_h}\left((M-D)g,y_h\right)_{L^2(F)}}=\sum_{i=1}^K\left(\sum_{T \in \mathcal{T}_{h, i}}  (f,y_h)_{L^2(T)} + {\tfrac{1}2 \sum_{F \in \mathcal{F}^b_{h, i}}\left((M-D)g,y_h\right)_{L^2(F)}}\right),
\end{align}
where we have defined the internal subsets $\mathcal{T}_{h, i}=\mathcal{T}_{h}\cap\Omega_i$, $\mathcal{F}^i_{h, i}=\mathcal{F}^i_{h}\cap\mathring{\Omega}_i$ and $\mathcal{F}^b_{h, i}=\mathcal{F}^b_{h}\cap\mathring{\Omega}_i$, $\forall i=1,\dots,K$ and the interfaces subsets $\mathcal{F}^i_{h, i, j}=\mathcal{F}^i_{h}\cap\overline{\Omega}_i\cap\overline{\Omega}_j$ and $\mathcal{F}^b_{h, i, j}=\mathcal{F}^b_{h}\cap\overline{\Omega}_i\cap\overline{\Omega}_j$, $\forall i=1,\dots,K$. We remark that the computational domain is always decomposed such that the cuts of the subdomains $\{\partial\Omega_i\}_{i=1}^K$ fall on the interfaces of the triangulation $\mathcal{F}^i_h\cup \mathcal{F}^b_h$.

We define the bilinear and linear operators in $V_h^*$,
\begin{align}
	\begin{split}
    \mathcal{A}_{ii} = & \sum_{T \in \mathcal{T}_{h, i}}  (\bullet,\tilde{A}\bullet)_{L^2(T)} +{\tfrac{1}2 \sum_{F \in \mathcal{F}^b_{h, i}}\left((\mathcal{M}+\mathcal{D})\bullet,\bullet\right)_{L^2(F)}} +{\sum_{F \in \mathcal{F}^i_{h, i}} \left(\mathcal{D}_F\lbrace\!\lbrace \bullet \rbrace \!\rbrace,[\![\bullet]\!]  \right)_{L^2(F)}}+\\
    &\sum_{F\in \mathcal{F}^{b}_{h, i}} (S^b_F\bullet,\bullet)_{L^2(F)} + \sum_{F\in \mathcal{F}^i_{h, i}} (S^i_h [\![\bullet]\!],[\![\bullet]\!])_{L^2(F)},\quad\forall i = 1,\dots, K,
	\end{split}\\
    \mathcal{A}_{ij} = \mathcal{A}_{ji} = & \sum_{F \in \mathcal{F}^i_{h, i, j}} \left(\mathcal{D}_F\lbrace\!\lbrace \bullet \rbrace \!\rbrace,[\![\bullet]\!]  \right)_{L^2(F)}+\sum_{F\in \mathcal{F}^i_{h, i, j}} (S^i_h [\![\bullet]\!],[\![\bullet]\!])_{L^2(F)},\quad\forall j,i = 1,\dots, K,\  i\neq j,\\
    \mathcal{F}_i = &\sum_{T \in \mathcal{T}_{h, i}}  (f,\bullet)_{L^2(T)} + {\tfrac{1}2 \sum_{F \in \mathcal{F}^b_{h, i}}\left((M-D)g,\bullet\right)_{L^2(F)}},\quad\forall i = 1,\dots, K,
\end{align}
and their matrix representation in the discontinuous Galerkin basis of $V_h$,
\begin{equation}
    \left.(\mathcal{A}_{ii})\right|_{{V_h}} = A_{ii}, \qquad\left.\mathcal{F}_i\right|_{V_h} = F_i,\quad\forall i = 1,\dots, K,\qquad \left.(\mathcal{A}_{ij})\right|_{{V_h}} = \left.(\mathcal{A}_{ji})\right|_{{V_h}} = A_{ij} = A_{ji},\quad\forall j,i = 1,\dots, K,\  i\neq j,
\end{equation}
and in the local reduced basis $V_i = \lbrace \psiRB_{j, i}\rbrace_{j=1}^{\NRB}\subset V_h(\Omega_i),\ i=1,\dots,K$,
\begin{equation}
    \left.(\mathcal{A}_{ii})\right|_{{V_{RB}}} = B_{ii},\qquad\left.\mathcal{F}_i\right|_{V_{RB}} = L_i,\quad\forall i = 1,\dots, K,\qquad \left.(\mathcal{A}_{ij})\right|_{{V_{RB}}} = \left.(\mathcal{A}_{ji})\right|_{{V_{RB}}} = B_{ij} = B_{ji},\quad\forall j,i = 1,\dots, K,\  i\neq j .
\end{equation}

As anticipated, in our test cases the subdomains interface penalties are naturally included inside $\{\mathcal{A}_{ij}\}_{i,j=1,\dots,K}$. In practice, additional penalty terms could be implemented:
\begin{equation}
    \mathcal{S}_{ij} = \sum_{F\in\mathcal{F}^i_{h, i, j}}(S[\![\bullet]\!], [\![\bullet]\!])_{L^2(F)},\qquad \left.\mathcal{S}_{ij}\right|_{V_h} = S_{ij}\qquad \left.\left(\mathcal{A}_{ij}+\mathcal{S}_{ij}\right)\right|_{V_{RB}} = B_{ij},\quad \forall j,i = 1,\dots, K,\  i\neq j.
\end{equation}

\begin{figure}[h]
    \centering
    \includegraphics[width=1\textwidth]{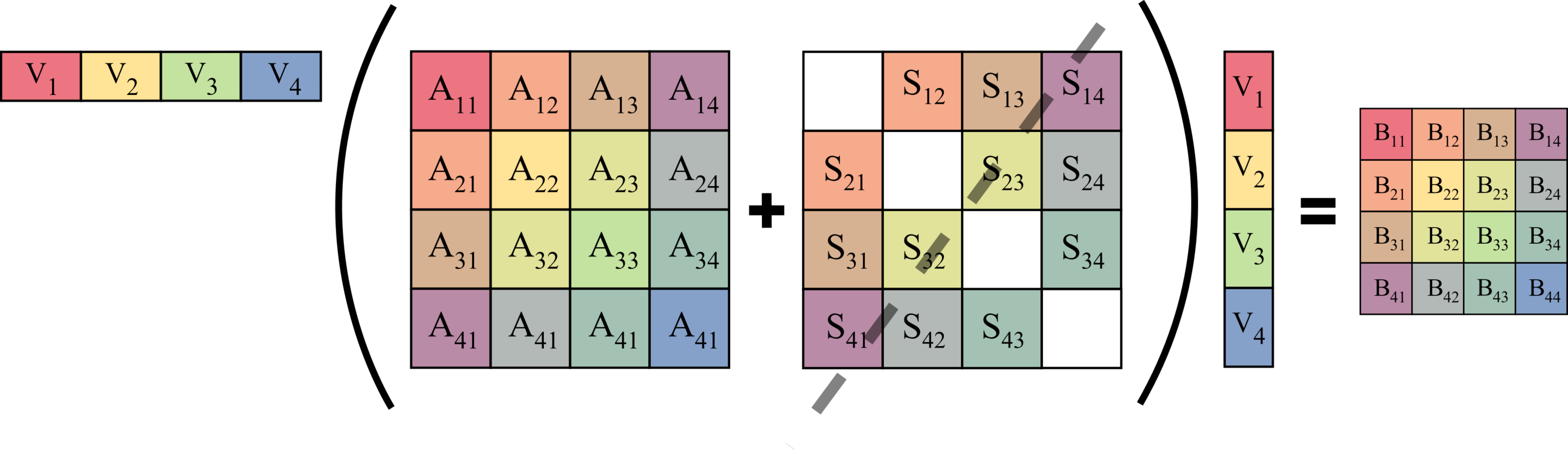}
    \caption{Assembly of the reduced block matrix $\{B_{i,j}\}_{i,j=1}^{4}$ through the projection onto the local reduced basis $\{V_{i}\}_{i=1}^{4}$ of the full-order partitioned matrix $A=\{A_{i,j}\}_{i,j=1}^{4}$ when considering $4$ subdomains. The natural DG penalty terms are included in the matrix $A$ without the need for additional penalty terms $\{S_{i,j}\}_{i\neq j,\ ,\\i,j=1}^{4}$ to impose stability at the reduced level.}
    \label{fig:block_FS}
\end{figure}

A matrix representation of the projection of the full-order block matrix $(A_{ij})_{i, j=1}^K\in\mathbb{R}^{d\times d}$ into the reduced order block matrix $(B_{ij})_{i, j=1}^K\in\mathbb{R}^{K \NRB \times K \NRB}$ is shown in Figure~\ref{fig:block_FS} for $K=4$. We remark that, differently from continuous Galerkin formulations, the DG penalization on jumps across the interfaces is already enough to couple the subdomains and there is no need of further stabilization, as shown in Figure~\ref{fig:block_FS}. Nonetheless, additional interface penalties terms can be easily introduced, taking also into account DG numerical fluxes. The reduced dimension is the number of subdomains $K$ times the local reduced basis dimensions $\{\NRB_i\}_{i=1}^K$, here supposed equal $\NRB=\NRB_i,\ i=1,\dots,K$, but in general can be different.

\subsection{Repartitioning strategy}
A great number of subdomains can pollute the efficiency of the developed DD-ROMs at the online stage since the reduced dimension would be $\sum_{i=1}^K \NRB_i$ that scales linearly with the number of cores if the local reduced dimensions $\NRB_i$ are equal. In order to keep the computational savings in the assembly of the affine decomposition at the offline stage, we may want to preserve the distributed property of our ROM. One possible solution is to fix a reduced number of subdomains $k\ll K$ such that $\sum_{i=1}^k \NRB_i$ is small enough to achieve a significant speedup with respect to the FOM. The additional cost with respect to a monodomain ROM is associated to the evaluation of the $k$ local reduced basis with SVD and the assembly of the affine decomposition operators. The new $k$ reduced subdomains do not need to be agglomerations of the FOM subdomains, hence, different strategies to assemble the new $k$ reduced subdomains can be investigated.

The number of subdomains $K$ was kept the same as the FOM since it is necessary to collect the snapshots efficiently at the full-order level through \textbf{\textit{p4est}}. However, if we decide to repartition our computational domain, we can develop decomposition strategies that reduce $\sum_{i=1}^K \NRB_i$. Ideally, having in mind the Eckhart-Young theorem, a possible strategy is to lump together all the dofs of the cells that have a fast decaying Kolmogorov n-width, and focus on the remaining ones. We test this procedure in the practical case $k=2,\ K=4$ to perform numerical experiments in~\cref{subsec:numerical_experiments}.

To solve the classification problem of partitioning the elements of the mesh into $k$ subdomains, we describe here two scalar indicators that will be used as metrics. For $k=2$ subdomains, it will be sufficient to choose the percentage of cells $P_l$ corresponding to the lowest values of the chosen scalar indicator. Other strategies for $k>2$ may also involve clustering algorithms and techniques to impose connectedness of the clusters, as done for local dimension reduction in parameter spaces in~\cite{romor2021local}. A first crude and cheap indicator to repartition the computational domain is the cellwise variance of the training snapshots, as it measures how well, in mean squared error, the training snapshots are approximated by their mean, $\forall T\in\mathcal{T}_h$.
\begin{definition}[Cellwise variance indicator]
    We define the cellwise variace indicator $I_{var}:\mathcal{T}_h\rightarrow \mathbb{R}^+$,
    \label{def:ind_var}
    \begin{equation}
        I_{var}(T)=\int_{T}\lVert\text{Var}(\{\mathbf{z}(\boldsymbol{\rho}_i)\}_{i=1}^n)\rVert_{L^2(\mathbb{R}^m)}\ d\mathbf{x},\qquad \left(\text{Var}(\{\mathbf{z}(\boldsymbol{\rho}_i)\}_{i=1}^n)\right)_l=\tfrac{1}{n}\sum_{i=1}^n \left|\mathbf{z}_l(\boldsymbol{\rho}_i)-\tfrac{1}{n}\sum_{j=1}^n \mathbf{z}_l(\boldsymbol{\rho}_j)\right|^2,\quad l=1,\dots,m,
    \end{equation}
    where $n>0$ is the number of training DG solutions $\{\mathbf{z}(\boldsymbol{\rho}_i)\}_{i=1}^n$ with $\mathbf{z}(\boldsymbol{\rho}_i): \Omega\subset\R^d\to\mathbb{R}^{m},\ \forall i\in\{1,\dots, n\}$.
\end{definition}
Note that the indicator is a scalar function on the set of elements of the triangulation $\mathcal{T}_h$. This is possible thanks to the assumption that boundaries of the subdomains belong to the interfaces of the elements of $\mathcal{T}_h$. When this hypothesis is not fulfilled, we would need to evaluate additional operators to impose penalties at the algebraical interfaces between subdomains that are not included in the set $\mathcal{F}^i_h \cup\mathcal{F}^b_h$, not to degrade the accuracy.

The cellwise variance indicator is effective for all the test cases for which there is a relatively large region that is not sensitive to the parametric instances, as in our advection diffusion reaction test case in Section~\ref{subsubsec:adr_test}. Common examples are all the CFD numerical simulations that have a far field with fixed boundary conditions. However, the variance indicator may be blind to regions in which the snapshots can be spanned by a one or higher dimensional linear subspace and are not well approximated by a constant field, as in the compressible linear elasticity test case in Section~\ref{subsubsection:cle_test}.

In these cases, a valid choice is represented by a cellwise Grassmannian dimension indicator. We denote with $D_{\text{T}}$ the number of degrees of freedom associated to each element $T$, assumed constant in our test cases.

\begin{definition}[Cellwise Grassmannian dimension indicator]
    \label{def:ind_grassmannian}
    Fixed $1\leq \NRB_T\in\mathbb{N}$, and $1\leq n_{\text{neig}}\in\mathbb{N}$, we define the cellwise Grassmannian dimension indicator $I_{G}:\mathcal{T}_h\rightarrow \mathbb{R}^+$,
    \begin{equation}
        I_{G}(T)=\lVert X_T - U_T U_T^T X_T\rVert_F,
    \end{equation}
    where $X_T\in\mathbb{R}^{n_{\text{neig}}D_{T}\times n}$ is the snapshots matrix restricted to the cell $T$ and its $n_{\text{neig}}$ nearest neighbours, and $U_T\in\mathbb{R}^{n_{\text{neig}}D_{T}\times r_T}$ are the modes of the truncated SVD of $X_T$ with dimension $r_T$.
\end{definition}
The cellwise Grassmannian dimension indicator $I_G$ is a measure of how well the training snapshots restricted to a neighbour of each cell are approximated by a $\NRB_T$ dimensional linear subspace. Employing this indicator, we recover an effective repartitioning of the computational subdomain of the compressible linear elasticity test case, see Section~\ref{subsubsection:cle_test}. The Grassmannian indicator has two hyper-parameters that we fix for each test case in section~\ref{subsec:numerical_experiments}: the number of nearest-neighbour cells is $n_{\text{neigh}}=3$ and the number of reduced local dimension used to evaluate the $L^2$ reconstruction error is $\NRB_T=1$. The number of nearest-neighbour is chosen to deal with critical cases at the boundaries and the closest neighbouring cells are chosen based on the distance of barycenters. The reduced local  dimension $\NRB_T=1$ is chosen very small as the computations must be done on very few cells.

We remark that both indicators do not guarantee that the obtained subdomains belong to the same connected components and, though this might be a problem in terms of connectivity and computational costs for the FOM, at the reduced level this does not affect the online computational costs.
Nevertheless, in the tests we perform, the obtained subdomains are connected.

\begin{figure}[h]
    \centering
    \includegraphics[width=1\textwidth]{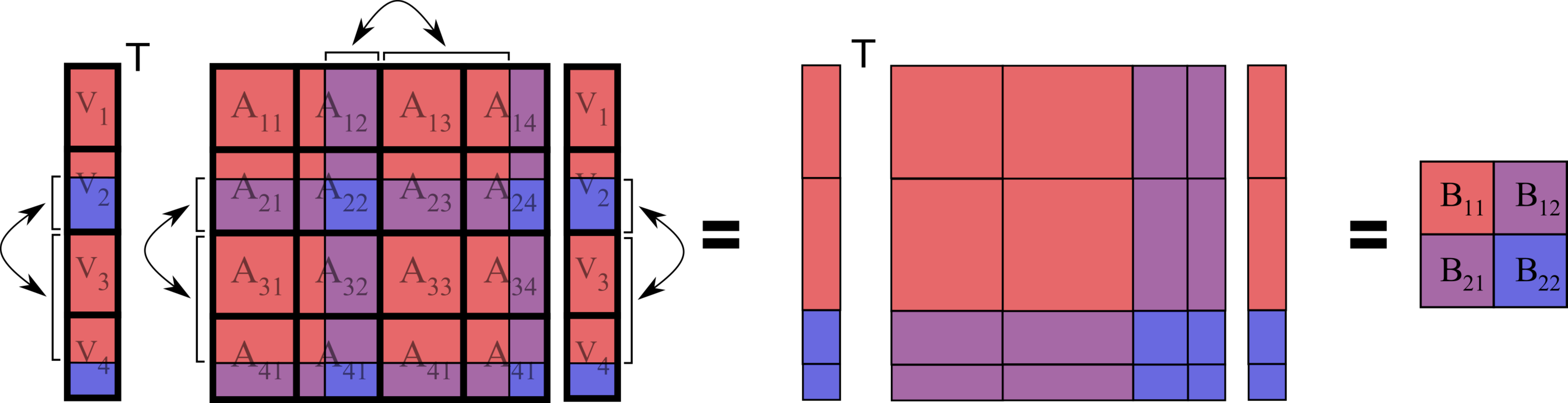}
    \caption{Repartitioning of the reduced block matrix shown in Figure~\ref{fig:block_FS} from $K=4$ subdomains to $k=2$ repartitioned subdomains. The projection from the full-order matrix $A=\{A_{i,j}\}_{i,j=1}^{4}$ to the reduced matrix $B=\{B_{i,j}\}_{i,j=1}^{2}$ is sketched. It is performed locally in a distributed memory setting, the re-ordering shown by the arrows is reported only to visually see the which block structure of the full-order matrix $A$ would correspond to the blocks of the reduced matrix $B$.}
    \label{fig:repartitioning}
\end{figure}

Now, the assembly of the affine decomposition proceeds as explained in Section~\ref{subsec:decoposableROMS} with the difference that at least one local reduced basis and reduced operator is split between at least $2$ subdomains/cores. A schematic block matrix representation of the procedure is shown in Figure~\ref{fig:repartitioning}.

\subsection{Numerical experiments}
\label{subsec:numerical_experiments}
In this section, we test the presented methodology for different linear parametric partial differential equations: the Maxwell equations in stationary regime in~\cref{subsubsec:max} (\textbf{MS}), the compressibile linear elasticity equations in~\cref{subsubsection:cle_test} (\textbf{CLE}) and the advection diffusion reaction equations in~\cref{subsubsec:adr_test} (\textbf{ADR}). We study two different parametrizations for the test cases \textbf{MS} and \textbf{CLE}: one with parameters that affect the whole domain \textbf{MS1} and \textbf{CLE1}, and one with parameters that affect independently different subdomains \textbf{MS2} and \textbf{CLE2}. We show a case in which DD-ROMs work effectively \textbf{MS2} and a case \textbf{CL2} in which the performance is analogous to single domain ROMs, even if the parameters have a local influence.

We test the effectiveness of the \textit{a posteriori} error estimates introduced in~\cref{sec:certified}, the accuracy of DD-ROMs for $K=4$ and the results of repartitioning strategies with $k=2$ subdomains. 
When performing a repartition of the computational domain $\Omega$ in subdomains $\{\Omega_i\}_{i=1}^{k}$ with reduced dimensions $\{r_{\Omega_i}\}_{i=1}^{k}$, we call the subdomains with lower values of the variance indicator $I_{\text{var}}$, see \cref{def:ind_var}, \textit{low variance} regions and with lower values of the Grassmannian indicator $I_{G}$, see \cref{def:ind_grassmannian}, \textit{low Grassmannian reconstruction error}. The complementary subdomains are the \textit{high variance} and \textit{high Grassmannian reconstruction error} regions, respectively. We show a case (\textbf{CLE1}) in which the Grassmannian indicator detects a better partition in terms of local reconstruction error with respect to the variance indicator.

We will observe that the relative errors in $R$-norm and energy norm and the $L^2$ relative error estimator and $L^2$ relative energy norm estimator are the most affected by the domain partitions.%, see figures~\ref{fig:a_posteriori_torus_disc} and~\ref{fig:a_posteriori_bar_disc}.

The open--source software library employed for the implementation of the full-order Friedrichs' systems discontinuous Galerkin solvers is \texttt{deal.II}~\cite{dealII93} and we have used piecewise $\mathbb P ^2$ basis functions in all simulations. The partition of the computational domain is performed in \texttt{deal.II} through the open--source \texttt{p4est} package~\cite{BursteddeWilcoxGhattas11}. The distributed affine decomposition data structures are collected in the offline stage and exported in the sparse \texttt{NumPy} format~\cite{harris2020array}. The reduced order models and the repartition of the computational domains are implemented in Python with MPI-based parallel distributed computing \texttt{mpi4py}~\cite{9439927} and \texttt{petsc4py}~\cite{petsc-web-page} for solving the linear full-order systems through \texttt{MUMPS}~\cite{mumps}, a sparse direct solver.

\subsubsection{Maxwell equations in stationary regime (MS)}
\label{subsubsec:max}
We consider the parametric Maxwell equations in the stationary regime in $d=3$ spatial dimensions, with $m=6$ equations, on a torus $\Omega\subset\mathbb{R}^3$  with inner radius $r=0.5$ and outer radius $R=2$ centered in $\boldsymbol{0}$ and lying along the $(x,z)$ plane:
\begin{equation}
    \left(
    \begin{array}{c}
        \mu\Hb+\nabla\times \Eb\\
        \sigma\Eb-\nabla\times \Hb
    \end{array}
    \right)
    = 
    \left(
    \begin{array}{c}
        \mathbf{g} \\
        \mathbf{f}
    \end{array}
    \right),
    \quad \forall\x\in\Omega,
\end{equation}
the tangential homogeneous boundary conditions $\mathbf{n}\times\Eb = \mathbf{0}$ are applied with the boundary operator~\eqref{eq:boundary_max}. We vary the parameters in the interval $\boldsymbol{\rho}=(\mu, \sigma)\in[0.5, 2]\times [0.5, 3]\subset\mathbb{R}^2$, leading to $\mu_0 = \min(\mu, \sigma)$.

\begin{figure}[h]
    \centering
    \includegraphics[width=0.48\textwidth, trim={300 60 400 0}, clip]{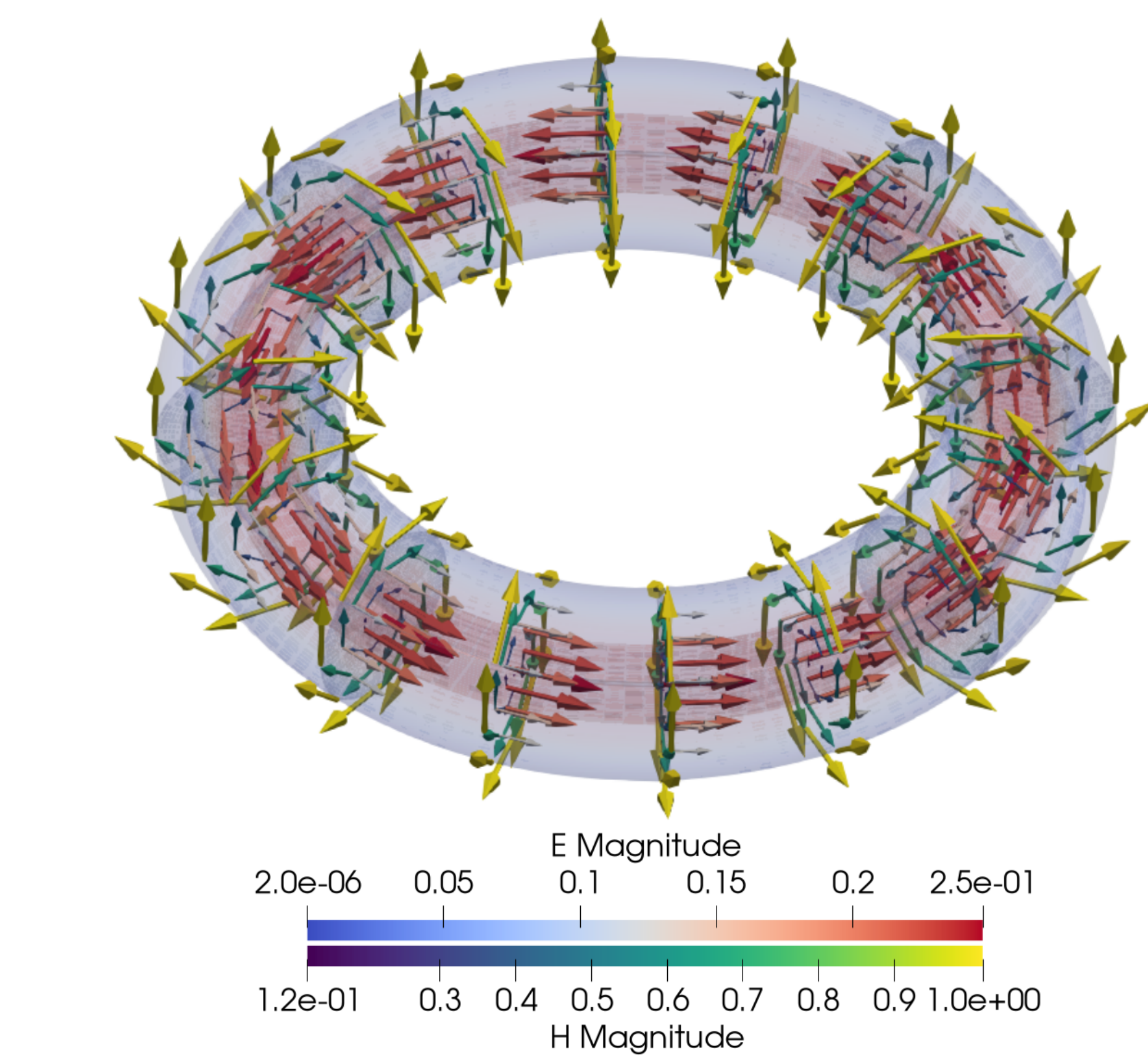}
    \includegraphics[width=0.48\textwidth, trim={300 60 400 0}, clip]{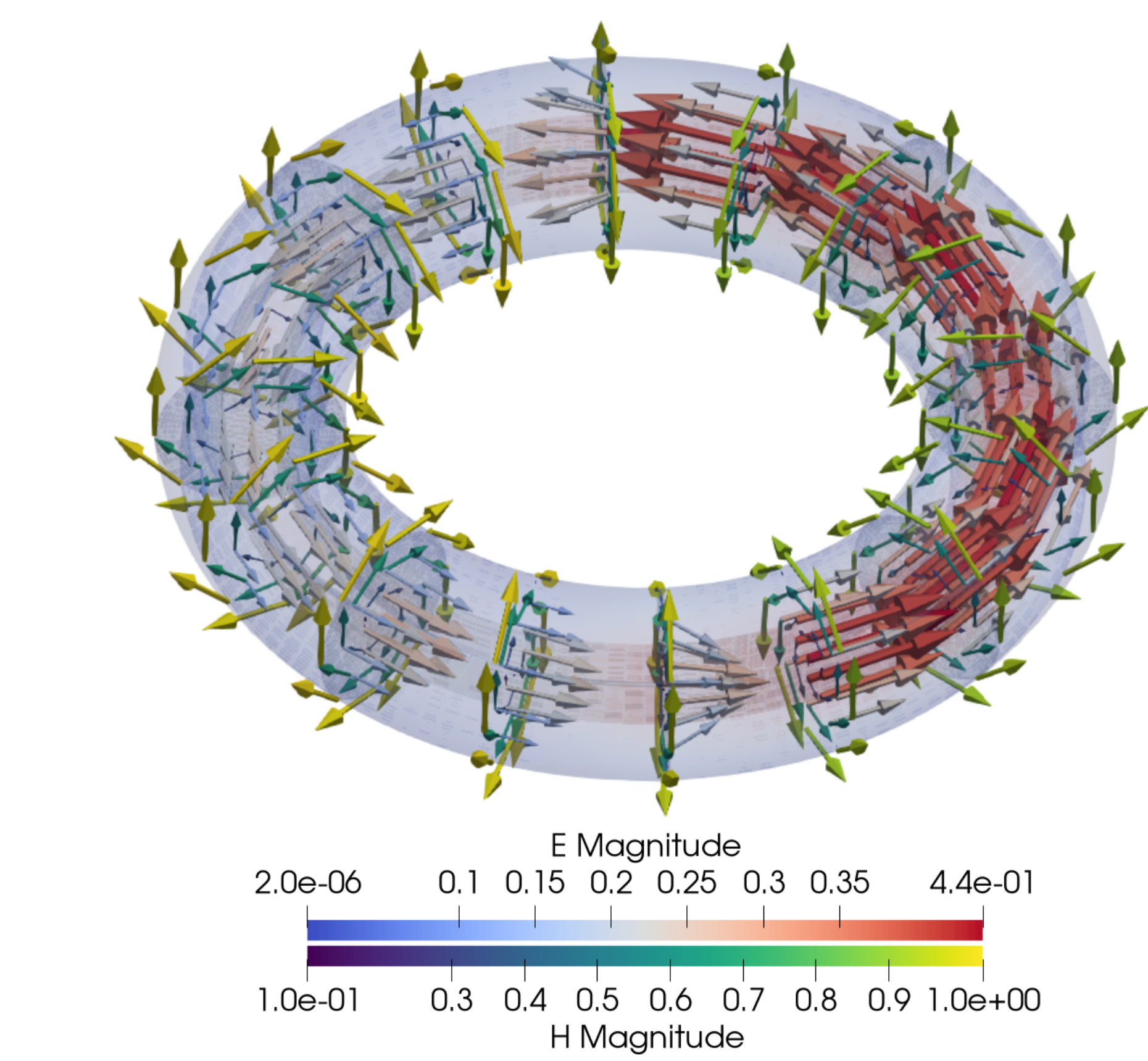}
    \caption{\textbf{MS}. Electric and magnetic fields of the Maxwell equations in stationary regime with Dirichlet homogeneous boundary conditions $\mathbf{n}\times\mathbf{E}=\mathbf{0}$. The vectors of the magnetic and electric fields are scaled by $0.5$ and $2$ of their magnitude respectively.  \textbf{Left}: \textbf{MS1}, $\mu=\sigma=1$, test case errors shown in Figure~\ref{fig:a_posteriori_torus}. \textbf{Right}: \textbf{MS2}, $\mu=\sigma=1$ in $\{x<0\}\cap\Omega$ and $\mu=\sigma=2$ in $\{x\geq0\}\cap\Omega$, test case errors shown in Figure~\ref{fig:a_posteriori_torus_disc}.}
    \label{fig:tori}
\end{figure}

\begin{figure}[h]
    \centering
    \includegraphics[width=0.32\textwidth, trim={120 20 120 150}, clip]{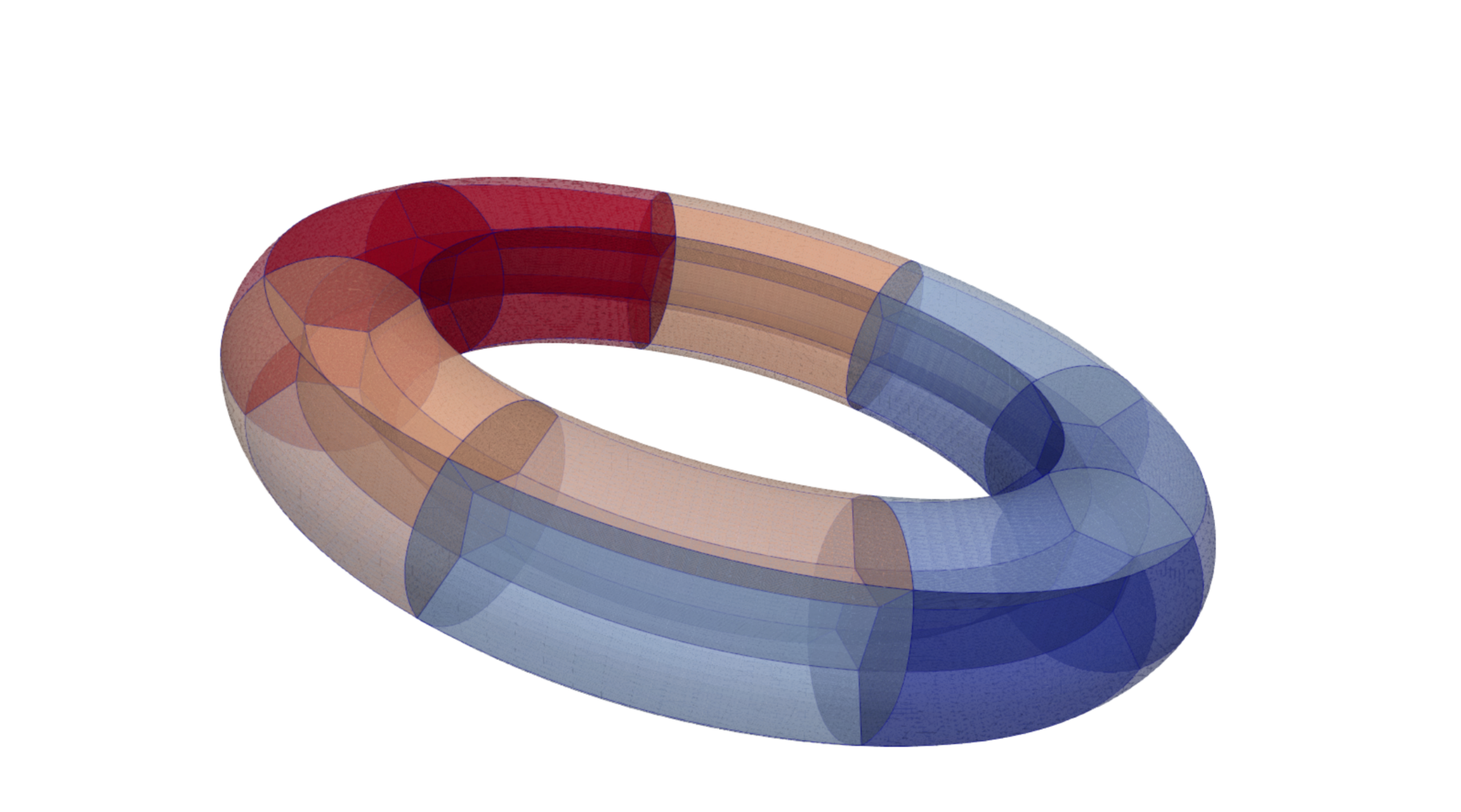}
    \includegraphics[width=0.32\textwidth, trim={120 20 120 150}, clip]{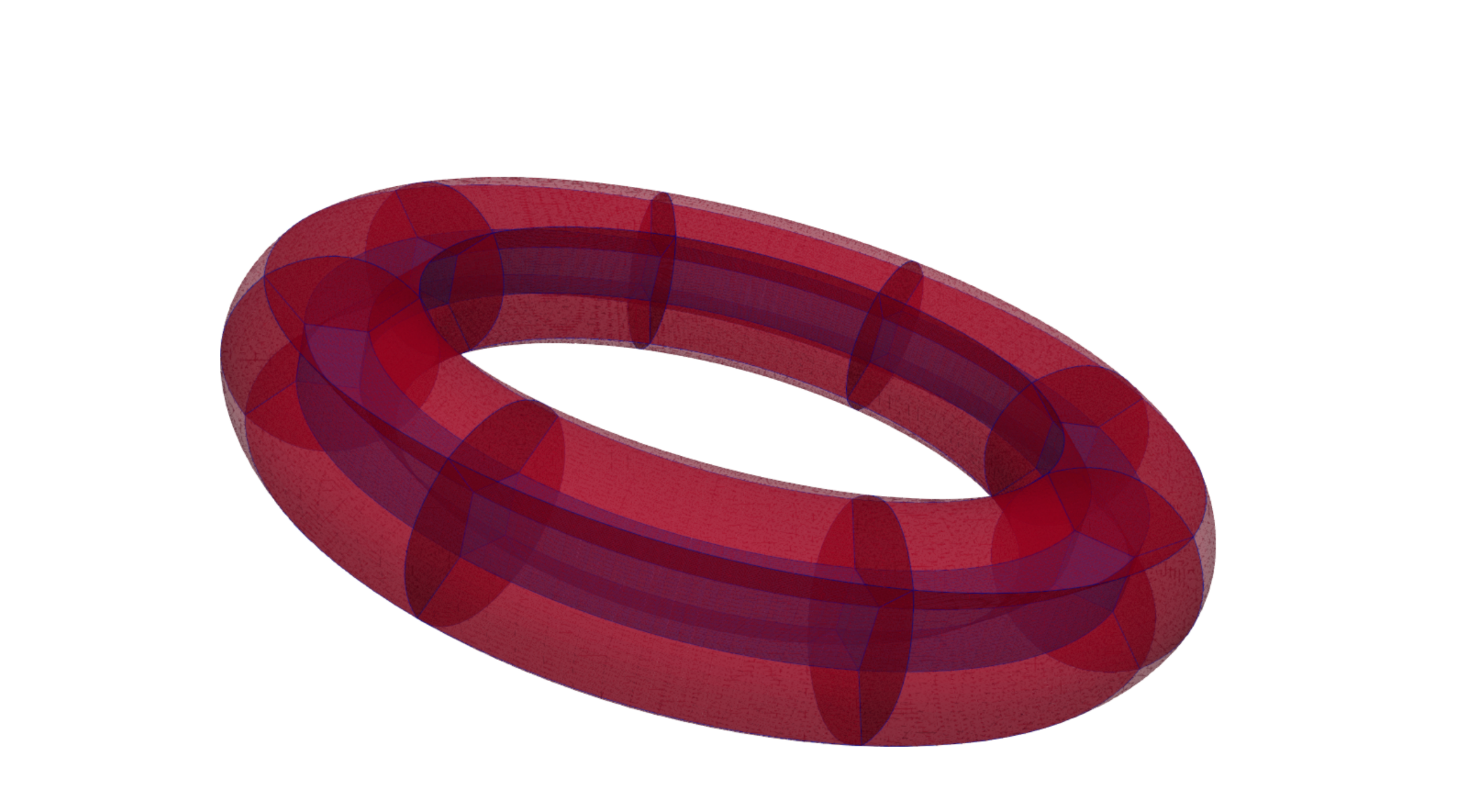}
    \includegraphics[width=0.31\textwidth, trim={120 30 120 150}, clip]{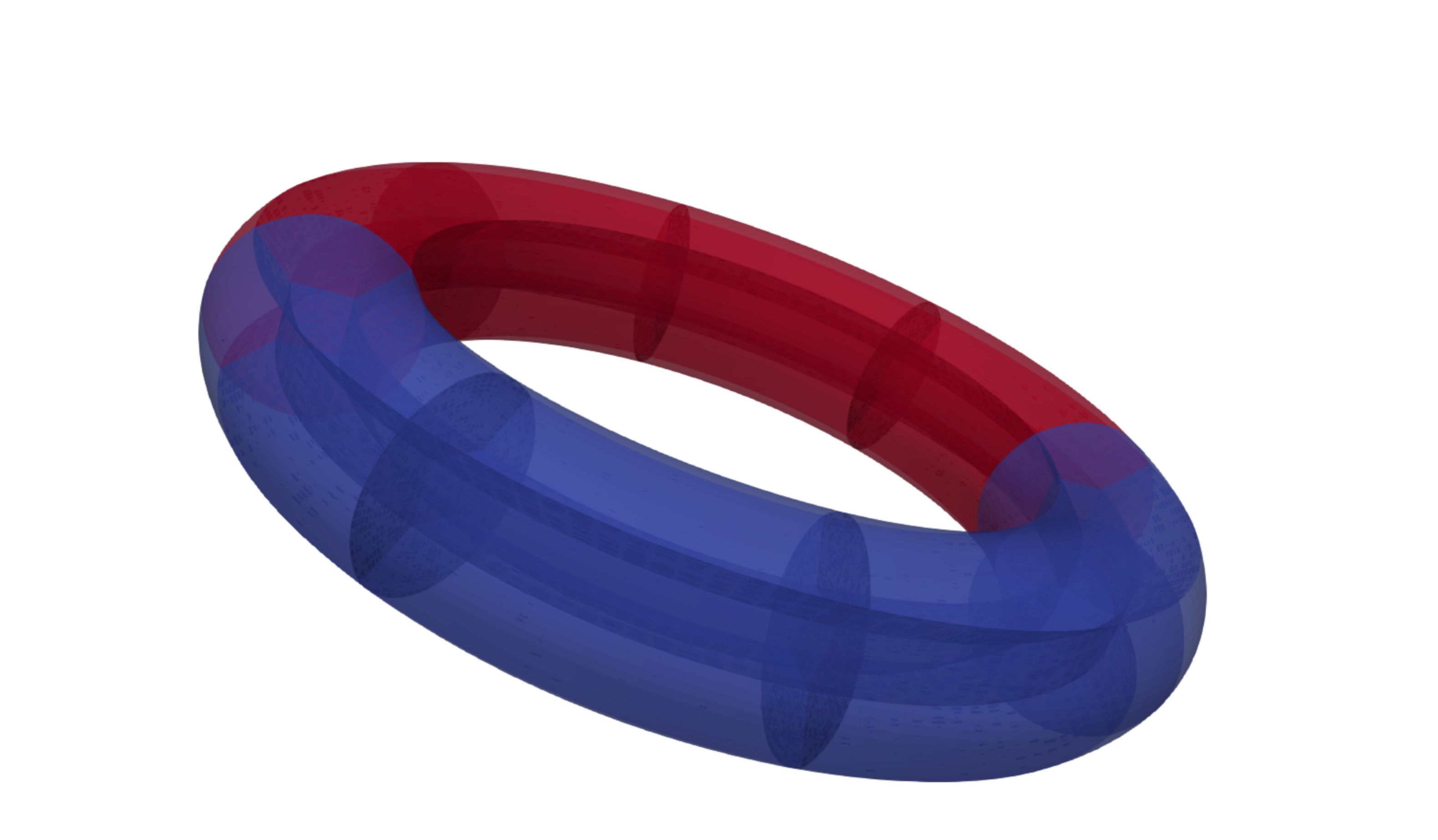}
    \caption{\textbf{MS}. \textbf{Left}: FOM computational domain partitioned in $K=4$ subdomains inside \texttt{deal.II}. \textbf{Center}: \textbf{MS1}, DD-ROM repartition of the computational subdomain $k=2$ with the cellwise variance indicator $I_{\text{var}}$, \cref{def:ind_var}: $20\%$ of the cells belong to the \textit{low variance} part, represented in blue inside the torus, and the other $80\%$ belong to the \textit{high variance} part, represented in red. \textbf{Right}: \textbf{MS2}, DD-ROM repartition with variance indicator $P_l=50\%$. The computational domain is exactly split at the interfaces that separate the subdomains $\{x<0\}\cap\Omega$ and $\mu=\sigma=2$ in which the parameters $\mu$ and $\sigma$ are constant.}
    \label{fig:subdomains_maxwell}
\end{figure}

We consider the exact solutions
\small
\begin{align*}
    &\mathbf{H}_{\text{exact}}(\x)=-\frac{1}{\mu}\left(\frac{2xy}{\sqrt{x^2+z^2}}, \frac{-4y^2\sqrt{x^2+z^2}+\sqrt{x^2+z^2}(-12(x^2+z^2)-15)+32(x^2+z^2)}{4(x^2+z^2)}, \frac{2xy}{\sqrt{x^2+z^2}}\right), \\
    &\mathbf{E}_{\text{exact}}(\x)=\left(\frac{z}{\sqrt{x^2+z^2}}, 0, -\frac{x}{\sqrt{x^2+z^2}}\right)\cdot\left(r^2-y^2-\left(R-\sqrt{x^2+z^2}\right)^2\right).
\end{align*}
\normalsize
We remark that the exact solutions can be approximated with a linear reduced subspace of dimension $1$, if we obtained the reduced basis with a partitioned SVD on the fields $(\mathbf{H}, \mathbf{E})$ separtely. We do not choose this approach and perform a monolithic SVD to test the convergence of the approximation with a DD-ROMs with respect to the local reduced dimensions.
The source terms are defined consequently as
\begin{equation}
    \mathbf{g}(\x)=0,\qquad \mathbf{f}(\x)=\sigma \mathbf{E}_{\text{exact}}-\nabla\times \mathbf{H}_{\text{exact}}.
\end{equation}
We consider two parametric spaces:
\begin{subequations}
\begin{align}
    \boldsymbol{\rho}=(\mu, \sigma)\in[0.5, 2]\times [0.5, 3]=\mathcal{P}_1\subset\mathbb{R}^2,\qquad (\mathbf{MS1})\\
    \boldsymbol{\rho}=(\mu_1, \sigma_1,\mu_2, \sigma_2)\in[0.5, 2]\times [0.5, 3]\times[0.5, 2]\times [0.5, 3]=\mathcal{P}_2\subset\mathbb{R}^4,\qquad (\mathbf{MS2}) 
\end{align}
\end{subequations}
where in the second case, the parameters $\mu$ and $\sigma$ are now piecewise constant:
\begin{equation}
    \mu=
    \begin{cases}
        \mu_1, \quad x<0,\\
        \mu_2, \quad x\geq 0,
    \end{cases}\quad
    \sigma=
    \begin{cases}
        \sigma_1, \quad x<0,\\
        \sigma_2, \quad x\geq 0,
    \end{cases}
\end{equation}
where $\mathbf{x}=(x,y,z)\in\Omega\subset\mathbb{R}^3$. In Figure~\ref{fig:tori}, we show solutions for $\mu=\sigma=1$ and for discontinuous values of the parameters: $\mu_1=\sigma_1=1$ in $\{x<0\}\cap\Omega$ and $\mu_2=\sigma_2=2$ in $\{x\geq0\}\cap\Omega$. The FOM partitioned and DD-ROM repartitioned subdomains are shown in Figure~\ref{fig:subdomains_maxwell}. For \textbf{MS1}, we choose the variance indicator to repartition the computational subdomain in two subsets: $P_l=20\%$ of the cells for the \textit{low variance} part and $80\%$ for the \textit{high variance} part. For \textbf{MS2}, we split the computational domain in two parts with the Grassmannian indicator and $P_l=50\%$.

\begin{figure}[!htpb]
    \centering
    \includegraphics[width=\textwidth]{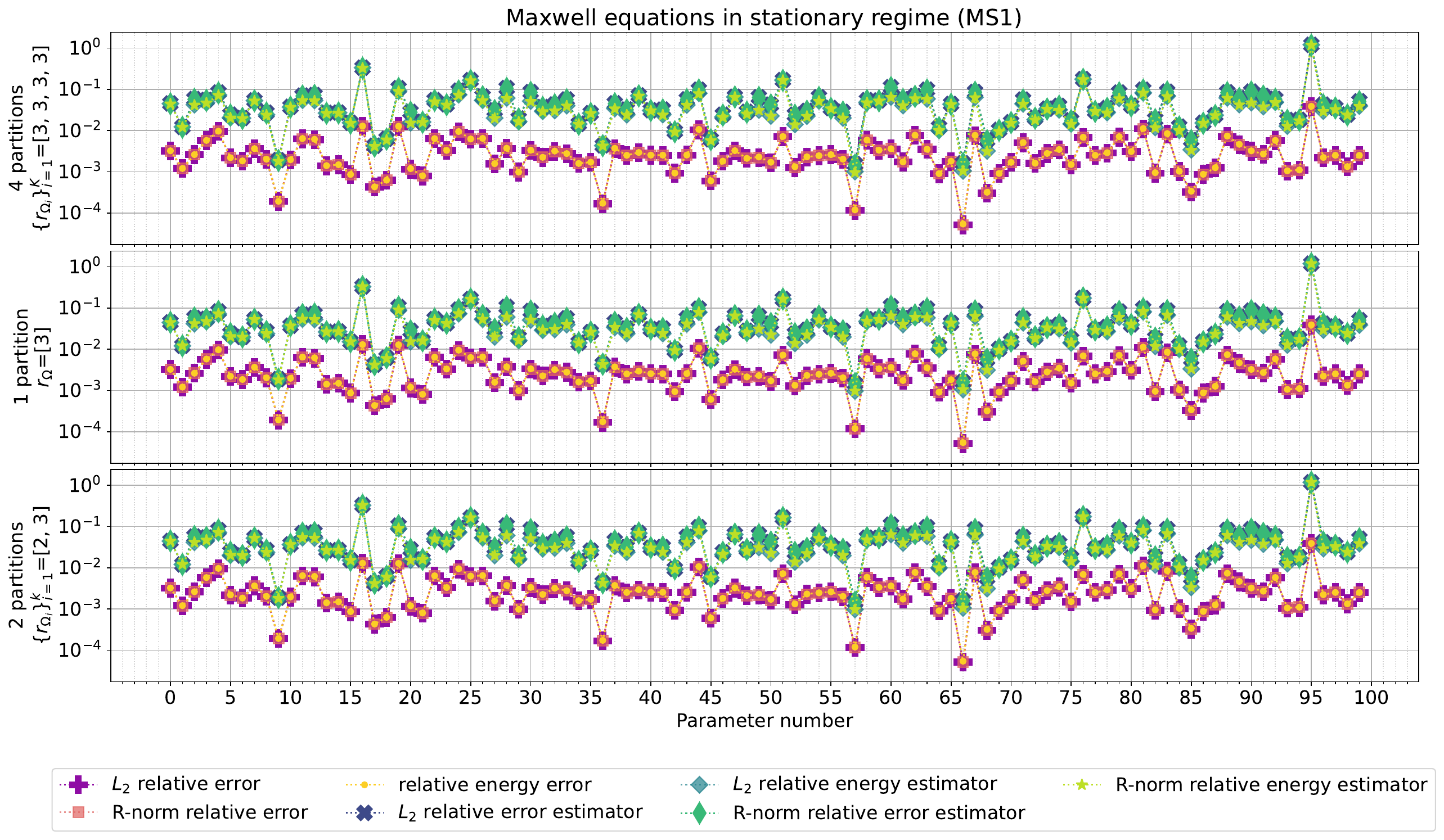}
    \caption{\textbf{MS1}. Errors and estimators for Maxwell equations corresponding to the $n_{\text{train}}=20$  uniformly sampled training snapshots corresponding to the abscissae $0,5,10,\dots,95$, and $n_{\text{test}}=80$ uniformly sampled test snapshots, corresponding to the other abscissae. The reduced dimensions of the ROMs are $\{\NRB_{\Omega_i}\}_{i=1}^K = [3, 3, 3, 3]$ for $K=4$ partitions, $\NRB_{\Omega} = 3$ for $k=1$ partition, and $\{\NRB_{\Omega_i}\}_{i=1}^k = [2, 3]$ for $k=2$ partitions. For the case $k=2$ we employed the cellwise variance indicator $I_G$, Definition~\ref{def:ind_var}, with $P_l=20\%$. It can be seen that even reducing the local dimension from $3$ to $2$ of one of the $k=2$ repartitioned subdomains, the accuracy of the predictions does not decrease sensibly.}
    \label{fig:a_posteriori_torus}
\end{figure}

\begin{figure}[!htpb]
    \centering
    \includegraphics[width=\textwidth]{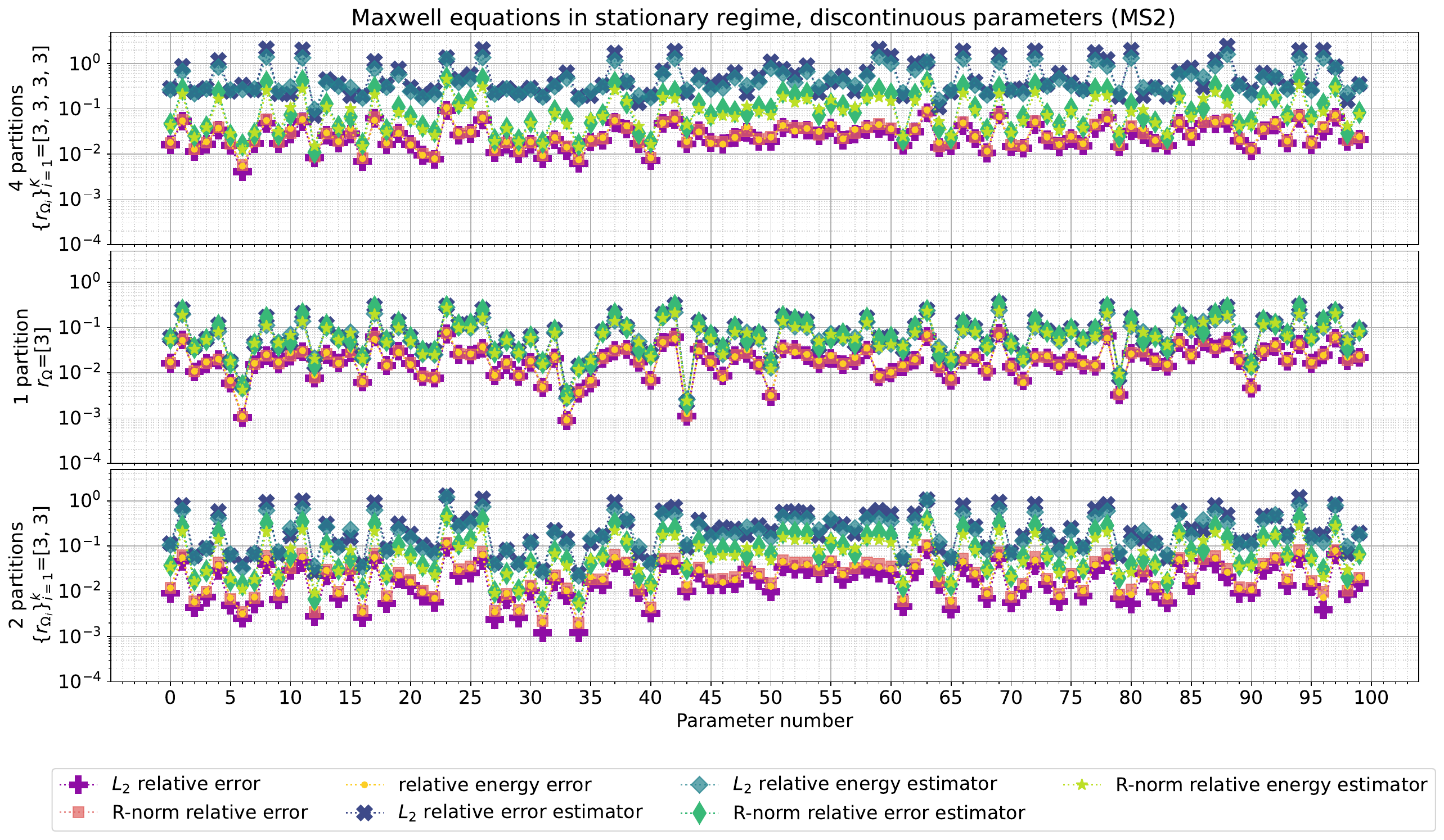}
    \caption{\textbf{MS2}. Errors and estimators for Maxwell equations with discontinuous $\mu$ and $\sigma$ corresponding to the $n_{\text{train}}=20$  uniformly sampled training snapshots corresponding to the abscissae $0,5,10,\dots,95$, and $n_{\text{test}}=80$ uniformly sampled test snapshots, corresponding to the other abscissae. The reduced dimensions of the ROMs are $\{\NRB_{\Omega_i}\}_{i=1}^K = [3, 3, 3, 3]$ for $K=4$ partitions, $\NRB_{\Omega} = 3$ for $k=1$ partition, and $\{\NRB_{\Omega_i}\}_{i=1}^k = [3, 3]$ for $k=2$ partitions. For the case $k=2$ we employed the cellwise local Grassmannian dimension indicator $I_G$, Definition~\ref{def:ind_var}, with $P_l=50\%$. The subdivisions detected exactly match the subdomains $\{x< 0\}\cap\Omega$ and $\{x\geq 0\}\cap\Omega$ in which the parameters are constant. An improvement of the predictions can be appreciated for some test parameters when employing $k=2$ repartitions.}
    \label{fig:a_posteriori_torus_disc}
\end{figure}

At the end of this subsection a comparison of the effectiveness of DD-ROMs with and without discontinuous parameters will be performed, the associated error plots are reported in Figure~\ref{fig:a_posteriori_torus} and Figure~\ref{fig:a_posteriori_torus_disc}. We will see that, for this simple test case \textbf{MS2}, there is an appreciable improvement of the accuracy when the computational domain subdivisions match the regions $\{x<0\}\cap\Omega$ and $\{x\geq 0\}\cap\Omega$ in which $\mu$ and $\sigma$ are constant. Such subdivision is detected by the Grassmannian indicator with $P_l=50\%$, as shown in Figure~\ref{fig:subdomains_maxwell} on the right. This is the archetypal
case in which DD-ROMs are employed successfully, in comparison with \textbf{MS1} for which there is no significant improvement with respect to classical global linear reduced basis.

\begin{figure}[!htpb]
    \centering
    \includegraphics[width=1\textwidth]{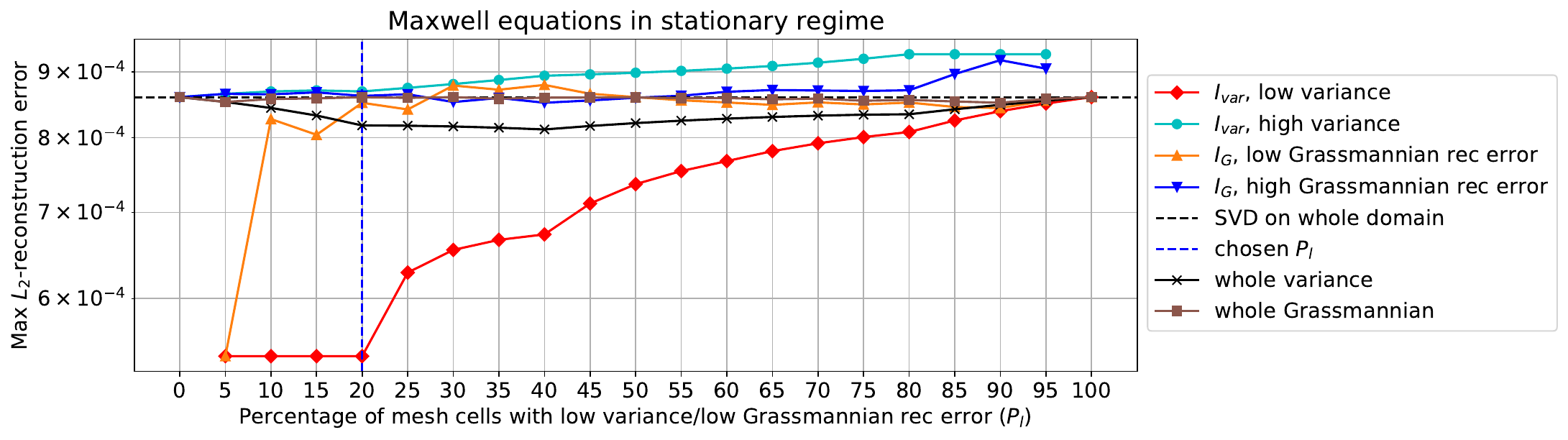}
    \caption{\textbf{MS1}. Local relative $L^2$-reconstruction errors of the snapshots matrix restricted to the two subdomains of the repartitioning performed with the indicator $I_{\text{var}}$ (in red and light-blue), Definition~\ref{def:ind_var}, and $I_{\text{G}}$ (in orange and blue), Definition~\ref{def:ind_grassmannian}. The relative $L^2$-reconstruction error attained on the whole domain is shown in black for the indicator $I_{\text{var}}$ and in brown for the indicator $I_{\text{G}}$. The local reduced dimensions used to evaluate the local reconstruction errors is $\NRB_{\Omega_i}=3,\ i=1,2$.}
    \label{fig:repartitioning_maxwell}
\end{figure}

In Figure~\ref{fig:repartitioning_maxwell}, we show how the different thresholds applied to the two indicators can affect the reconstruction error on a reduced space with $\NRB_{\Omega_i}=3$. All the lines plot the local relative error computed on different subdomains (either one of the $k$ DD-ROM subdomains or on the whole domain). 
On the $x$-axis it is shown the percentage of cells that are grouped into the low variance or low Grassmannian DD-ROM subdomain.
We observe that the cellwise variance indicator is a good choice for the purpose of repartitioning the subdomain from $K=4$ to $k=2$. 
Indeed, it is possible to build a low variance subdomain (value of the abscissa $20\%$ in Figure~\ref{fig:repartitioning_maxwell}) with a low local relative reconstruction error ($5\cdot 10^{-4}$) with respect to the global one ($8.6\cdot 10^{-4}$). This means that choosing the threshold $P_l = 20\%$ for the low variance subdomain, we should be able to use less reduced basis functions for that subdomain without affecting too much the global error.

Test case \textbf{MS1}. We evaluate $n_{\text{train}}=20$ training full-order solutions and $n_{\text{test}}=80$ test full-order solutions, corresponding to a uniform independent sampling from the parametric domain $\mathcal{P}_1\subset\mathbb{R}^2$. Figure~\ref{fig:a_posteriori_torus} shows the result relative to the relative $L^2$-error and relative errors in energy norm, with associated \textit{a posteriori} estimators. The numberd abscissae $0, 5, 10, \dots, 95$ represents the train parameters $n_{\text{train}}=20$ while the others $n_{\text{test}}=80$ parameters are the test set. 
For these studies, we have fixed the local reduced dimensions to $\NRB_{\Omega_i}=3,\ i=1,\dots,K$ for $K=4$, $\NRB_{\Omega}=3$ for the whole computational domain and $\NRB_{\Omega_1}=2,\ \NRB_{\Omega_2}=3$ for the DD-ROM repartitioned case with $k=2$. 
This choice of repartitioning with the $20\%$ of \textit{low variance} cells and local reduced dimension $\NRB_{\Omega_1}=2$ does not deteriorate significantly the accuracy and the errors almost coincide for all approaches.
However, unless the parameters $\sigma,\mu$ assume different discontinuous values in the computational domain $\Omega$, DD-ROMs are not advisable for this test case if the objective is improving the predictions' accuracy.

%For a comparison, see the results with the advection diffusion reaction test case in Section~\ref{subsubsec:adr_test}: in Figure~\ref{fig:repartitioning_maxwell} the gain in local $L^2$-reconstruction error is relatively small compared to the one in Figure~\ref{fig:repartitioning_channel}, where it reaches different orders of magnitude.

Test case \textbf{MS2}. Similarly to the previous case, we evaluate $n_{\text{train}}=20$ training full-order solutions and $n_{\text{test}}=80$ test full-order solutions, corresponding to a uniform independent sampling from the parametric domain $\mathcal{P}_2\subset\mathbb{R}^4$. As mentioned above, if we vary the parameters $\boldsymbol{\rho}=(\mu,\sigma)$ discontinuously on the subdomains $\{x\geq 0\}\cap\Omega$ and $\mathbf{x}\in\{x<0\}\cap\Omega$, we obtain the results shown in Figure~\ref{fig:a_posteriori_torus_disc}. It can be seen that repartitioning $\Omega$ in $k=2$ DD-ROM subdomains with the local Grassmannian indicator $I_G$ and $P_l=50\%$ produces effective DD-ROMs compared to the case of a single reduced solution manifold for the whole computational domain and for the DD-ROM with $k=4$ for which the subdomains do not match $\{x< 0\}\cap\Omega$ and $\{x\geq 0\}\cap\Omega$. In this case, we kept the local dimension of DD-ROM repartitioned case with $k=2$ equal $\NRB_{\Omega_1}=\NRB_{\Omega_2}=3$. For this simple test case, there is an appreciable improvement for some test parameters in the accuracy for $k=2$ instead of $K=4$ or a classical global linear basis ROM.

In Table~\ref{tab:comp_time_maxwell}, we list the computational times and speedups for a simulation with the different methods. For an error convergence analysis with respect to the size of the reduced space, we refer to Appendix~\ref{app:convergence_ROM}.
\begin{table}[!htbp]
	\centering
    \caption{\textbf{MS1}. Average computational times and speedups for ROM and DD-ROM approaches for Maxwell equations. The speedup is computed as the FOM computational time over the ROM one. The FOM runs in parallel with $4$ cores, so ``FOM time'' refers to wallclock time.}\label{tab:comp_time_maxwell}	
	\begin{tabular}{|c|c|c|c|c|c|c|c|}
      \hline
      \multicolumn{2}{|c|}{FOM}&\multicolumn{3}{c|}{ROM}&\multicolumn{3}{c|}{DD-ROM}\\\hline 
      $N_h$ & time & $\NRB$ & time & speedup & $\NRB_i$ & time& speedup \\ \hline\hline
		6480 & 254.851 [ms] & 3 & 51.436 [$\mu$s] & $\sim 495$ & [3, 3, 3, 3] &  62.680 [$\mu$s] & $\sim 406$ \\
        \hline
	\end{tabular}
\end{table}

\subsubsection{Compressibile linear elasticity (CLE)}
\label{subsubsection:cle_test}
Next, we consider the parametric compressible linear elasticity system in $d=3$ physical dimensions with a cylindrical shell along the z-axis as domain: the inner radius is $1$, outer radius $3$ and height $10$, and the base centered in $\boldsymbol{0}$. The $m=12$ equations of the FS are
\begin{equation}
    \left(
    \begin{array}{c}
        \boldsymbol{\sigma} - \mu_1(\nabla\cdot \ub)\mathbb{I}_3 -2\mu_2\frac{\left(\nabla \ub + (\nabla \ub)^t\right)}{2}\\
        -\frac{1}{2}\nabla\cdot\left(\boldsymbol{\sigma}+\boldsymbol{\sigma}^t\right) +\mu_3 \ub
    \end{array}
    \right)
    = 
    \left(
    \begin{array}{c}
        0 \\
        \mathbf{f}
    \end{array}
    \right),
    \quad \forall\x\in\Omega\subset\mathbb{R}^3,
\end{equation}
where $\boldsymbol{\rho}=(\mu_1, \mu_2, \mu_3)\in[100, 1000]^2\times [1,1]=\mathcal{P}\subset\mathbb{R}^3$ and $\mathbf{f} = (0, -1, 0)$. The system can be rewritten as FS as in \eqref{eq:elastic_FS}.
%:
%\begin{equation}
%    \left(
%    \begin{array}{c}
%        \boldsymbol{\sigma} - \frac{\mu_1}{2\mu_2+3\mu_1}\text{tr}(\boldsymbol{\sigma})\mathbb{I}_3 -\frac{\left(\nabla \ub + (\nabla \ub)^t\right)}{2}\\
%        -\frac{1}{2}\nabla\cdot\left(\boldsymbol{\sigma}+\boldsymbol{\sigma}^t\right) +\frac{\mu_3}{2\mu_2} \ub
%    \end{array}
%    \right)
%    = 
%    \left(
%    \begin{array}{c}
%        0 \\
%        \mathbf{f}
%    \end{array}
%    \right),
%    \quad \forall\x\in\Omega\subset\mathbb{R}^3.
%\end{equation}
We define the boundaries
\begin{equation}
    \Gamma_D = \partial\Omega\cap\{z=0\},\qquad \Gamma_N = \partial\Omega\setminus \Gamma_D.
\end{equation}
Mixed boundary conditions are applied with the boundary operator~\eqref{eq:boundary_bar}: homogeneous Dirichlet boundary conditions are imposed on $\Gamma_D$ and homogeneous Neumann boundary conditions on $\Gamma_N$.

We consider two parametric spaces:
\begin{subequations}
\begin{align}
    \boldsymbol{\rho}=(\mu_1, \mu_2)\in[100, 1000]^2=\mathcal{P}_1\subset\mathbb{R}^2,\qquad (\mathbf{CLE1})\\
    \boldsymbol{\rho}=(\mu_1, \mu_2, f_1, f_2)\in[100, 1000]^2\times [-2,2]^2=\mathcal{P}_2\subset\mathbb{R}^4,\qquad (\mathbf{CLE2}) 
\end{align}
\end{subequations}
where in the second case, the source term $\mathbf{f}$ is now piecewise constant:
\begin{equation}
    \mathbf{f}=
    \begin{cases}
        f_1\cdot (0, -1, 0), \quad z<5,\\
        f_2\cdot (0, -1, 0), \quad z\geq 5.
    \end{cases}
\end{equation}

\begin{figure}[ht!]
    \centering
    \includegraphics[width=0.49\textwidth, trim={100 50 120 70}, clip]{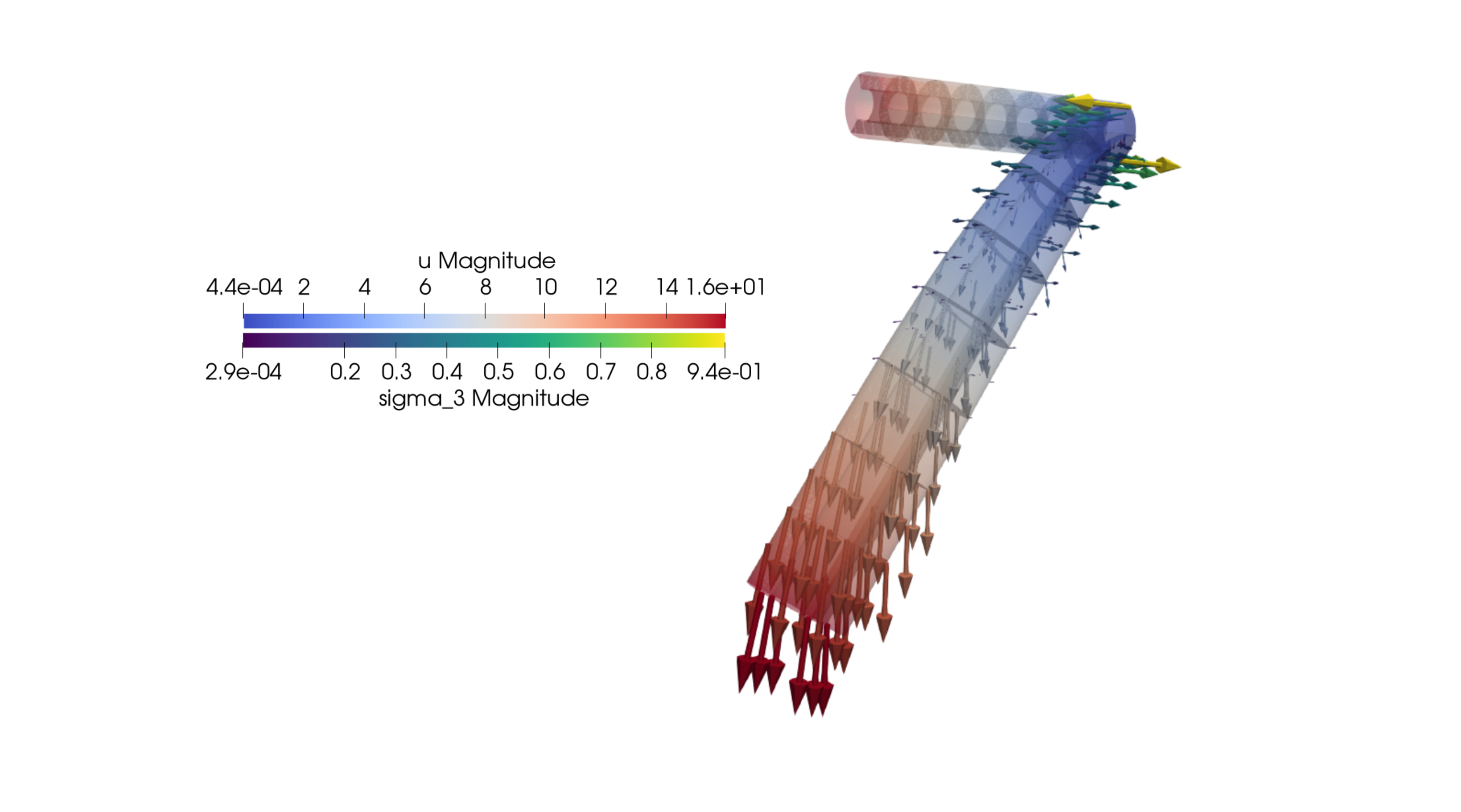}
    \includegraphics[width=0.44\textwidth, trim={0 120 300 40}, clip]{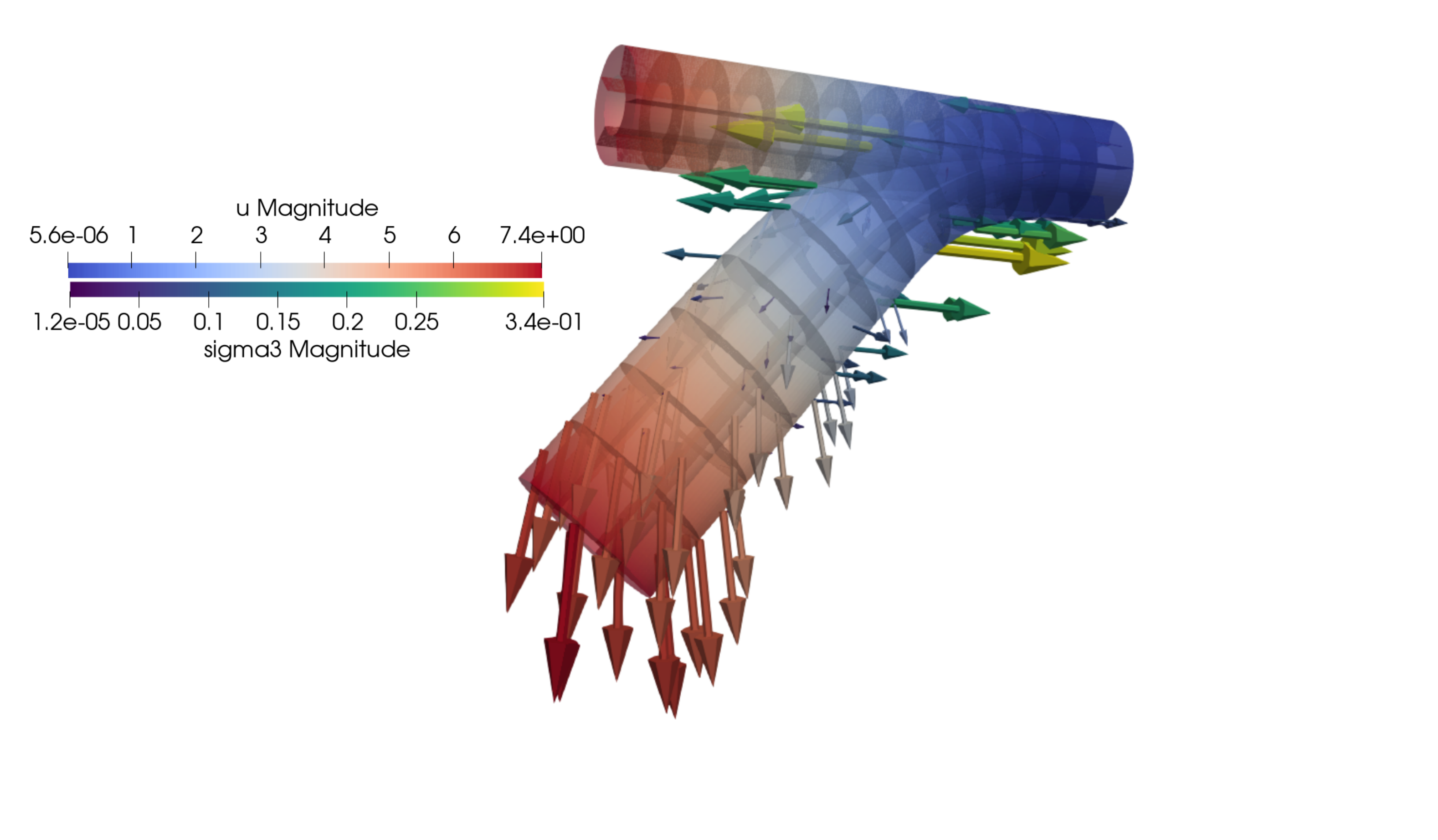}
    \caption{\textbf{CLE}. \textbf{Left: } solution of the compressible linear elasticity FS \textbf{CLE1} with parameter values $\mu_1=\mu_2=1000$. The cylindrical shell displacement $\mathbf{u}$, and with a different colorbar also the field $\sigma\mathbf{e_z}$, named \texttt{sigma\_3}, are shown. At the extremity close to $z=0$ homogeneous Dirichlet boundary conditions are imposed. \textbf{Right:} solution of the test case \textbf{CLE2} with discontinuous values of the source terms along the computational domain $\{z<5\}\cap\Omega$ and $\{z\geq 5\}\cap\Omega$: $\mu_1=\mu_2=1000$, $f_1=1$ and $f_2=-1$.}
    \label{fig:elasticity_sample}
\end{figure}

\begin{figure}[h]
    \centering
    \includegraphics[width=0.32\textwidth, trim={300 200 300 150}, clip]{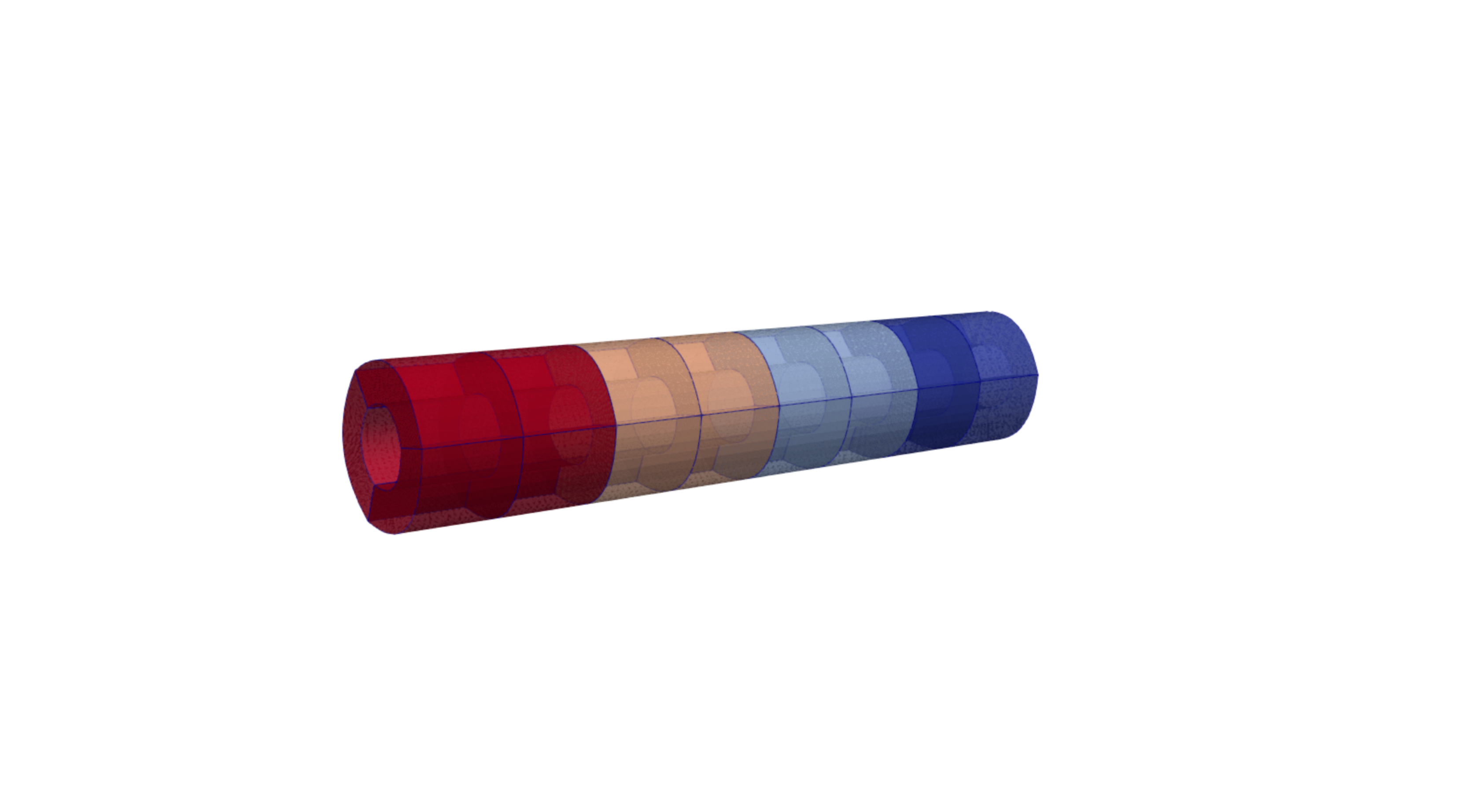}
    \includegraphics[width=0.32\textwidth, trim={300 200 300 150}, clip]{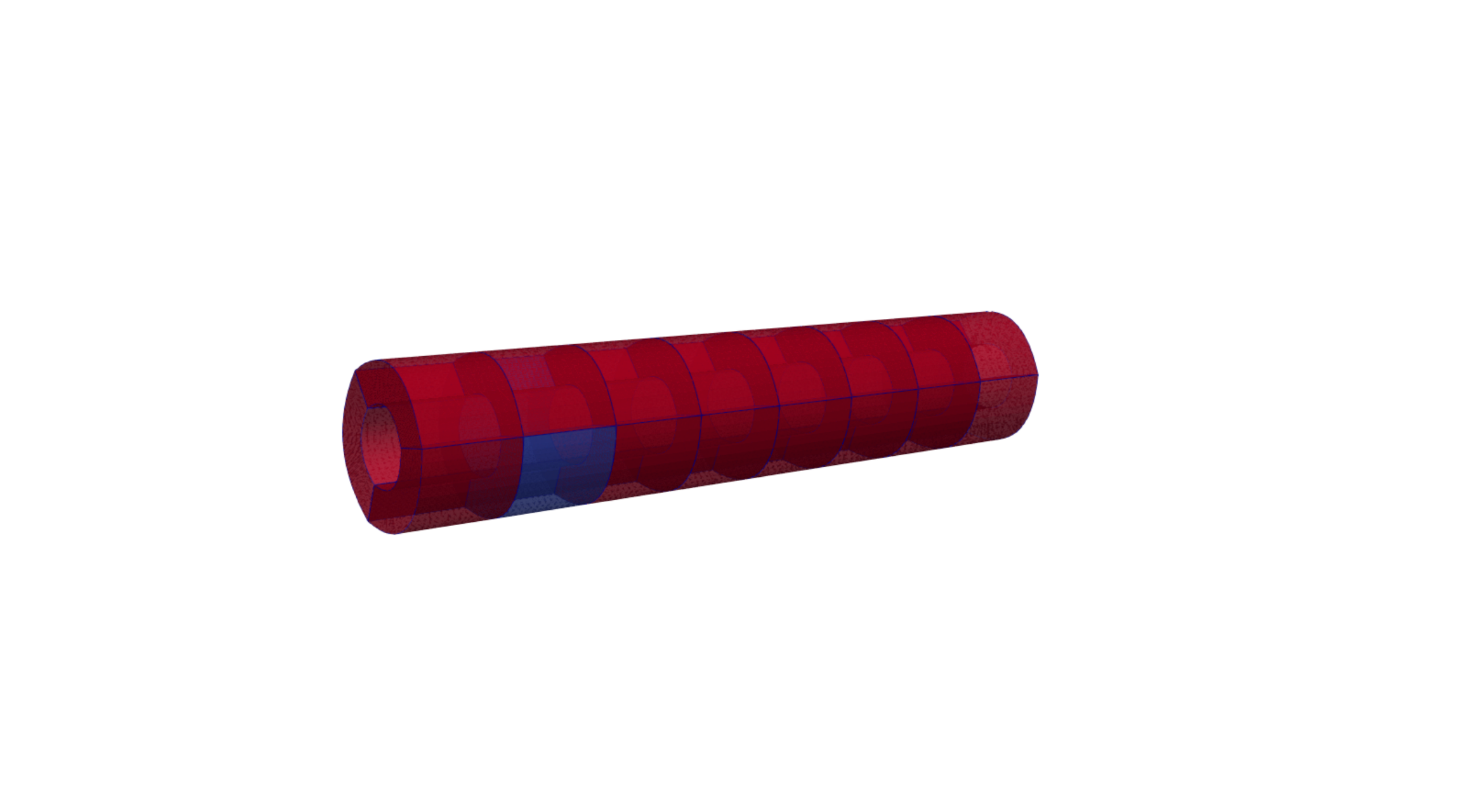}
    \includegraphics[width=0.30\textwidth, trim={300 200 300 150}, clip]{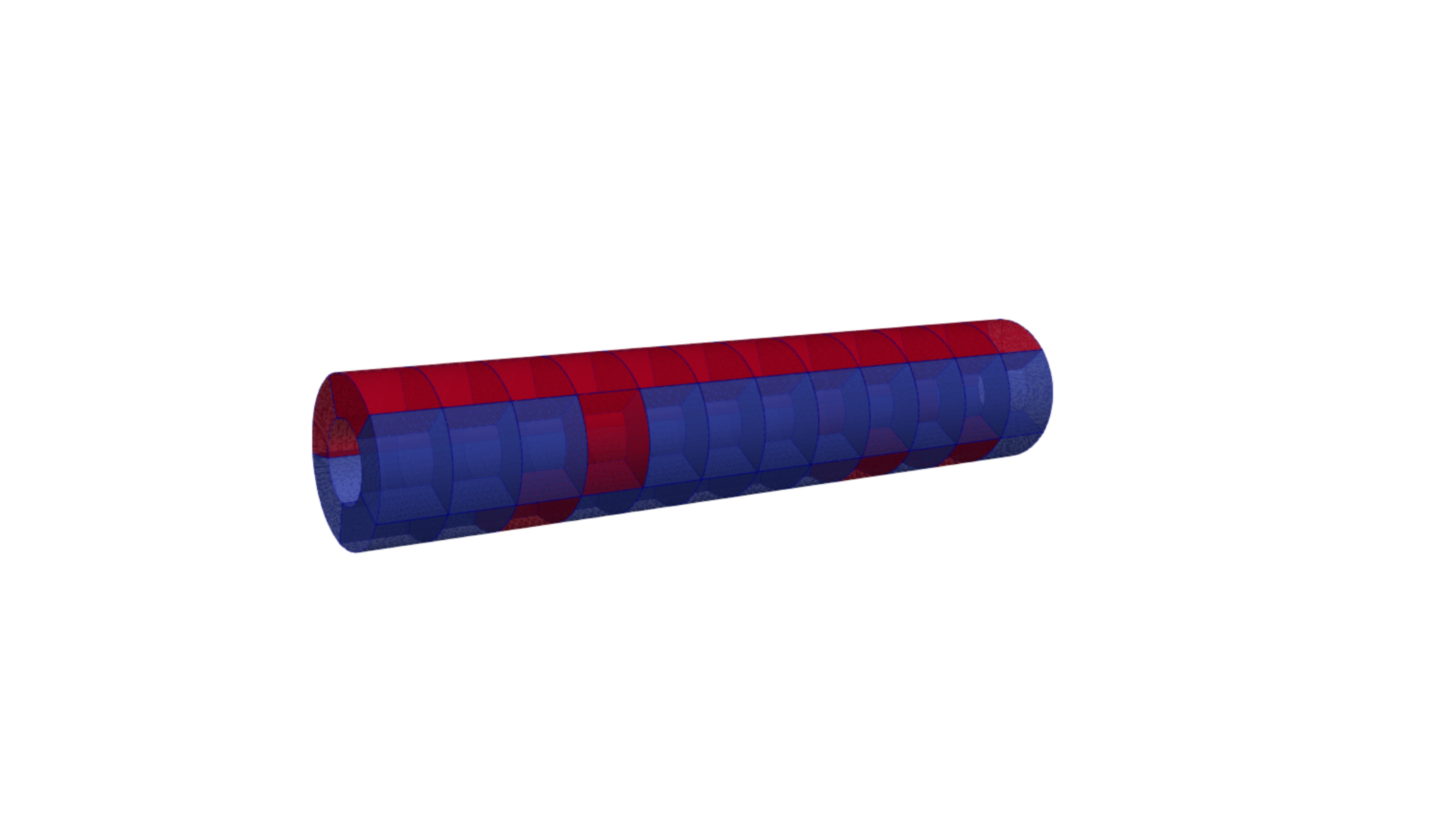}
    \caption{\textbf{CLE.} \textbf{Left:} computational subdomains partitioned in $K=4$ subdomains by \texttt{petsc4py} inside \texttt{deal.II}. \textbf{Center:} test case \textbf{CLE1} repartition of the computational subdomain $k=2$ with the cellwise Grassmannian dimension indicator $I_{\text{G}}$, Definition~\ref{def:ind_grassmannian}: $12\%$ of the cells belong to the \textit{low local Grassmannian dimension} part, represented in blue inside the torus, and the other $88\%$ belong to the \textit{high local Grassmannian dimension} part, represented in red. \textbf{Right:} test case \textbf{CLE2} repartition of the computational subdomain $k=2$ with the cellwise Grassmannian dimension indicator $I_{\text{G}}$ and $P_l=50\%$.}
    \label{fig:subdomains_bar}
\end{figure}

We show two sample solutions for $\mu_1=\mu_2=1000$ in Figure~\ref{fig:elasticity_sample} for \textbf{CLE1} and $\mu_1=\mu_2=1000$, $f_1=1$ and $f_2=-1$ for \textbf{CLE2}, on the left and on the right, respectively. The partitioned and repartitioned subdomains are shown in Figure~\ref{fig:subdomains_bar}. For the first case \textbf{CLE1} we employ a mesh of $24$ cells and $7776$ dofs, for the second \textbf{CLE2} a mesh of $60$ cells and $19440$ dofs.

\begin{figure}[!htpb]
    \centering
    \includegraphics[width=1\textwidth]{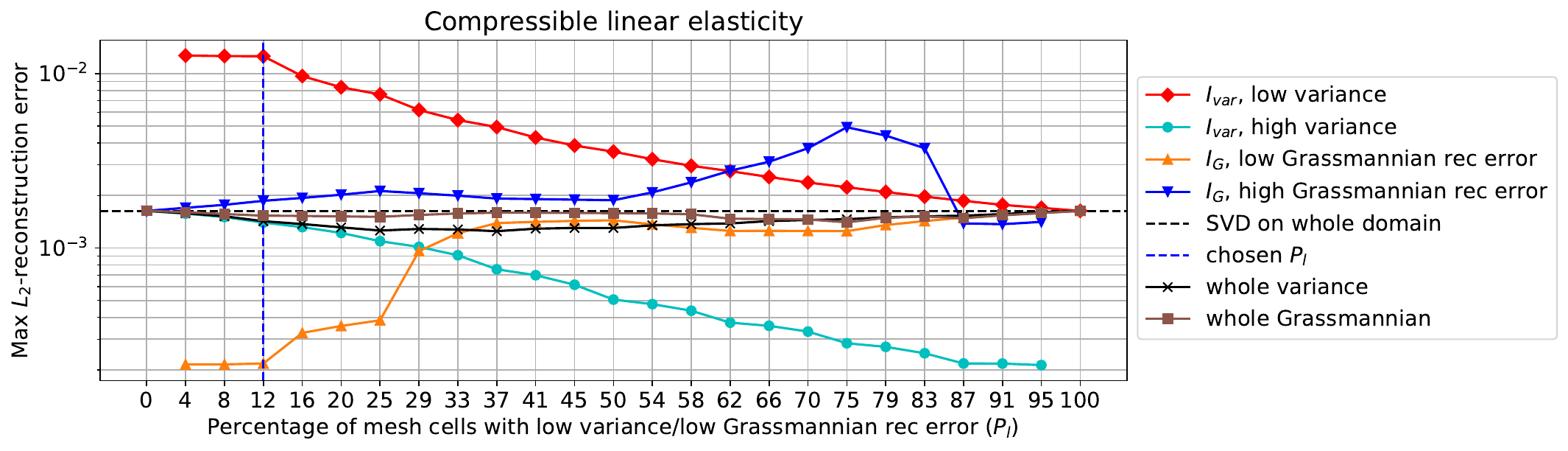}
    \caption{\textbf{CLE1.} Local relative $L^2$-reconstruction errors of the snapshots matrix for elasticity equations restricted to the two subdomains of the repartitioning performed with the indicator $I_{\text{var}}$ (in red and light-blue), Definition~\ref{def:ind_var}, and $I_{\text{G}}$ (in orange and blue), Definition~\ref{def:ind_grassmannian}. The relative $L^2$-reconstruction error attained on the whole domain is shown in black for the indicator $I_{\text{var}}$ and in brown for the indicator $I_{\text{G}}$. The local reduced dimensions used to evaluate the local reconstruction errors is $\NRB_{\Omega_i}=3,\ i=1,2$.}
    \label{fig:repartitioning_elasticity}
\end{figure}

Test case \textbf{CLE1}. This test case presents no region for which the restricted solutions are more or less approximable with a constant field, as would be detected by the variance indicator: as shown in Figure~\ref{fig:repartitioning_elasticity}, the local relative $L^2$-reconstruction error in the region with \textit{low variance}, assigned by $I_{\text{var}}$, deteriorates from the value $2 \cdot 10^{-3}$ of the abscissae $0\%$ and $100\%$ to $1\cdot 10^{-2}$ of the abscissae $4\%$. 
Nonetheless, despite the parametric solutions are not approximabile efficiently with a constant field, they are well represented by a one dimensional linear subspace in the region located by the cellwise Grassmannian dimension indicator $I_G$, for $P_l=12\%$. The associated \textit{low local Grassmannian dimension} region for $P_l=12\%$ is shown in Figure~\ref{fig:subdomains_bar} in blue.

\begin{figure}[!htpb]
    \centering
    \includegraphics[width=1\textwidth]{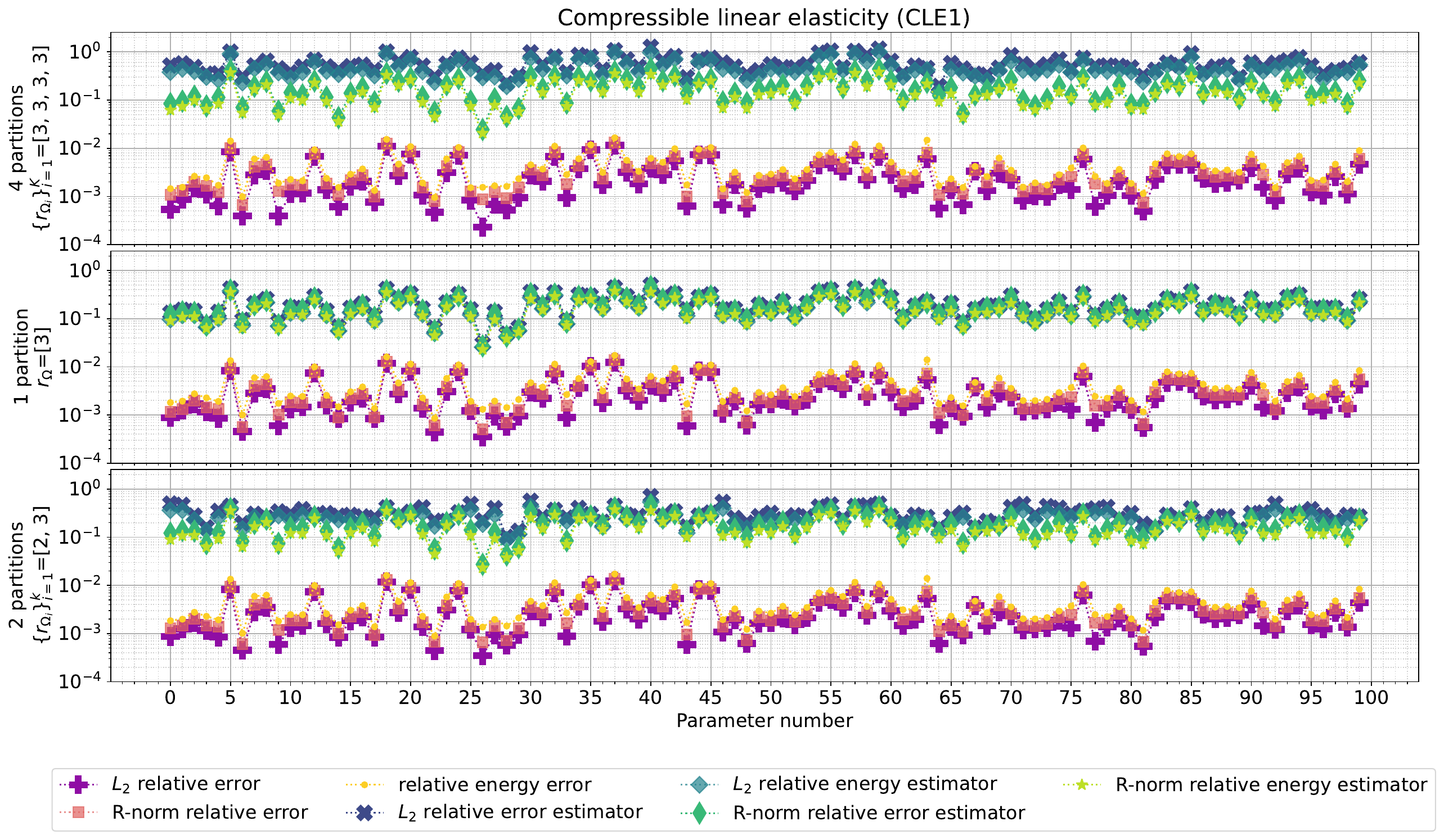}
    \caption{\textbf{CLE1.} Errors and estimators for elasticity equations corresponding to the $n_{\text{train}}=20$ uniformly sampled training snapshots corresponding to the abscissae $0,5,10,\dots,95$, and $n_{\text{test}}=80$ uniformly sampled test snapshots, corresponding to the other abscissae. The reduced dimensions of the ROMs are $\{\NRB_{\Omega_i}\}_{i=1}^K = [3, 3, 3, 3]$ for $K=4$ partitions, $r_{\Omega} = 3$ for $k=1$ partition, and $\{\NRB_{\Omega_i}\}_{i=1}^k = [2, 3]$ for $k=2$ partitions. For the case $k=2$ we employed the cellwise local Grassmannian dimension indicator $I_G$, Definition~\ref{def:ind_grassmannian}, with $P_l=12\%$.}
    \label{fig:a_posteriori_bar}
\end{figure}

Also in this test case, the employment of DD-ROMs is not advisable, since there are little gains in the local relative $L^2$-reconstruction error for the \textit{low local Grassmannian dimensional} region (values around $3\cdot 10^{-3}$, in orange for the abscissa $P_l=12\%$, in Figure~\ref{fig:repartitioning_elasticity}).
The choice of local reduced dimensions $\NRB_{\Omega_1}=2$ and $\NRB_{\Omega_2}=3$ does not affect greatly the errors shown in Figure~\ref{fig:a_posteriori_bar}. Also in this case, we evaluate $n_{\text{train}}=20$ training full-order solutions and $n_{\text{test}}=80$ test full-order solutions, corresponding to a uniform independent sampling from the parametric domain $\mathcal{P}\subset\mathbb{R}^3$. Also for these studies, we have fixed the local dimensions to $\NRB_{\Omega_i}=3,\ i=1,\dots,K$ for $K=4$, $\NRB_{\Omega}=3$ for the whole computational domain and $\NRB_{\Omega_1}=2,\ \NRB_{\Omega_2}=3$ for the repartitioned case with $k=2$.

\begin{figure}[!htpb]
    \centering
    \includegraphics[width=1\textwidth]{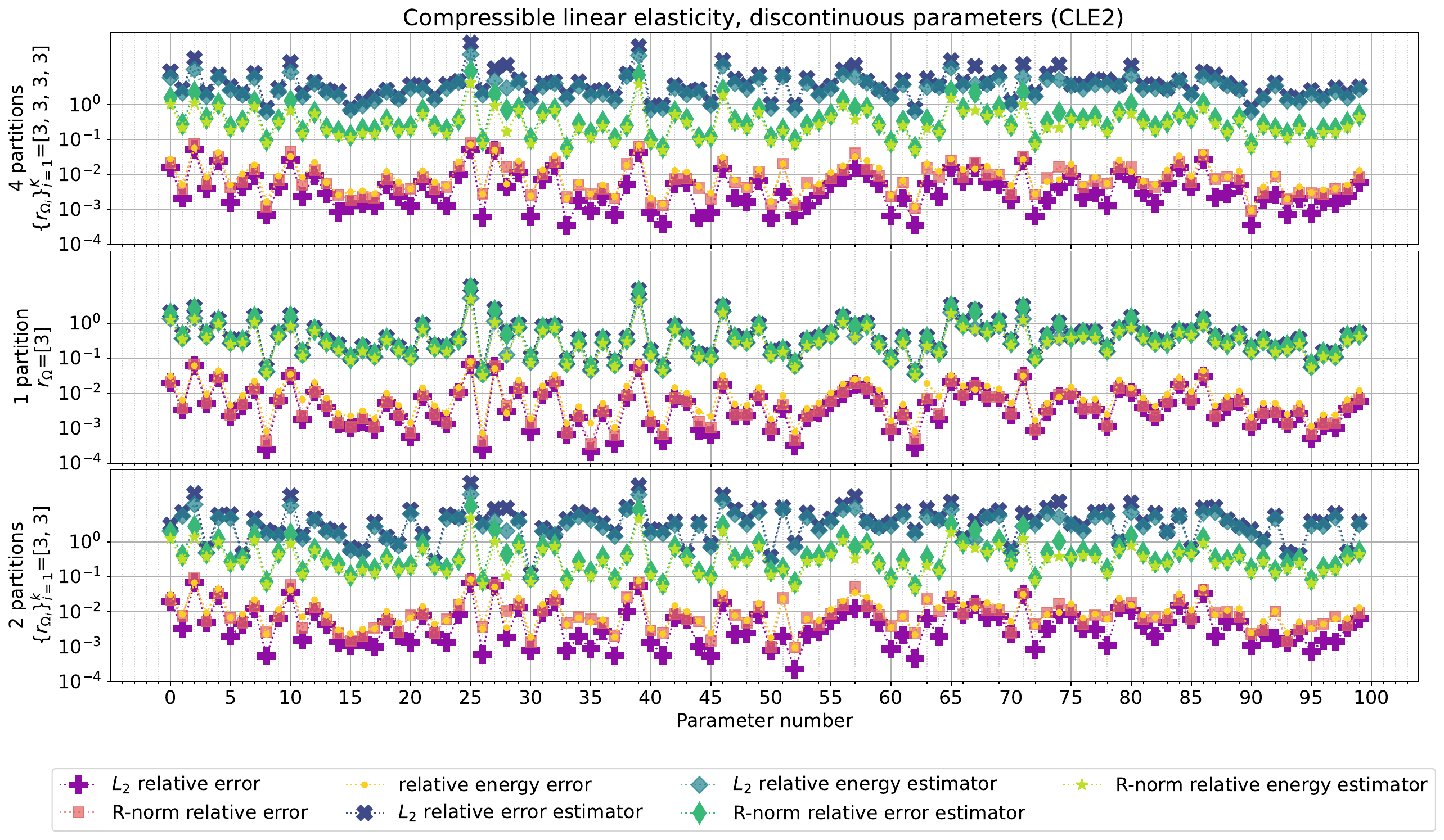}
    \caption{\textbf{CLE2.} Errors and estimators corresponding to the $n_{\text{test}=80}$ uniformly sampled test snapshots corresponding to the abscissae $0,5,10,\dots,95$, and $n_{\text{train}=20}$ uniformly sampled training snapshots, corresponding to the other abscissae. The reduced dimensions of the ROMs are $\{\NRB_{\Omega_i}\}_{i=1}^K = [3, 3, 3, 3]$ for $K=4$ partitions, $\NRB_{\Omega} = 3$ for $k=1$ partition, and $\{\NRB_{\Omega_i}\}_{i=1}^k = [3, 3]$ for $k=2$ partitions. For the case $k=2$ we employed the cellwise variance dimension indicator $I_{\text{var}}$, Definition~\ref{def:ind_var}, with $P_l=50\%$.}
    \label{fig:a_posteriori_bar_disc}
\end{figure}

Test case \textbf{CLE2}. Similarly to the previous case, we evaluate $n_{\text{train}}=20$ training full-order solutions and $n_{\text{test}}=80$ test full-order solutions, corresponding to a uniform independent sampling from the parametric domain $\mathcal{P}_2\subset\mathbb{R}^4$. This time, if we vary the parameters $f_1$ and $f_2$ inside different subdomains $\{z\geq 5\}\cap\Omega$ and $\{z < 5\}\cap\Omega$, we obtain the results shown in Figure~\ref{fig:a_posteriori_bar_disc}. It can be seen that repartitioning $\Omega$ in $k=2$ DD-ROM subdomains with the local Grassmannian indicator $I_G$ and $P_l=50\%$ does not produce more accurate DD-ROMs compared to the case of a single reduced solution manifold for the whole computational domain and for the DD-ROM with $k=4$. In this case, we kept the local dimension of DD-ROM repartitioned case with $k=2$ equal $\NRB_{\Omega_1}=\NRB_{\Omega_2}=3$. For this simple test case, there is not an appreciable improvement for some test parameters in the accuracy for $k=2$ instead of $K=4$ or a classical global linear basis ROM. The reason is that even if the parameters $f_1$ and $f_2$ affect different subdomains of $\Omega$, the solutions on the whole domain are still well correlated. Differently from the previous test case \textbf{MS2} from section~\ref{subsubsec:max}, this is a typical case for which DD-ROMs are not effective, even if the parametrization affects independently two regions of the whole domain $\Omega$. 

In Table~\ref{tab:comp_time_elastic}, we list the computational times and speedups for a simulation with the different methods. For an error analysis with respect to the size of the reduced space, we refer to Appendix~\ref{app:convergence_ROM}.
\begin{table}
	\centering
    \caption{\textbf{CLE.} Average computational times and speedups for ROM and DD-ROM approaches for Maxwell equations. The speedup is computed as the FOM computational time over the ROM one. The FOM runs in parallel with $4$ cores, so ``FOM time'' refers to wallclock time. The first row correspond to test case \textbf{CLE1}, the second to test case \textbf{CLE2}.}\label{tab:comp_time_elastic}	
	\begin{tabular}{|c|c|c|c|c|c|c|c|}\hline
		\multicolumn{2}{|c|}{FOM}&\multicolumn{3}{c|}{ROM}&\multicolumn{3}{c|}{DD-ROM}\\\hline 
		$N_h$ & time & $\NRB$ & time & speedup & $\NRB_i$ & time& speedup \\ \hline\hline
		7776 & 411.510 [ms] & 3 & 80.444 [$\mu$s] & $\sim 5115$ & [3, 3, 3 ,3] & 85.108 [$\mu$s] & $\sim 4835$\\ \hline
        19440 & 2.080 [s] & 3 & 69.992 [$\mu$s] & $\sim 29718$ & [3, 3, 3 ,3] & 94.258 [$\mu$s] & $\sim 22067$\\ \hline
	\end{tabular}
\end{table}

\subsubsection{Scalar concentration advected by an incompressible flow (ADR)}
\label{subsubsec:adr_test}
\begin{figure}[!htpb]
    \centering
    \includegraphics[width=0.7\textwidth]{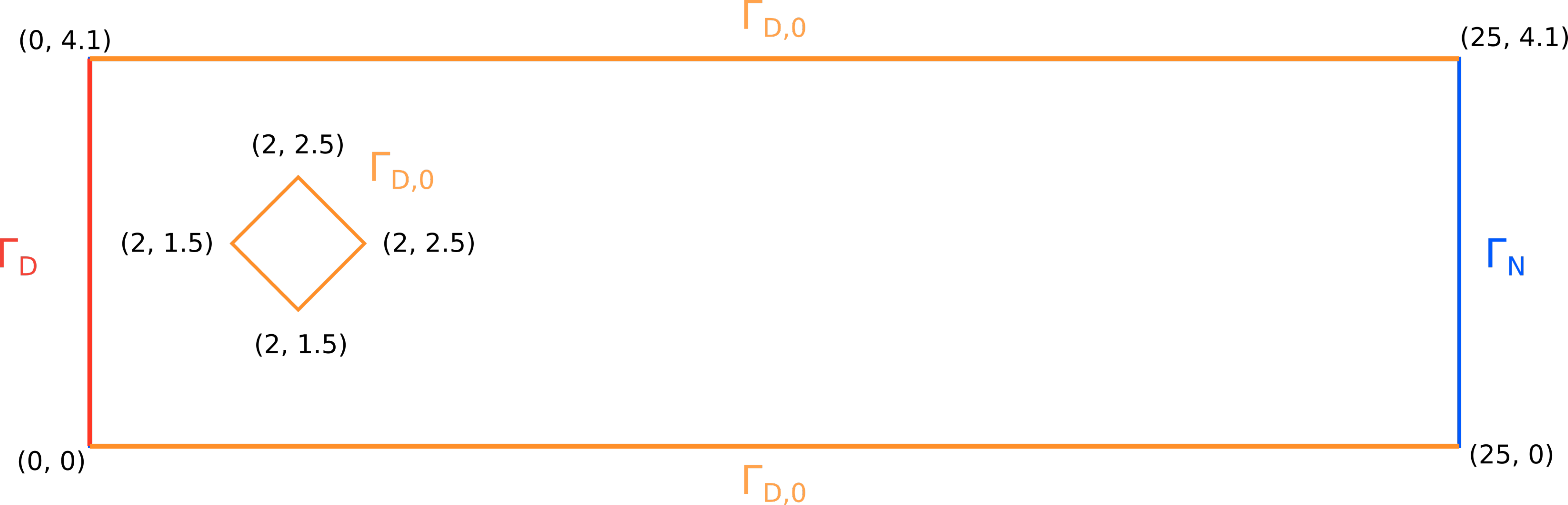}    
    \caption{\textbf{ADR.} Computational domain of the advection diffusion reaction equation FS~\eqref{eq:adr_vv} and the incompressible Navier-Stokes equations~\eqref{eq:INS}. The boundary conditions specified for each system are reported in the text.}  
    \label{fig:adr_domain}
\end{figure}

\begin{figure}[h]
    \centering
    \begin{minipage}{0.49\textwidth}
    	\centering
    	Scalar concentration $u$
    	\includegraphics[width=\textwidth, trim={0 330 0 250}, clip]{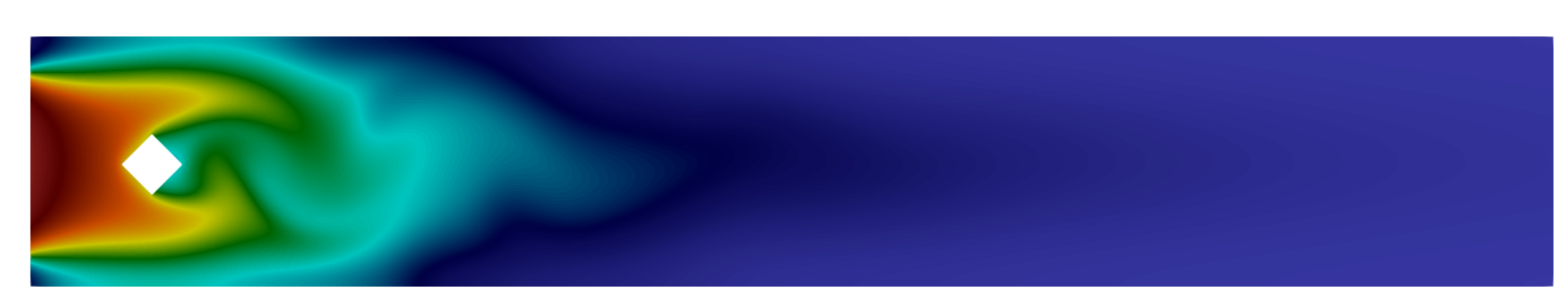}\\
    	\includegraphics[width=\textwidth, trim={360 230 450 590}, clip]{figures/scalar_concentration_u_red.pdf}
    \end{minipage}
    \begin{minipage}{0.49\textwidth}
    \centering
    Magnitude of the advection velocity $\vb$
    \includegraphics[width=\textwidth, trim={0 330 0 250}, clip]{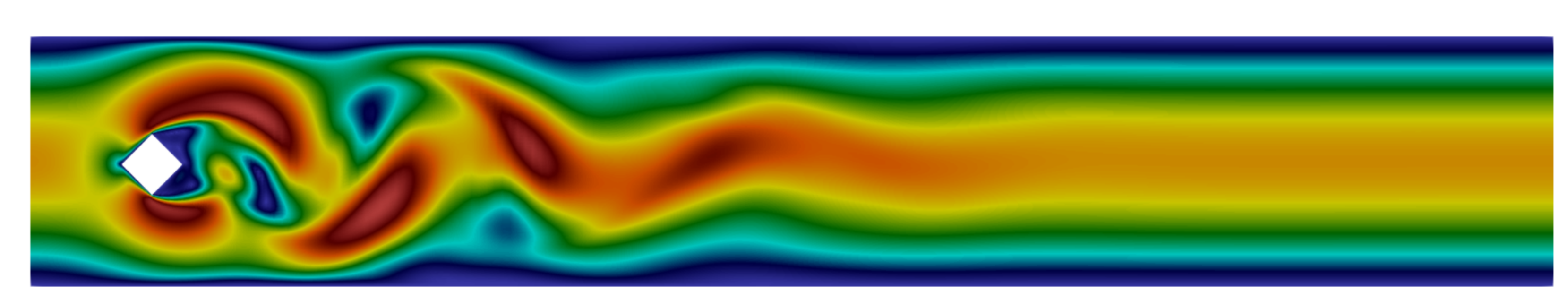}\\
    \includegraphics[width=\textwidth, trim={360 215 450 600}, clip]{figures/advection_velocity_v_red.pdf}
    \end{minipage}
    \caption{\textbf{ADR.} \textbf{Left}: scalar concentration $u$ of the advection diffusion reaction equations~\eqref{eq:adr_vv}, with $\mu_i = 0,\ i=0, 98$ and $\mu_{99}=1$, $\kappa=0.05$. \textbf{Right}: advection velocity employed for the FS~\eqref{eq:adr_vv}, obtained as the velocity $\mathbf{v}$ from the INS~\eqref{eq:INS} at $t=2\text{s}$.}
    \label{fig:sample_adr}
\end{figure}

We consider the parametric semi-linear advection diffusion reaction equation in $d=2$ dimensions, with $m=3$ equations, rewritten in mixed form:
\begin{align}
    \label{eq:adr_vv}
    \begin{cases}
        \kappa^{-1}\sigma + \nabla u = 0,&\quad \text{in}\ \Omega,\\
        \nabla\cdot\sigma +\mathbf{v}\cdot\nabla u + u = f,&\quad \text{in}\ \Omega,\\
        \sigma\cdot\mathbf{n} = 0,&\quad\text{on}\ \Gamma_N\cup\Gamma_{D, 0},\\
        u = \sum_{i=1}^{P}\mu_i \chi_{I_i},&\quad\text{on}\ \Gamma_D,\\
    \end{cases}
\end{align}
where $\kappa=0.05$ is fixed for this study,
\begin{equation}
    \label{eq:parametrization}
    \boldsymbol{\rho}=(\mu_1,\dots,\mu_P)\in\mathcal{P}\subset\mathbb{R}^{P}, \qquad \mathcal{P}=\{\boldsymbol{\rho}\in \{0, 1\}^P\vert \mu_i = 1,\ \mu_{j}=0,\ \forall j\in \{0,\dots,99\}\backslash\{i\}\},
\end{equation}
and $\{\chi_{I_i}\}_{i=0}^{\Npar}$ are the characteristic functions of the symmetric intervals $I_i = {0}\times[-i0.01+1.5, i0.01+2.5]$, with $\Npar=99$. The domain is shown in Figure~\ref{fig:adr_domain}. The advection velocity $\textbf{v}$ is obtained from the following incompressible Navier-Stokes equation at $t=2\text{s}$:
\begin{align}
    \label{eq:INS}
    \begin{cases}
        \partial_t \mathbf{v} + \mathbf{v}\cdot\nabla \mathbf{v} -\nu \Delta \mathbf{v} + \nabla p = \mathbf{0},&\quad \text{in}\ \Omega\\
        \nabla\cdot\mathbf{v} =0,&\quad \text{in}\ \Omega\\
        \mathbf{v}\times\mathbf{n} = 0,\ p = 0,&\quad \text{on}\ \Gamma_N\\
        \mathbf{v} =0,&\quad\text{on}\ \Gamma_{D, 0}\\
        \mathbf{v}(t=0) = \mathbf{v}_b,&\quad\text{on}\ \Gamma_D\\
    \end{cases}
\end{align}
with initial conditions on the boundary $\Gamma_D$, $\mathbf{v}_b=\mathbf{v}(x, y, t=0)=(6y(4.1-y)/4.1^2, 0)\in\mathbb{R}^2$ and $\nu\in\mathbb{R}$ such that the Reynolds number is $Re=100$. The implementation is the one of \texttt{step-35} of the tutorials of the \texttt{deal.II} library~\cite{dealII93}.
Homogeneous Neumann boundary on $\Gamma_N\cup\Gamma_{D, 0}$ and Dirichlet non-homogeneous boundary conditions on $\Gamma_D$ are applied with the boundary operator~\eqref{eq:boundary_adr}.
%:
%\begin{equation}
%	\langle M(\sigma, u), (\tau, v)\rangle_{V^{'}, V} = 
%	\langle \sigma\cdot \n, v\rangle_{-\frac{1}{2}, \frac{1}{2}, \Gamma_D}
%	+ \langle  \tau\cdot \n, u\rangle_{-\frac{1}{2}, \frac{1}{2}, \Gamma_D}
%	-\langle \sigma\cdot \n, v\rangle_{-\frac{1}{2}, \frac{1}{2}, \Gamma_N\cup\Gamma_{D, 0}}
%	- \langle  \tau\cdot \n, u\rangle_{-\frac{1}{2}, \frac{1}{2}, \Gamma_N\cup\Gamma_{D, 0}}.
%\end{equation}
A sample solution is shown in Figure~\ref{fig:sample_adr} for $\mu_i = 0,\ i=0,\dots,98$ and $\mu_{99}=1$, $\kappa=0.05$. We remark that, for the moment, we consider only fixed values of $\kappa=0.05$. For a convergence of ROMs to vanishing viscosity solutions with graph neural networks, see Section~\ref{sec:vv}.

\begin{figure}[h]
    \centering
    \includegraphics[width=0.48\textwidth, trim={0 300 0 300}, clip]{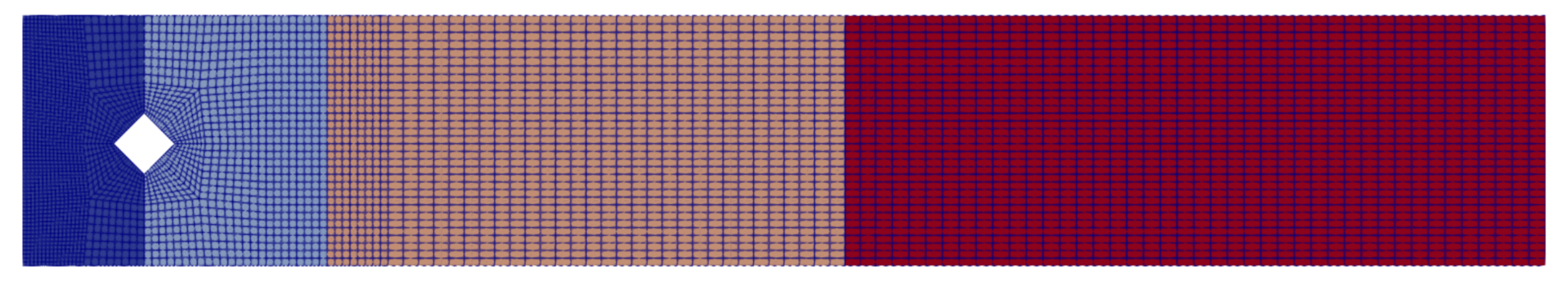}
    \includegraphics[width=0.48\textwidth, trim={0 300 0 300}, clip]{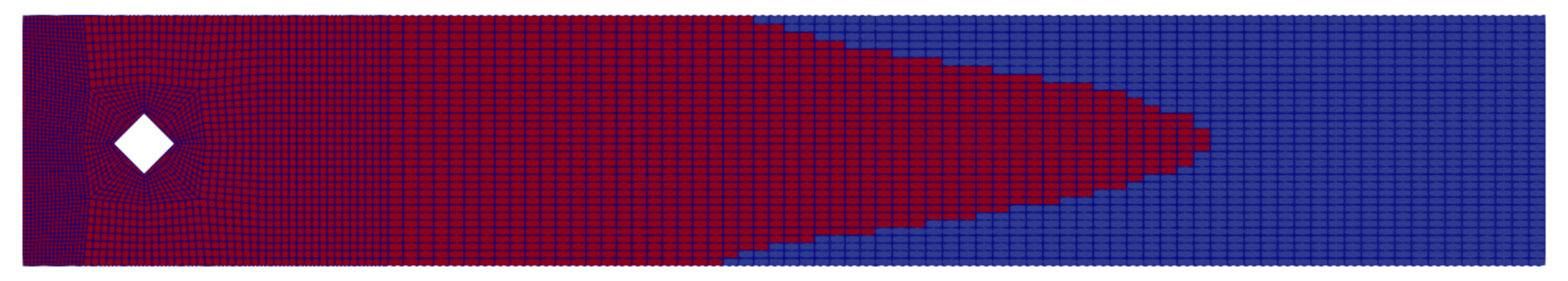}
    \caption{\textbf{ADR.} Domain of advection diffusion reaction equation. \textbf{Left}: computational subdomains partitioned in $K=4$ subdomains by \texttt{petsc4py} inside \texttt{deal.II}. \textbf{Right}: DD-ROM repartition of the computational subdomain $k=2$ with the cellwise variance indicator $I_{\text{var}}$, Definition~\ref{def:ind_var}: $21\%$ of the cells belong to the \textit{low variance} part, represented in blue, and the other $79\%$ belong to the \textit{high variance} part in red.}
    \label{fig:subdomains_adr}
\end{figure}

\begin{figure}[!htpb]
    \centering
    \includegraphics[width=1\textwidth]{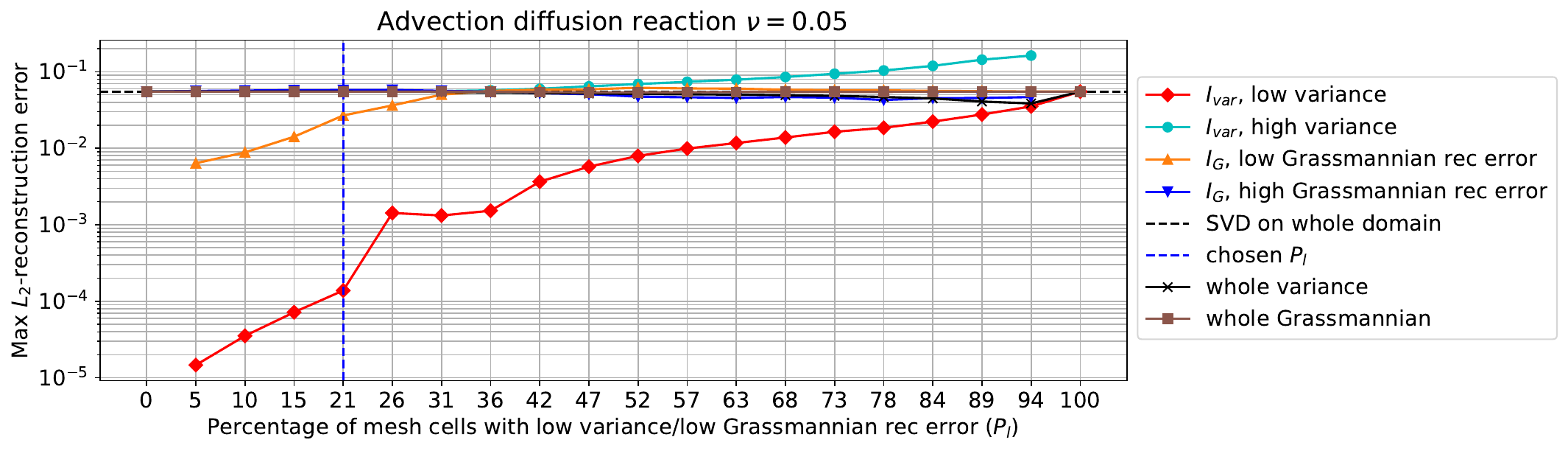}
    \caption{\textbf{ADR.} Local relative $L^2$-reconstruction errors of the snapshots matrix of the advection diffusion reaction equation restricted to the two subdomains of the repartitioning performed with the indicator $I_{\text{var}}$ (in red and light-blue), Definition~\ref{def:ind_var}, and $I_{\text{G}}$ (in orange and blue), Definition~\ref{def:ind_grassmannian}. The relative $L^2$-reconstruction error attained on the whole domain is shown in black for the indicator $I_{\text{var}}$ and in brown for the indicator $I_{\text{G}}$. The local reduced dimensions used to evaluate the local reconstruction errors is $\NRB_{\Omega_i}=3,\ i=1,2$.}
    \label{fig:repartitioning_channel}
\end{figure}

The FOM partitioned and DD-ROM repartitioned subdomains are shown in Figure~\ref{fig:subdomains_adr}. We choose the variance indicator to repartition the computational subdomain in two subset: $21\%$ of the cells for the \textit{low variance} part and $79\%$ for the \textit{high variance} part. With respect to the previous test cases, now it is evident the change in the order of magnitude of the local relative $L^2$-reconstruction error in Figure~\ref{fig:repartitioning_channel}, especially for the cellwise variance indicator $I_{\text{var}}$. We expect that lowering the local reduced dimension of the \textit{low variance} repartitioned region will not affect sensibly the accuracy. 

\begin{figure}[!htpb]
    \centering
    \includegraphics[width=1\textwidth]{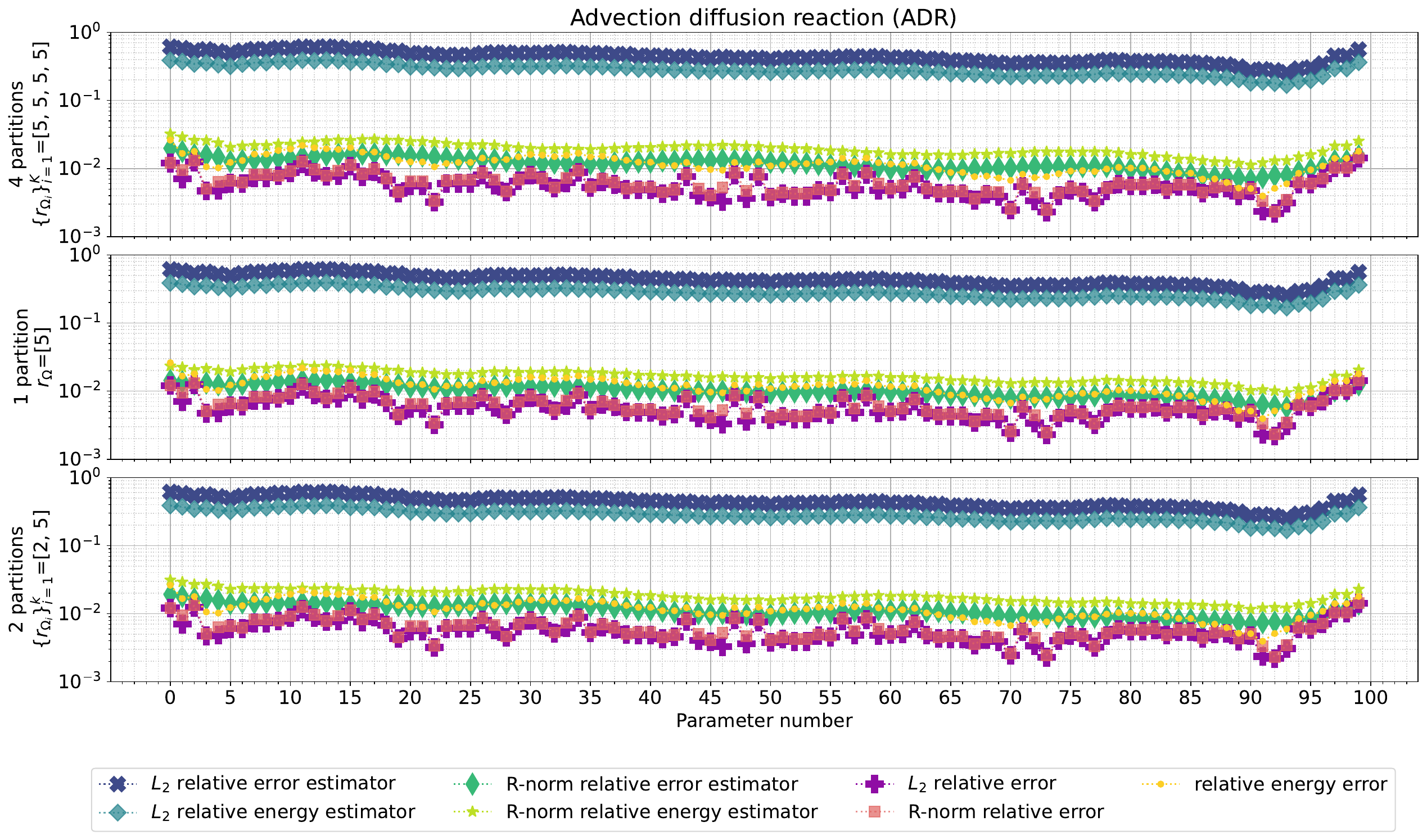}
    \caption{\textbf{ADR.} Errors and estimators for advection diffusion reaction equation corresponding to the $n_{\text{train}=20}$ uniformly sampled train snapshots corresponding to the abscissae $0,5,10,\dots,95$, and $n_{\text{test}=80}$ uniformly sampled test snapshots, corresponding to the other abscissae. The reduced dimensions of the ROMs are $\{\NRB_{\Omega_i}\}_{i=1}^K = [5, 5, 5, 5]$ for $K=4$ partitions, $r_{\Omega} = 5$ for $k=1$ partition, and $\{\NRB_{\Omega_i}\}_{i=1}^k = [2, 5]$ for $k=2$ partitions. For the case $k=2$ we employed the cellwise variance indicator $I_G$, Definition~\ref{def:ind_var}, with $P_l=21\%$.}
    \label{fig:a_posteriori_channel}
\end{figure}

We use for the monodomain approach $\NRB_\Omega=5$ reduced basis as well as $\NRB_{\Omega_i}=5$ for $i=1,\dots,K$ for the FOM partitioned subdomains.
In the DD-ROM approach, we can use even $\NRB_{\Omega_1}=2$ and $\NRB_{\Omega_2}=5$ for the lower and higher variance subdomains, respectively, without affecting the error of the ROM solution, as we see 
in Figure~\ref{fig:a_posteriori_channel}. Indeed, the accuracy in terms of $L^2$ and energy norms is essentially identical for all approaches, even with so little number of basis functions for the DD-ROM one.

Again, we evaluate $n_{\text{train}}=20$ training full-order solutions and $n_{\text{test}}=80$ test full-order solutions, corresponding to the parameter choices $\mu_i=1$ and $\mu_{\bar{i}}=0$, for $i=0,\dots,99$, with fixed viscosity $\kappa=0.05$, where $\bar{i}$ represents all the indices in $\{0,\dots,99\}$ except from $i$. So, the training snapshots correspond to $i=0,5,10,\dots,95$. For these studies we have fixed the local dimensions to $\NRB_{\Omega_i}=5,\ i=1,\dots,K$ for $K=4$, $\NRB_{\Omega}=5$ for the whole computational domain and $\NRB_{\Omega_1}=2,\ r_{\Omega_2}=5$ for the repartitioned case with $k=2$, as mentioned.

In Table~\ref{tab:comp_time_adr}, we list the computational times and speedups for a simulation with the different methods.
\begin{table}
	\centering
    \caption{\textbf{ADR.} Average computational times and speedups for ROM and DD-ROM approaches for Maxwell equations. The speedup is computed as the FOM computational time over the ROM one. The FOM runs in parallel with $4$ cores, so ``FOM time'' refers to wallclock time.}\label{tab:comp_time_adr}	
	\begin{tabular}{|c|c|c|c|c|c|c|c|}\hline
		\multicolumn{2}{|c|}{FOM}&\multicolumn{3}{c|}{ROM}&\multicolumn{3}{c|}{DD-ROM}\\\hline 
		$N_h$ & time & $\NRB$ & time & speedup & $\NRB_i$ & time& speedup \\ \hline\hline
		131328 & 3.243 [s] & 5 & 79.112 [$\mu$s]& $\sim 40992$ & [5, 5, 5, 5] & 59.912 [$\mu$s] & $\sim 54129$\\ \hline
	\end{tabular}
\end{table}

\section{Graph Neural Networks approximating Vanishing Viscosity solutions}
\label{sec:vv}
In this section, we want to highlight how the well-known concept of vanishing viscosity solutions can be related to FS.
In hyperbolic problems, the uniqueness of the weak solution is not guaranteed, already for very simple problems, e.g. inviscid Burgers' equations. 
In order to filter out the physically relevant solution, the concept of vanishing viscosity solution has been introduced, \textit{inter alia} \cite{godlewski1991hyperbolic}, and, consequently, vanishing viscosity methods have been developed, e.g.\ \cite{diperna1983convergence,maday1989analysis}. 

We will consider the topic of vanishing viscosity solutions from the different perspective of model order reduction. 
It is known that slow decaying Kolmogorov n-width solution manifolds result in ineffective linear reduced order models. 
The origin of this problem rests theoretically on the regularity of the parameter to solution map~\cite{cohen2015approximation,cohen2016kolmogorov}, and with less generality on the nature of some PDEs (e.g. advection dominated PDEs, nonlinearities, complex dynamics), 
on the size of the parameter space, and on the smoothness of the parametric initial data or parametric boundary conditions~\cite{arbes2023kolmogorov}, mainly. 
A possible way to obtain more approximable solution manifolds is through regularization or filtering~\cite{xie2018numerical, wells2017evolve}, e.g. adding artificial viscosity. 
Heuristically, the objective is to smoothen out the parametric solutions of the PDEs, for example removing sharp edges, local features, complex patterns, with the aim of designing more efficient ROMs for the filtered solution manifolds.
Then, the linear ROMs will be applied to different levels of \textit{regularization}, still remaining in the regime where they have good approximation properties.
Finally, the original (vanishing viscosity) solutions will be recovered with a regression method from the succession of filtered linear ROMs. 
This is realized without the need to directly reduce with a linear reduced manifold the original solution manifold, thus avoiding the problem of its approximability with a linear subspace and the slow Kolmogorov n-width decay.

In our case, we consider regularization by viscosity levels: the vanishing viscosity solutions $u_{\nu}$ with viscosity $0\leq\nu \ll 1$, will be recovered as the limit $\lim_{i\to\infty}u_{\nu_i} = u_{\nu}$ of a potentially infinite succession of viscosity levels $\{\nu_i\}_{i=0}^{\infty},\ \nu_0>\nu_1>\dots>0$, each associated to its efficient reduced order model. In practice, $\{\nu_i\}_{i=0}^{\infty}\approx\{\nu_i\}_{i=0}^{q}$, where $q$ is the number of additional viscosity ROMs. It is clear the connection with multi-fidelity and super-resolution methods~\cite{fukami2019super, kochkov2021machine}. The rationale of the approach is supported by the proofs of convergence to vanishing viscosity solutions of hyperbolic PDEs under various hypotheses \cite{oleinik1963discontinuous,kruvzkov1970first,diperna1982convergence,goodman1992viscous}.

The framework is general and can be applied in particular to FS. We will achieve this for the advection--diffusion--reaction problem changing the viscosity constant $\mathbb{R}\ni\kappa>0$ in~\eqref{eq:adr_vv}. While this choice is specific for the model we are considering, a more general approach could consist in adding a viscous dissipative term to the generic FS obtaining another FS:
\begin{equation}
    \begin{cases}
        Au = f + \boldsymbol{\kappa}\boldsymbol{\Delta}\mathbf{u},\quad &\text{in}\ \Omega \\
        (\mathcal{D}-\mathcal{M})(u-g) = 0,\quad &\text{on}\ \partial\Omega
    \end{cases}
    \rightarrow 
    \begin{cases}
        \begin{cases}
            \kappa^{-1}\sigma + \nabla u = 0\\
            \nabla\cdot\sigma + Au = f
        \end{cases}
        ,\quad \text{in}\ \Omega \\
        (\mathcal{D}-\mathcal{M})(u-g) = 0,\quad \text{on}\ \partial\Omega,
    \end{cases}
\end{equation}
recalling that the additional degrees of freedom are needed only for the high viscosity ROMs and FOMs (to collect the snapshots) and not the full-order vanishing viscosity solutions. This is only an example of how the procedure could be applied to any FS. In fact, the methodology is not designed specifically for FS.

The overhead of the methodology is related to the evaluation of the snapshots, the assembling of each level of viscosity $\{\nu_i\}_{i=0}^{q}$, and the computational costs of the regression method. We remark that the full matrices of the affine decomposition of each $\{ROM_{\nu_i}\}_{i=0}^{q}$ are the same. This is the price necessary to tackle the realization of reduced order models of parametric PDEs affected by a slow Kolmogorov n-width decay with our approach. 

With respect to standard techniques for nonlinear manifold approximation, the proposed one is more interpretable as a mathematical limit of a succession of solutions to the vanishing viscosity one. Moreover, it has a faster training stage relying on the efficiency of the $\{ROM_{\nu_i}\}_{i=0}^{q}$. To the authors' knowledge, cheap analytical ways to obtain the vanishing viscosity solution from a finite succession of high viscosity ones are not available, so we will rely on data-driven regression methods.

\subsection{Graph neural networks augmented with differential operators}
Generally, machine learning (ML) architectures are employed in surrogate modelling to approximate nonlinear solution manifolds, otherwise linear subspaces are always preferred. The literature is vast on the subject and there are many frameworks that develop surrogate models with ML architectures. They promise to define data-driven reduced order models that infer solutions for new unseen parameters provided that there are enough data to train such architectures. This depends crucially on the choice of the encoding and inductive biases employed to represent the involved datasets: the training computational time and the amount of training data can change drastically.

On this matter, convolutional autoencoders (CNN) are one of the most efficient architectures to approximate nonlinear solution manifolds \cite{lee2020model} for data structured on Cartesian grids, mainly thanks to their shift-equivariance property. For fields on unstructured meshes the natural choice are Graph neural networks (GNNs). Since their employment, GNNs architectures from the ML community have been enriched with physical inductive biases and other tools from numerical analysis. We want to test one of the first implementations and modifications of GNNs~\cite{tencer2021tailored}. We also want to remark that in the literature, there are still very few test cases of ROMs that employ GNNs with more than $\geq 50000$ degrees of freedom. The difficulty arises when the training is performed on large meshes, thus the need for tailored approaches.

The majority of GNNs employed for surrogate modelling are included in autoencoders~\cite{franco2022learning, pichi2023graph} or are directly parametrized to infer the unseen solution with a forward evaluation. These architectures may become heavy, especially for non-academic test cases. One way to tackle the problem of parametric model order reduction of slow Kolmogorov n-width solution manifolds is to employ GNNs only to recover the high-fidelity solution in a multi-fidelity setting, through super-resolution. Since efficient ROMs are employed to obtain the lower levels of fidelity (high viscosity solutions in our case), the solution manifold dimension reduction is performed only at those levels, avoiding the costly and heavy in memory training of autoencoders of GNNs.

We describe the implementation of augmented GNNs as in~\cite{tencer2021tailored}, with the difference that we need to train only a map from a collection of DD-ROMs solutions to the full-order vanishing viscosity solution, and not an autoencoder with pooling and unpooling layers to perform dimension reduction. The GNN we will employ is rather thin with respect to autoencoder GNNs used to perform dimension reduction. Its details are reported in Table~\ref{tab:gnn}.

We represent with
\begin{equation}
    \mathcal{G}=(\mathcal{V}, \mathcal{E},\mathcal{W}),\quad \mathcal{V}\in\mathbb{R}^{n_{\text{nodes}}\times f},\quad \mathcal{E}\in\mathbb{N}^{n_{\text{edges}}\times 2},\quad \mathcal{W}\in\mathbb{R}^{n_{\text{attr}}\times d},
\end{equation}
 a graph with node features $\mathcal{V}$, edges $\mathcal{E}$ and edge attributes $\mathcal{W}$. The number $f$ represents the nodal features dimension. We denote with $\mathbf{e}_{ij}=(i,j)\in\mathbb{N}^2$ the edge between the nodes $\mathbf{n}_i,\mathbf{n}_j\in\mathbb{R}^f$: $\mathbf{e}_{ij}$ corresponds to a row of $\mathcal{E}$, and $\mathbf{n}_i,\mathbf{n}_j$ correspond to the $i$-th and $j$-th rows of $\mathcal{V}$, for $i,j=1,\dots,n_{\text{nodes}}$. Similarly, $\boldsymbol{\omega}_{ij}$ represents the edge attributes of edge $\mathbf{e}_{ij}$. We have $n_{\text{edges}}=n_{\text{attr}}$. For their efficiency, GNNs rely on a message passing scheme composed of propagation and aggregation steps. Supposing that the graph is sparsely connected their implementation is efficient.

When the graph is supported on a mesh, it is natural to consider the generalized support points of finite element spaces as nodes of the graph and the sparsity pattern of the linear system associated to the numerical model as the adjacency matrix of the graph. 
We employ only Lagrangian nodal basis of discontinuous finite element spaces, but the framework can be applied to more general finite element spaces. 
As edge attributes $\boldsymbol{\omega}_{ij}$, we will employ the difference $\boldsymbol{\omega}_{ij}=\mathbf{x}_i-\mathbf{x}_j\in\mathbb{R}^d$ between the corresponding spatial coordinates associated to the nodes $\mathbf{n}_i,\mathbf{n}_j\in\mathbb{R}^f$. The nodes adjacent to node $\mathbf{n}_i$ are represented with the set $\mathcal{N}_{\text{neigh}(i)}$ for all $i={1,\dots,n_{\text{nodes}}}$.

We consider only the two following types of GNN layers: a continuous kernel-based convolutional operator $l_{\text{NNconv}}$~\cite{gilmer2017neural,simonovsky2017dynamic} and the GraphSAGE operator $l_{\text{SAGEconv}}$~\cite{hamilton2017inductive},
\begin{align}
    &\mathcal{V}_{\text{out}} = l_{NNconv}(\mathcal{V}_{\text{inp}}, \mathcal{E},\mathcal{W}) = \mathcal{V}_{\text{inp}}W_3+ \text{Avg}_1(\mathcal{V}_{\text{inp}}, h(\mathcal{W}))+\mathbf{b}_3,\quad h(\mathcal{W})= \text{ReLU}(\mathcal{W}W_1+\mathbf{b_1})W_2 + \mathbf{b_2},\\
    &\mathcal{V}_{\text{out}} = l_{SAGEconv}(\mathcal{V}_{\text{inp}}, \mathcal{E},\mathcal{W}) = \mathcal{V}_{\text{inp}}W_6+ \text{Avg}_2(\text{ReLU}(\mathcal{V}_{\text{inp}}W_4+\mathbf{b_4}))W_5+\mathbf{b}_5,
\end{align}
with weight matrices dimensions,
\begin{align}
    &W_1\in\mathbb{R}^{2\times l},\ W_2\in\mathbb{R}^{l\times (f_{\text{inp}} \times f_{\text{out}}) },\ W_3\in\mathbb{R}^{f_{\text{inp}}\times f_{\text{out}}},\ W_4\in\mathbb{R}^{f_{\text{inp}}\times f_{\text{inp}}},\ W_5,W_6\in\mathbb{R}^{f_{\text{inp}}\times f_{\text{out}}},\\
    &\mathbf{b}_1\in\mathbb{R}^{l},\ \mathbf{b}_2\in\mathbb{R}^{(f_{\text{inp}} \times f_{\text{out}})},\ \mathbf{b}_3,\mathbf{b}_5\in\mathbb{R}^{f_{\text{out}}},\ \mathbf{b}_5\in\mathbb{R}^{f_{\text{inp}}},\\
    &h(\mathcal{W})\in\mathbb{R}^{n_{\text{edges}}\times (f_{\text{inp}} \times f_{\text{out}})},\quad \mathcal{W}=\{W^h_{s}\}_{s=1}^{n_{\text{edges}}},\ W^h_{s}\in\mathbb{R}^{f_{\text{inp}, s}\times f_{\text{out}, s}},\ \forall s=1,\dots,n_{\text{edges}},
\end{align}
with the following average operators used as aggregation operators,
\begin{equation}\label{eq:aggregation_operators}
    (\text{Avg}_1(\mathcal{V}, \{W^h_{s}\}_{s=1}^{n_{\text{edges}}}))_i = \frac{1}{\mathcal{N}_{\text{neigh}(i)}} \sum_{s\in\mathcal{N}_{\text{neigh}}(i)} \!\!W^h_{s} \mathbf{n}_s,\qquad 
    (\text{Avg}_2(\mathcal{V}))_i = \frac{1}{\mathcal{N}_{\text{neigh}(i)}} \sum_{s\in\mathcal{N}_{\text{neigh}}(i)}\mathbf{n}_s,\qquad\forall i = 1,\dots,n_{\text{nodes}},
\end{equation}
where $\mathcal{V}_{\text{inp}}\in\mathbb{R}^{n_{\text{nodes}}\times f_{\text{inp}}},\mathcal{V}_{\text{out}}\in\mathbb{R}^{n_{\text{nodes}}\times f_{\text{out}}}$ are the input and output nodes with feature dimensions $f_{\text{inp}},f_{\text{out}}$. We remark that, differently from graph neural networks with heterogeneous layers, i.e., with changing mesh structure between different layers, in this network the edges $\mathcal{E}$ and edge attributes $\mathcal{W}$ are kept fixed, only the node features change. The feed-forward neural network $h:\mathbb{R}^{n_{\text{edges}}\times d}\rightarrow \mathbb{R}^{f_{\text{inp}, s}\times f_{\text{out}, s}}$ defines a weight matrix $W^h_{s}\in\mathbb{R}^{f_{\text{inp}, s}\times f_{\text{out}, s}}$ for each edge $s=1,\dots,n_{\text{edges}}$. The number $l$ is the hidden layer dimension of $h$.

The aggregation operators are defined from the edges $\mathcal{E}$ that are related to the sparsity pattern of the linear system of the numerical model. So, the aggregation is performed on the stencils of the numerical scheme chosen for every layer of the GNN architecture in Table~\ref{tab:gnn}. Many variants are possible, in particular, we do not employ pooling and unpooling layers to move from different meshes: we always consider the same adapted mesh.

Since our GNNs work on the nodal features, a good strategy is to augment their dimensions as proposed in~\cite{tencer2021tailored}. In fact, in the majority of applications of GNNs for physical models the input features dimensions is the dimension of the physical fields considered and it is usually very small. Considering FS, the fields' dimension is $m$. To augment the input features, we will filter them with some differential operators discretized on the same mesh in which the GNN is supported. We consider the following differential operators
\begin{align}
    \Delta&: V_h(\Omega)\rightarrow V_h(\Omega),\qquad\text{(Laplace operator)}\\\
    \mathbf{v}\cdot\nabla&: V_h(\Omega)\rightarrow V_h(\Omega),\qquad\text{(Advection operator)}\\\
    \nabla_x&: V_h(\Omega)\rightarrow V_h(\Omega),\qquad\text{(Gradient x-component)}\\\
    \nabla_y&: V_h(\Omega)\rightarrow V_h(\Omega),\qquad\text{(Gradient y-component)}\
\end{align}
for a total of four possible feature augmentation operators, where, in our case, $\mathbf{v}$ is the advection velocity from the incompressible Navier-Stokes equations~\eqref{eq:INS}. We employ the representation of the previous differential operators with respect to the polynomial basis of Lagrangian shape functions, so they act on the vectors of nodal evaluations in $\mathbb{R}^{N_h}$. As in~\cite{tencer2021tailored}, we consider three sets of possible augmentations:
\begin{align}
    \mathcal{O}_1 &= \{\mathbb I_{\NDG}, \Delta, \mathbf{v}\cdot\nabla, \nabla_x, \nabla_y\},\\
    \mathcal{O}_2 &= \{\mathbb I_{\NDG}, \nabla_x, \nabla_y\}\\
    \mathcal{O}_3 &= \{\mathbb I_{\NDG}\}
\end{align}
where $\mathbb I_{\NDG}$ is the identity matrix in $\mathbb{R}^{N_h}$, $|\mathcal{O}_1|=5=n_{\text{aug}}$, $|\mathcal{O}_2|=3=n_{\text{aug}}$ and $|\mathcal{O}_3|=1=n_{\text{aug}}$.
We will reconstruct only the scalar concentration $u$ with the GNN, so, in our case, the field dimension is $1$, which is the output dimension. The input dimension depends on the number of high viscosity DD-ROMs employed that we denote with $q$. Given a single parametric instance $\boldsymbol{\rho}\in\mathbb{R}^P$ the associated solutions of $\{\text{D-ROM}_{\kappa_i}\}_{i=1}^{q}$ are $\{\mathbf{u}^{\RB}(\boldsymbol{\rho}_i)\}_{i=1}^{q}$.

We divide the snapshots $\{\mathbf{u}^{\RB}(\boldsymbol{\rho}_i)\}_{i=1}^{n_{\text{train}}+n_{\text{test}}}$ in training $\{\mathbf{u}^{\RB}(\boldsymbol{\rho}_i)\}_{i\in I_{n{\text{train}}}}$ and test snapshots $\{\mathbf{u}^{\RB}(\boldsymbol{\rho}_i)\}_{i\in I_{n{\text{test}}}}$, with $\vert I_{\text{train}}\vert = n_{\text{train}}$ and $\vert I_{\text{test}}\vert = n_{\text{test}}$. We have decided to encode the reconstruction of the vanishing viscosity solution $\mathbf{u}_{q+1}$ learning the difference $\mathbf{u}_{q+1}(\boldsymbol{\rho})-\mathbf{u}^{\RB}_q(\boldsymbol{\rho})-\overline{\mathbf{u}_{q+1}}_{\text{train}}$ with the mesh-supported-augmented GNN ($\text{MSA-GNN}$) described in Table~\ref{tab:gnn}:
\begin{equation}
    \mathbf{u}_{q+1}(\boldsymbol{\rho}) = ReLU\left(\mathbf{u}^{\RB}_q(\boldsymbol{\rho}) + \text{MSA-GNN}\left(\{O_a\{\{\mathbf{u}^{\RB}_i(\boldsymbol{\rho})\}_{i=1}^{q}, \{\mathbf{u}^{\RB}_i(\boldsymbol{\rho})-\mathbf{u}^{\RB}_{i-1}(\boldsymbol{\rho})\}_{i=1}^{q-1}\}\}_{a=2}^{n_{\text{aug}}}\right)+\overline{\mathbf{u}_{q+1}}_{\text{train}}\right),
\end{equation}
where $\overline{\mathbf{u}_{q+1}}_{\text{train}}=\tfrac{1}{n_{\text{train}}}\sum_{i=1}^{n_{\text{train}}}\mathbf{u}_{q+1}(\boldsymbol{\rho_i})$. Learning the difference instead of the solution itself helps in getting more informative features.
The input dimension is therefore $3n_{\text{aug}}=15$ for $\mathcal{O}_1$ and $3n_{\text{aug}}=9$ for $\mathcal{O}_2$.
\begin{table}[htp!]
    \centering
    \caption{Mesh supported augmented GNN}
    \label{tab:gnn}
    % \footnotesize
    \begin{tabular}{ l | c | c | c}
        \hline
        \hline
        Net  & Weights $[f_{\text{inp}}, f_{\text{out}}]$ & Aggregation  & Activation\\
        \hline
        \hline
        Input NNConv & [$3n_{\text{aug}}$, 18] & $\text{Avg}_1$ & ReLU\\
        \hline
        SAGEconv  & [18, 21] & $\text{Avg}_2$ & ReLU\\
        \hline
        SAGEconv  & [21, 24] & $\text{Avg}_2$ & ReLU\\
        \hline
        SAGEconv  & [24, 27] & $\text{Avg}_2$ & ReLU\\
        \hline
        SAGEconv  & [27, 30] & $\text{Avg}_2$ & ReLU\\
        \hline
        Output NNConv  & [30, 1] & $\text{Avg}_1$ & -\\
        \hline
        \hline
    \end{tabular}
    \vspace{1mm}

    \begin{tabular}{ l | c| c| c}
        \hline
        \hline
        NNConvFilters  & First Layer $[2, l]$& Activation & Second Layer $[l, f_{\text{inp}}f_{\text{out}}]$\\
        \hline
        \hline
        Input NNConv  & [2, 12] & ReLU & [12, $3n_{\text{aug}}\cdot 18$]\\
        \hline
        Output NNConv  & [2, 8] & ReLU & [8, 30]\\
        \hline
        \hline
    \end{tabular}
\end{table}

\subsection{Decomposable ROMs approximating vanishing viscosity (VV) solution through GNNs}
In this section, we test the proposed multi--fidelity approach that reconstruct the lowest viscosity level with the GNN. 
We consider the FS~\eqref{eq:adr_vv}, with three levels of viscosity, from highest to lowest: $\kappa_1=0.05$, $\kappa_2=0.01$ and $\kappa_3=0.0005$. 
We want to build a surrogate model that efficiently predicts the parametric solutions of the FS~\eqref{eq:adr_vv} for unseen values of $\boldsymbol{\rho}\in\mathbb{P}$ with fixed viscosity $\kappa_3=0.0005$.
These solutions will be referred to as \textit{vanishing viscosity} solutions. 
The other two viscosity levels are employed to build the $\text{D-ROM}_{\kappa_1}$ and $\text{D-ROM}_{\kappa_2}$ with viscosities $\kappa_1=0.05$ and $\kappa_2=0.01$, respectively. The parametrization affects the inflow boundary condition and is the same as the one described in \cref{subsubsec:adr_test}, see equation~\eqref{eq:parametrization}. We also employ the same number of training $20$ and test $80$ parameters.

\begin{figure}[!htpb]
    \centering    
    \includegraphics[width=0.49\textwidth, trim={0 200 0 200}, clip]{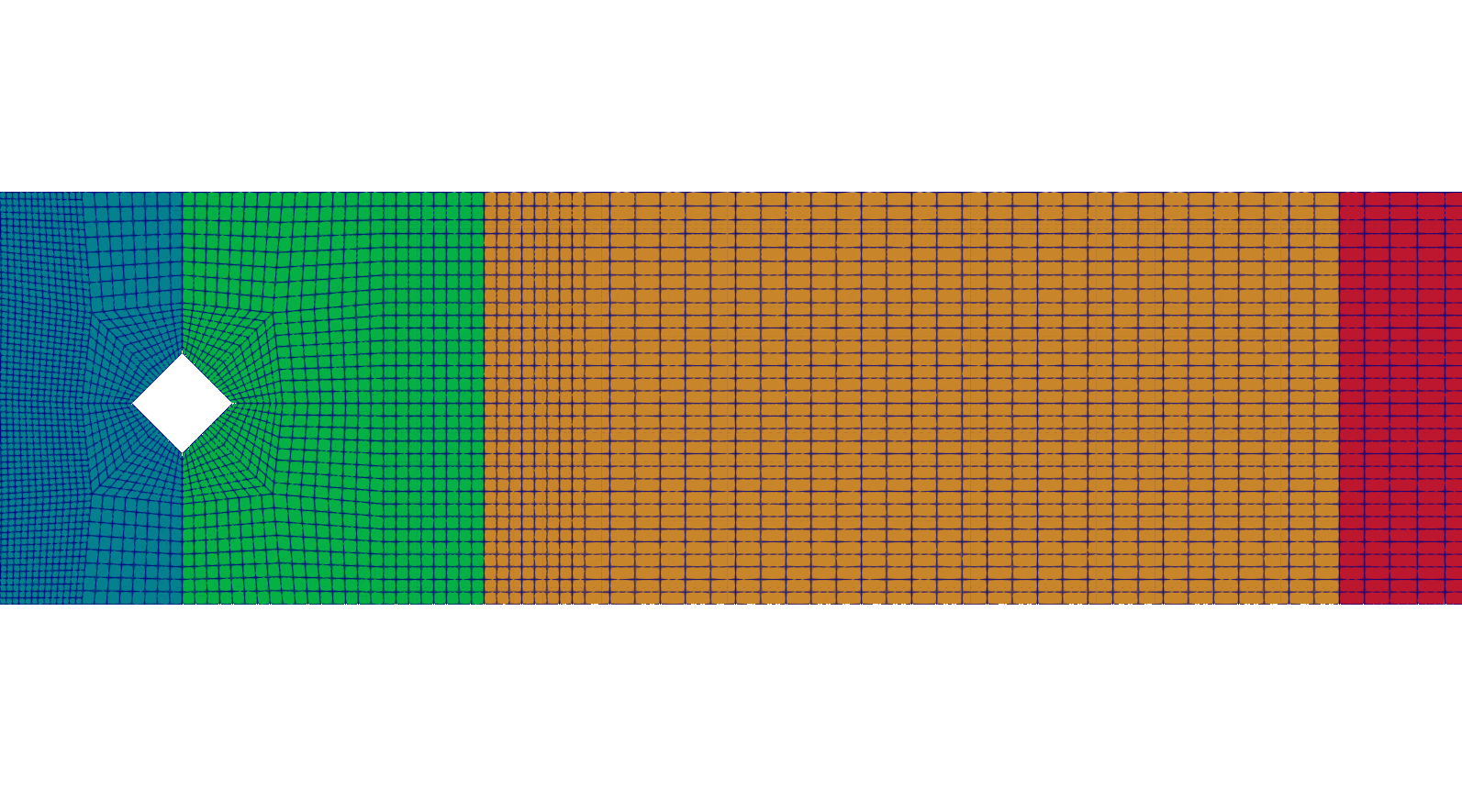}
    \includegraphics[width=0.49\textwidth, trim={0 200 0 200}, clip]{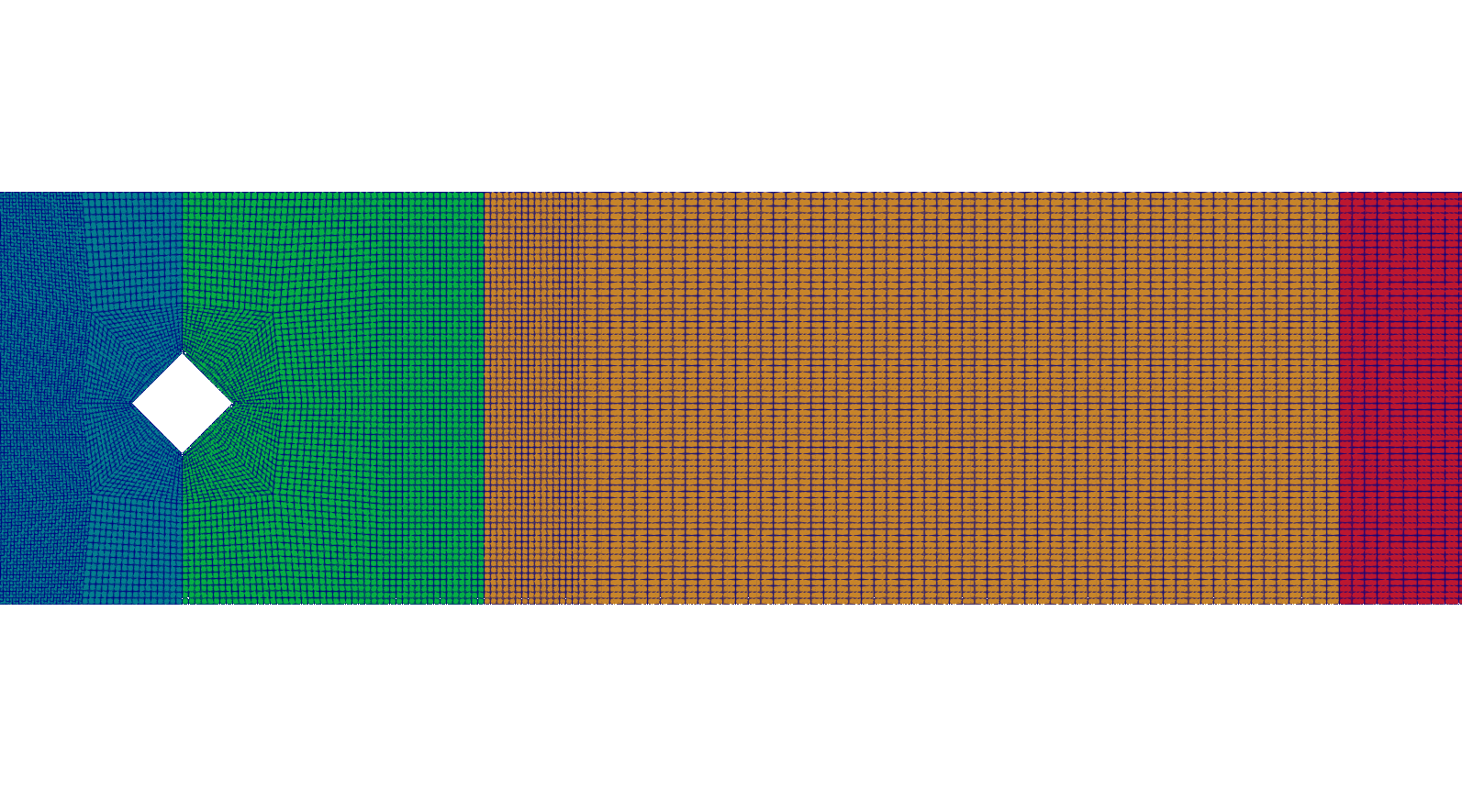}
    \caption{\textbf{VV.} Left part of the computational domain partitioned in $4$ for distributed parallelism: coarse mesh (\textbf{left}), fine mesh (\textbf{right}). The solution with viscosity $\kappa\in\{\textbf{0.05}, \textbf{0.01}\}$ are evaluated on the coarse mesh with ${4868}$ cells and ${43776}$ dofs, those with $\kappa=\textbf{0.0005}$ on the finer with ${19456}$ cells and ${175104}$ dofs.}
    \label{fig:domain_vv}
\end{figure}

The DD-ROMs provided for $\kappa_1=0.05$ and $\kappa_2=0.01$ can be efficiently designed with reduced dimensions $\{r_{\Omega_i}\}_{i=1}^{K}[5, 5, 5 ,5]$. 
To further reduce the cost, we employ an even coarser mesh for $\text{ROM}_{\kappa_1}$ and $\text{ROM}_{\kappa_2}$ and a finer mesh for the vanishing viscosity solutions. 
The former is represented on the left of Figure~\ref{fig:domain_vv}, the latter on the right.
The degrees of freedom related to the coarse mesh are $43776$, while the ones on the fine one are $175104$.

For the training of the GNN we use the open source software library PyTorch Geometric~\cite{Fey/Lenssen/2019}. The employment of efficient samplers that partition the graphs on which the training set is supported is crucial to lower the otherwise heavy memory burden~\cite{hamilton2017inductive}. We preferred samplers that partition the mesh with METIS~\cite{Karypis_1998} as it is often employed in this context. We decided to train the GNN with early stopping at $50$ epochs as our focus is also in the reduction of the training time of the NN architectures used for model order reduction. It corresponds on average to less than $60$ minutes of training time. The batch size is $100$ and we clustered the whole domain in $100000$ subgraphs in order to fit the batches in our limited GPU memory. Each augmentation strategy and additional fidelity level, do not affect the whole training time as they only increase the dimension of the input features from a minimum of $1$ ($1$ fidelity, no augmentation) to a maximum of $15$ (all augmentations $\mathcal{O}_3$, $2$ fidelities). As optimizer we use ADAM~\cite{kingma2014adam} stochastic optimizer. Every architecture is trained on a single GPU NVIDIA Quadro RTX 4000.

\begin{figure}[!htpb]
	\newcommand{\testNumb}{0}
    \centering    
    \begin{minipage}{0.325\textwidth}\centering
    	FOM $u$, $\kappa = 0.05$\\
    	\includegraphics[width=\textwidth,trim={10 314 10 300}, clip]{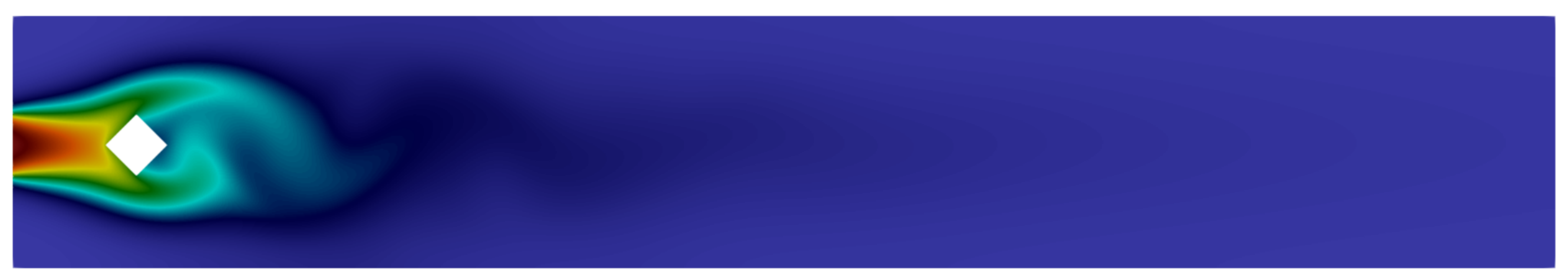}\\
    	\includegraphics[width=\textwidth,trim={400 200 400 606}, clip]{figures/gnn/05_fom_\testNumb_red.pdf}    	
    \end{minipage}
    \begin{minipage}{0.325\textwidth}\centering
    	ROM $u$, $\kappa = 0.05$\\
    	\includegraphics[width=\textwidth,trim={10 314 10 300}, clip]{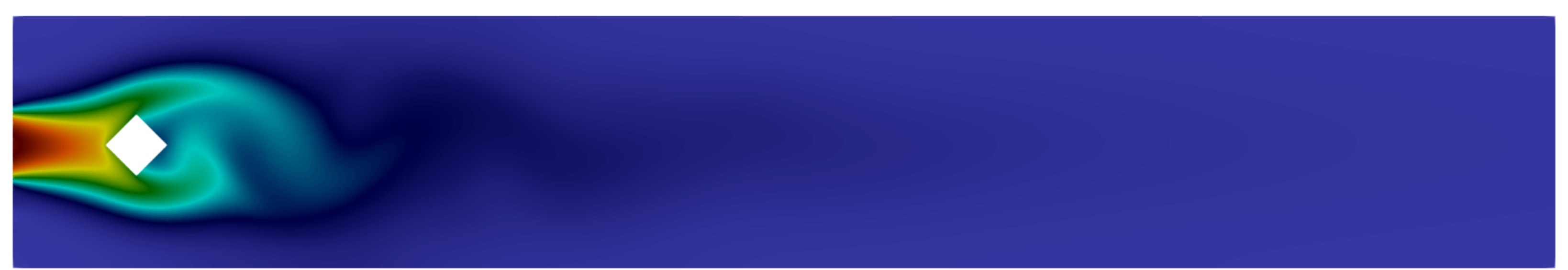}\\
    	\includegraphics[width=\textwidth,trim={400 200 400 606}, clip]{figures/gnn/05_rom_\testNumb_red.pdf}    	
    \end{minipage}
    \begin{minipage}{0.325\textwidth}\centering
    Difference FOM--ROM $u$ with $\kappa = 0.05$\\
    \includegraphics[width=\textwidth,trim={10 314 10 300}, clip]{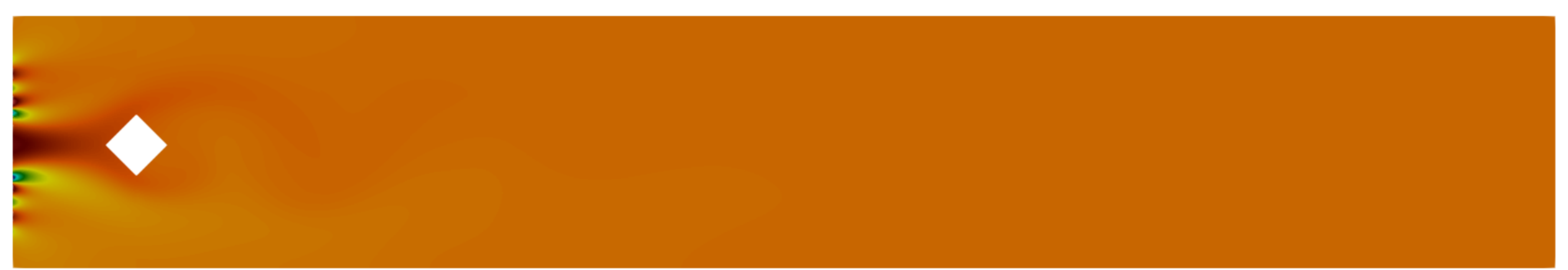}\\
    \includegraphics[width=\textwidth,trim={400 200 400 606}, clip]{figures/gnn/05_diff_\testNumb_red.pdf}    	
    \end{minipage}\\ \vspace*{2mm}
    \begin{minipage}{0.325\textwidth}\centering
    	FOM $u$, $\kappa = 0.01$\\
    	\includegraphics[width=\textwidth,trim={10 314 10 300}, clip]{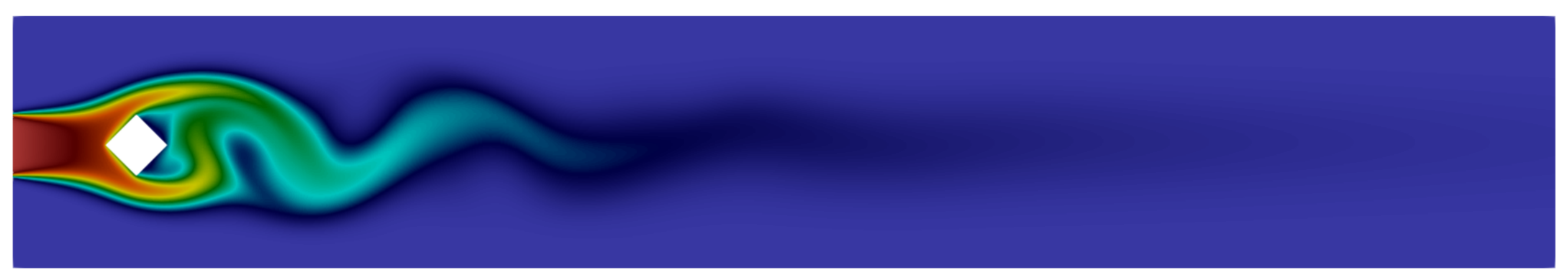}\\
    	\includegraphics[width=\textwidth,trim={400 200 400 606}, clip]{figures/gnn/01_fom_\testNumb_red.pdf}    	
    \end{minipage}
    \begin{minipage}{0.325\textwidth}\centering
    	ROM $u$,  $\kappa = 0.01$\\
    	\includegraphics[width=\textwidth,trim={10 314 10 300}, clip]{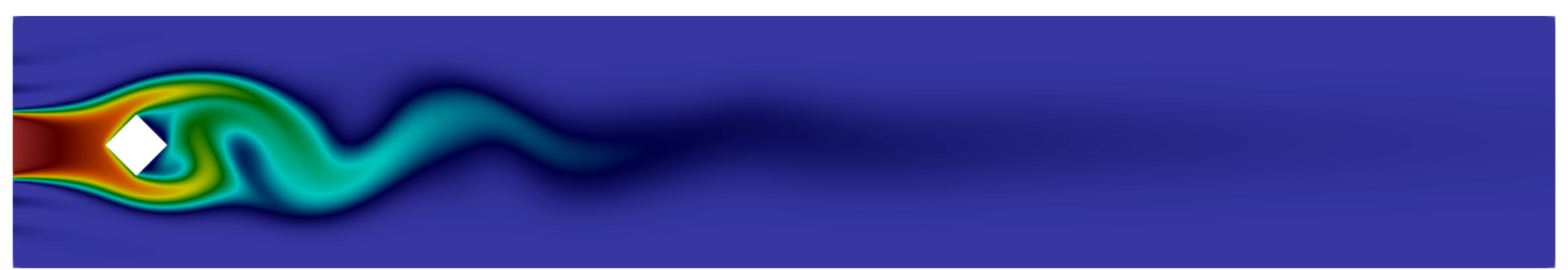}\\
    	\includegraphics[width=\textwidth,trim={400 200 400 606}, clip]{figures/gnn/01_rom_\testNumb_red.pdf}    	
    \end{minipage}
    \begin{minipage}{0.325\textwidth}\centering
    	Difference FOM--ROM $u$, $\kappa = 0.01$\\
    	\includegraphics[width=\textwidth,trim={10 314 10 300}, clip]{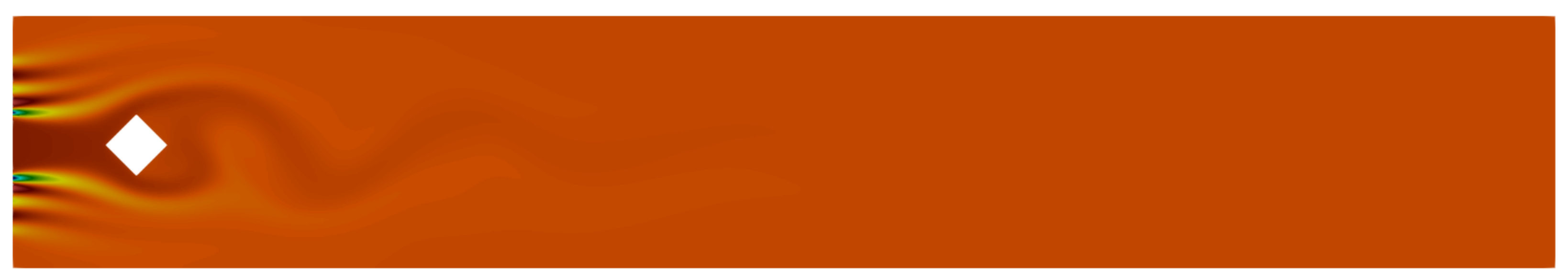}\\
    	\includegraphics[width=\textwidth,trim={400 200 400 606}, clip]{figures/gnn/01_diff_\testNumb_red.pdf}    	
    \end{minipage}\\ \vspace*{2mm}    
    \begin{minipage}{0.325\textwidth}\centering
    	FOM $u$, $\kappa = 0.0005$\\
    	\includegraphics[width=\textwidth,trim={10 314 10 300}, clip]{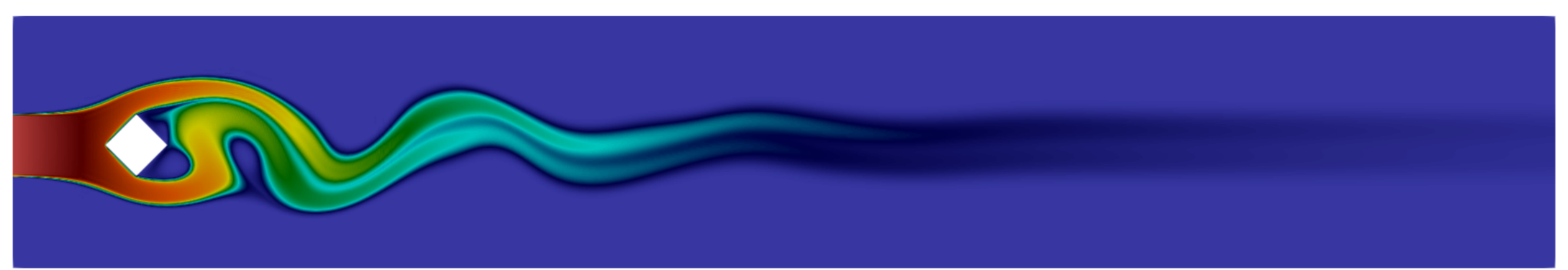}\\
    	\includegraphics[width=\textwidth,trim={400 200 400 606}, clip]{figures/gnn/0005_fom_\testNumb_red.pdf}    	
    \end{minipage}
    \begin{minipage}{0.325\textwidth}\centering
    	ROM $u$, $\kappa = 0.0005$\\
    	\includegraphics[width=\textwidth,trim={10 314 10 300}, clip]{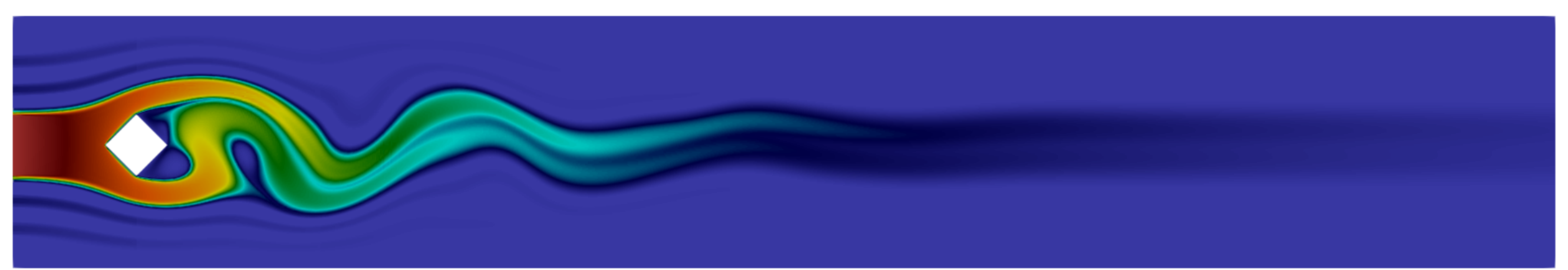}\\
    	\includegraphics[width=\textwidth,trim={400 200 400 606}, clip]{figures/gnn/0005_rom_\testNumb_red.pdf}    	
    \end{minipage}
    \begin{minipage}{0.325\textwidth}\centering
    	Difference FOM--ROM $u$, $\kappa = 0.0005$\\
    	\includegraphics[width=\textwidth,trim={10 314 10 300}, clip]{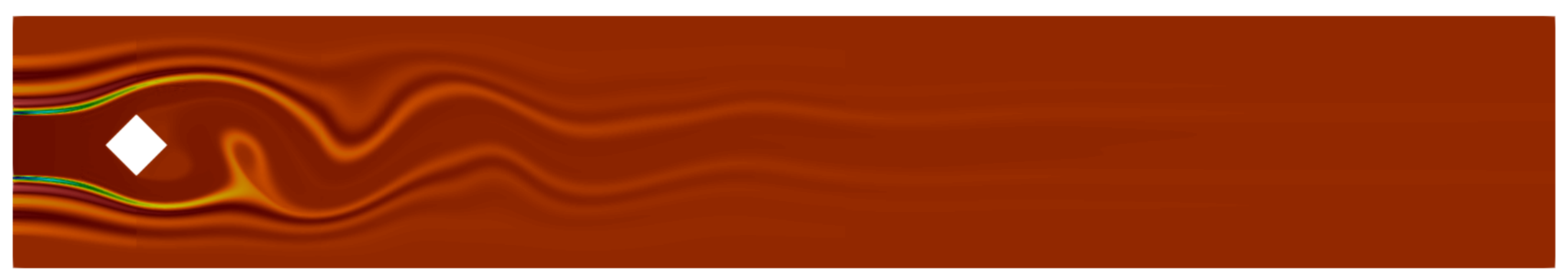}\\
    	\includegraphics[width=\textwidth,trim={400 200 400 606}, clip]{figures/gnn/0005_diff_\testNumb_red.pdf}    	
    \end{minipage}\\ \vspace*{2mm}
        \begin{minipage}{0.325\textwidth}\centering
    	FOM $u$, $\kappa = 0.0005$\\
    	\includegraphics[width=\textwidth,trim={10 314 10 300}, clip]{figures/gnn/0005_fom_\testNumb_red.pdf}\\
    	\includegraphics[width=\textwidth,trim={400 200 400 600}, clip]{figures/gnn/0005_fom_\testNumb_red.pdf}    	
    \end{minipage}
    \begin{minipage}{0.325\textwidth}\centering
    	GNN $u$, $\kappa = 0.0005$\\
    	\includegraphics[width=\textwidth,trim={10 314 10 300}, clip]{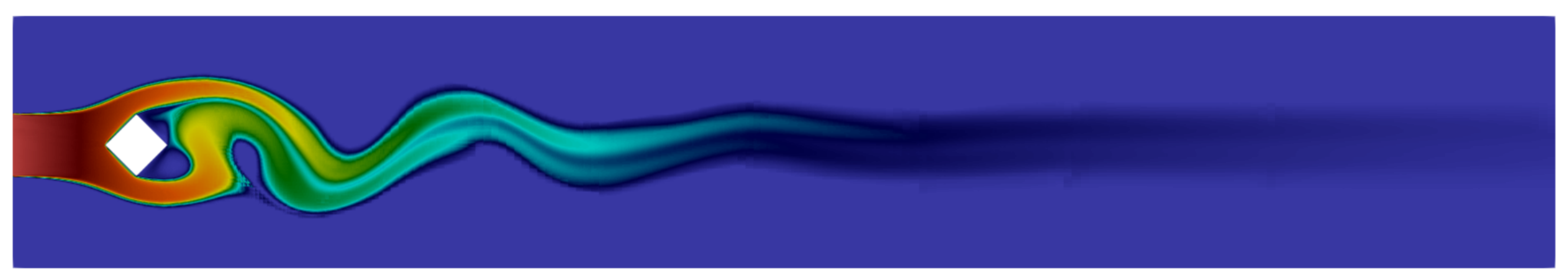}\\
    	\includegraphics[width=\textwidth,trim={400 200 400 600}, clip]{figures/gnn/0005_gnn_\testNumb_red.pdf}    	
    \end{minipage}
    \begin{minipage}{0.325\textwidth}\centering
    	Difference FOM--GNN $u$, $\kappa = 0.0005$\\
    	\includegraphics[width=\textwidth,trim={10 314 10 300}, clip]{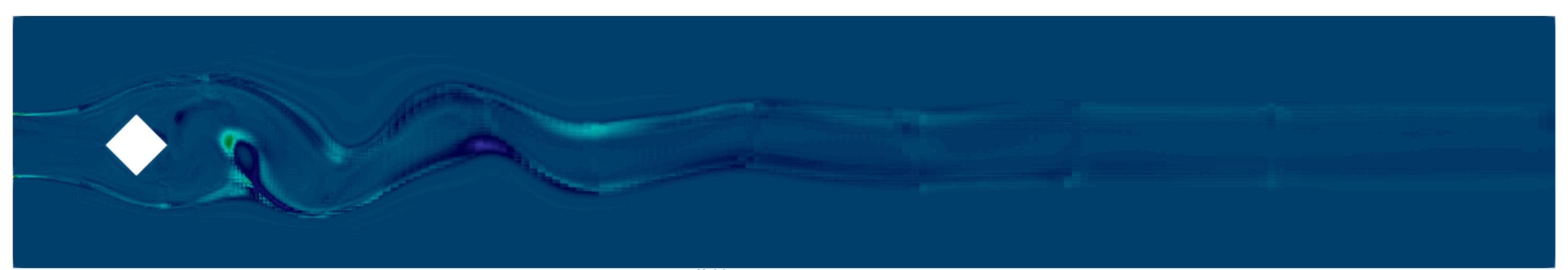}\\
    	\includegraphics[width=\textwidth,trim={400 200 400 600}, clip]{figures/gnn/0005_diff_gnn_\testNumb_red.pdf}    	
    \end{minipage}
    \caption{\textbf{VV.} Scalar concentration advected by incompressible flow for ${i}=\testNumb$. Comparison of ROM approach at different viscosity levels $\kappa\in \{0.05, 0.01, 0.0005\}$ and GNN for $\kappa=0.0005$. FOMs on the left, reduced solution at the center and error on the right.}
    \label{fig:show_vv_0}
\end{figure}

\begin{figure}[!htpb]
	\newcommand{\testNumb}{50}
	\centering    
	\begin{minipage}{0.325\textwidth}\centering
		FOM $u$, $\kappa = 0.05$\\
		\includegraphics[width=\textwidth,trim={10 314 10 300}, clip]{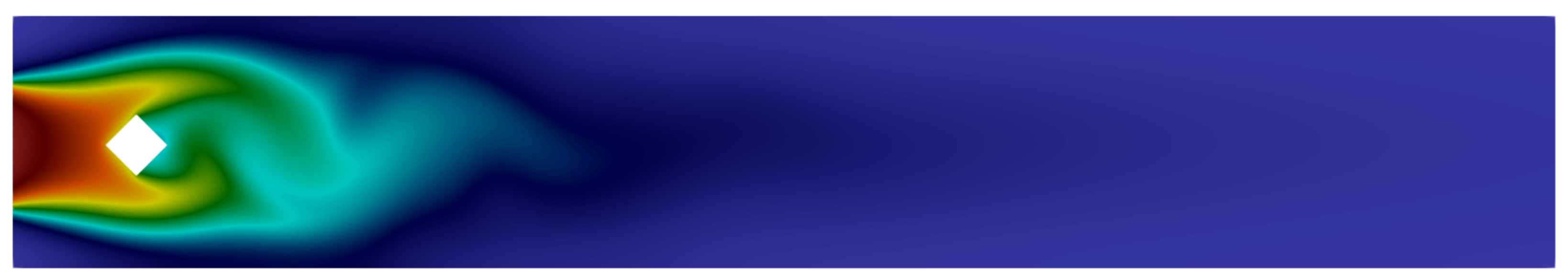}\\
		\includegraphics[width=\textwidth,trim={400 200 400 606}, clip]{figures/gnn/05_fom_\testNumb_red.pdf}    	
	\end{minipage}
	\begin{minipage}{0.325\textwidth}\centering
		ROM $u$, $\kappa = 0.05$\\
		\includegraphics[width=\textwidth,trim={10 314 10 300}, clip]{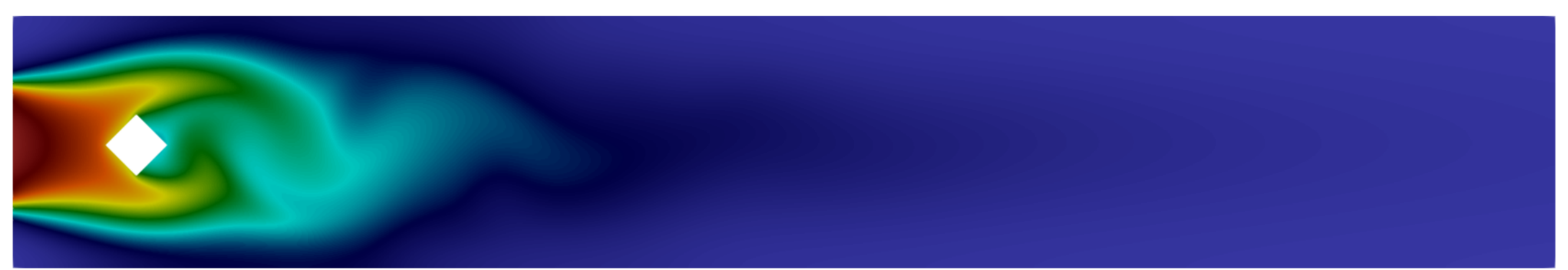}\\
		\includegraphics[width=\textwidth,trim={400 200 400 606}, clip]{figures/gnn/05_rom_\testNumb_red.pdf}    	
	\end{minipage}
	\begin{minipage}{0.325\textwidth}\centering
		Difference FOM--ROM $u$ with $\kappa = 0.05$\\
		\includegraphics[width=\textwidth,trim={10 314 10 300}, clip]{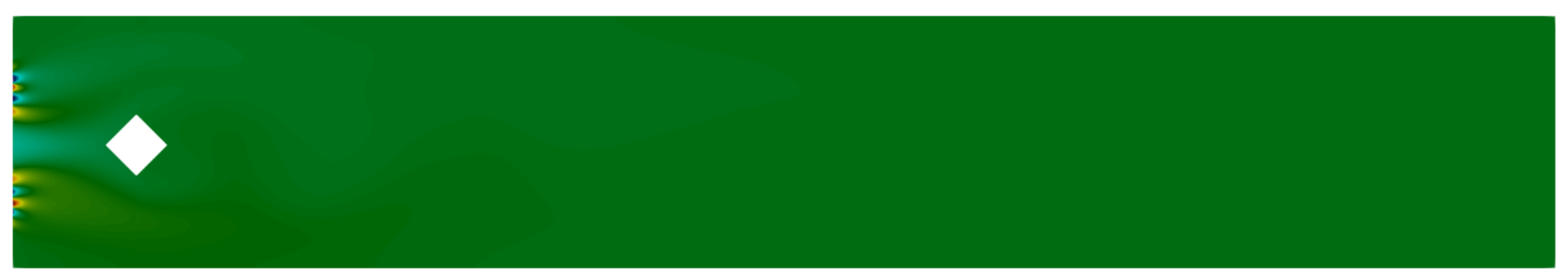}\\
		\includegraphics[width=\textwidth,trim={400 200 400 606}, clip]{figures/gnn/05_diff_\testNumb_red.pdf}    	
	\end{minipage}\\ \vspace*{2mm}
	\begin{minipage}{0.325\textwidth}\centering
		FOM $u$, $\kappa = 0.01$\\
		\includegraphics[width=\textwidth,trim={10 314 10 300}, clip]{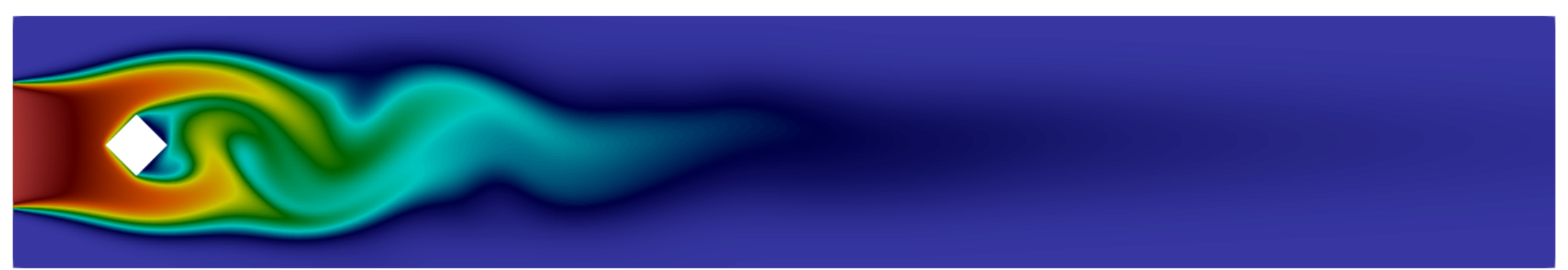}\\
		\includegraphics[width=\textwidth,trim={400 200 400 606}, clip]{figures/gnn/01_fom_\testNumb_red.pdf}    	
	\end{minipage}
	\begin{minipage}{0.325\textwidth}\centering
		ROM $u$,  $\kappa = 0.01$\\
		\includegraphics[width=\textwidth,trim={10 314 10 300}, clip]{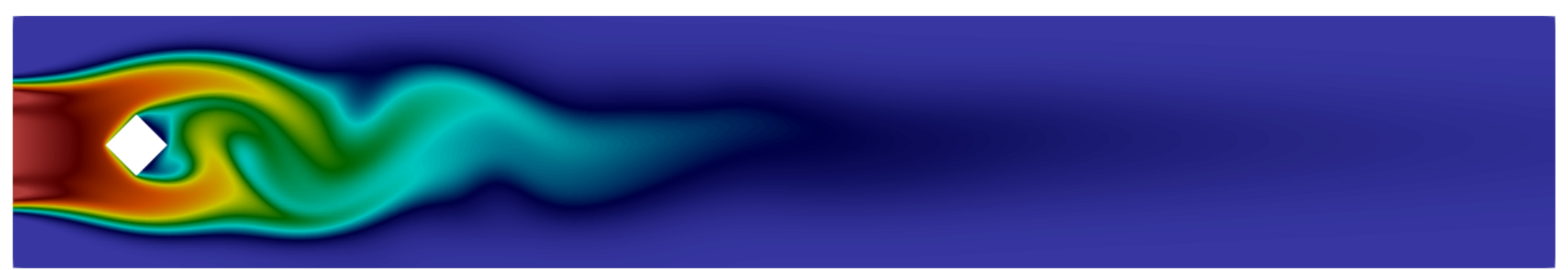}\\
		\includegraphics[width=\textwidth,trim={400 200 400 606}, clip]{figures/gnn/01_rom_\testNumb_red.pdf}    	
	\end{minipage}
	\begin{minipage}{0.325\textwidth}\centering
		Difference FOM--ROM $u$, $\kappa = 0.01$\\
		\includegraphics[width=\textwidth,trim={10 314 10 300}, clip]{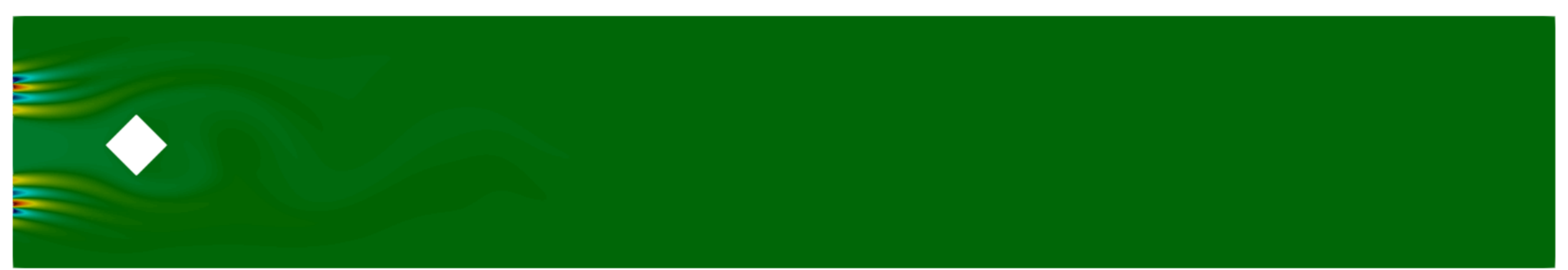}\\
		\includegraphics[width=\textwidth,trim={400 200 400 606}, clip]{figures/gnn/01_diff_\testNumb_red.pdf}    	
	\end{minipage}\\ \vspace*{2mm}    
	\begin{minipage}{0.325\textwidth}\centering
		FOM $u$, $\kappa = 0.0005$\\
		\includegraphics[width=\textwidth,trim={10 314 10 300}, clip]{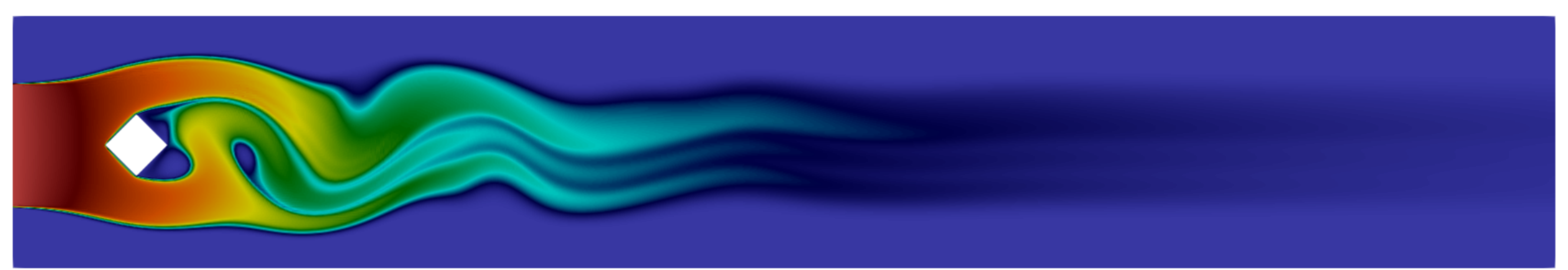}\\
		\includegraphics[width=\textwidth,trim={400 200 400 606}, clip]{figures/gnn/0005_fom_\testNumb_red.pdf}    	
	\end{minipage}
	\begin{minipage}{0.325\textwidth}\centering
		ROM $u$, $\kappa = 0.0005$\\
		\includegraphics[width=\textwidth,trim={10 314 10 300}, clip]{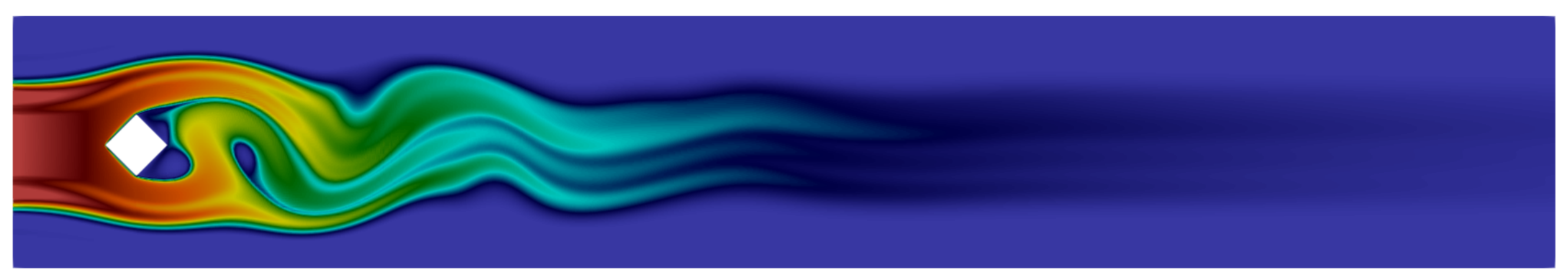}\\
		\includegraphics[width=\textwidth,trim={400 200 400 606}, clip]{figures/gnn/0005_rom_\testNumb_red.pdf}    	
	\end{minipage}
	\begin{minipage}{0.325\textwidth}\centering
		Difference FOM--ROM $u$, $\kappa = 0.0005$\\
		\includegraphics[width=\textwidth,trim={10 314 10 300}, clip]{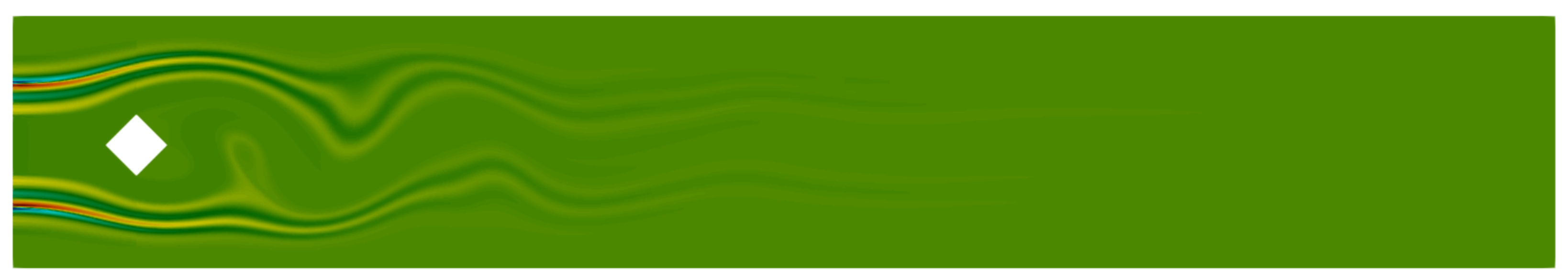}\\
		\includegraphics[width=\textwidth,trim={400 200 400 606}, clip]{figures/gnn/0005_diff_\testNumb_red.pdf}    	
	\end{minipage}\\ \vspace*{2mm}
	\begin{minipage}{0.325\textwidth}\centering
		FOM $u$, $\kappa = 0.0005$\\
		\includegraphics[width=\textwidth,trim={10 314 10 300}, clip]{figures/gnn/0005_fom_\testNumb_red.pdf}\\
		\includegraphics[width=\textwidth,trim={400 200 400 600}, clip]{figures/gnn/0005_fom_\testNumb_red.pdf}    	
	\end{minipage}
	\begin{minipage}{0.325\textwidth}\centering
		GNN $u$, $\kappa = 0.0005$\\
		\includegraphics[width=\textwidth,trim={10 314 10 300}, clip]{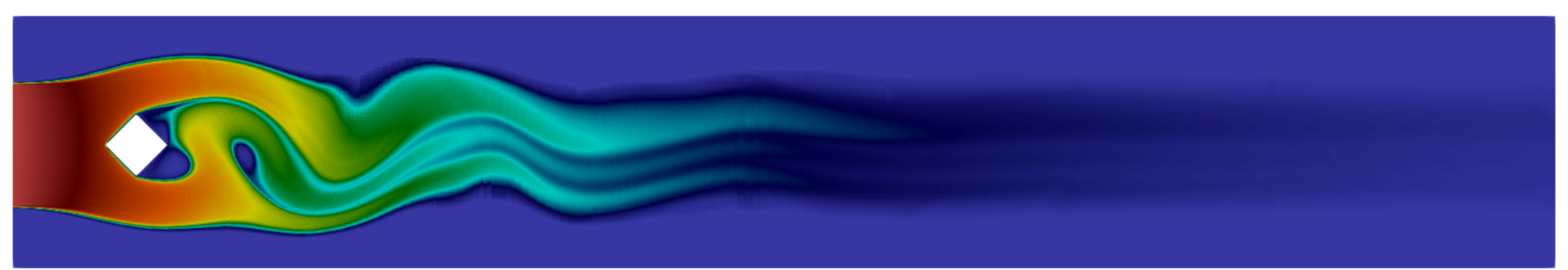}\\
		\includegraphics[width=\textwidth,trim={400 200 400 600}, clip]{figures/gnn/0005_gnn_\testNumb_red.pdf}    	
	\end{minipage}
	\begin{minipage}{0.325\textwidth}\centering
		Difference FOM--GNN $u$, $\kappa = 0.0005$\\
		\includegraphics[width=\textwidth,trim={10 314 10 300}, clip]{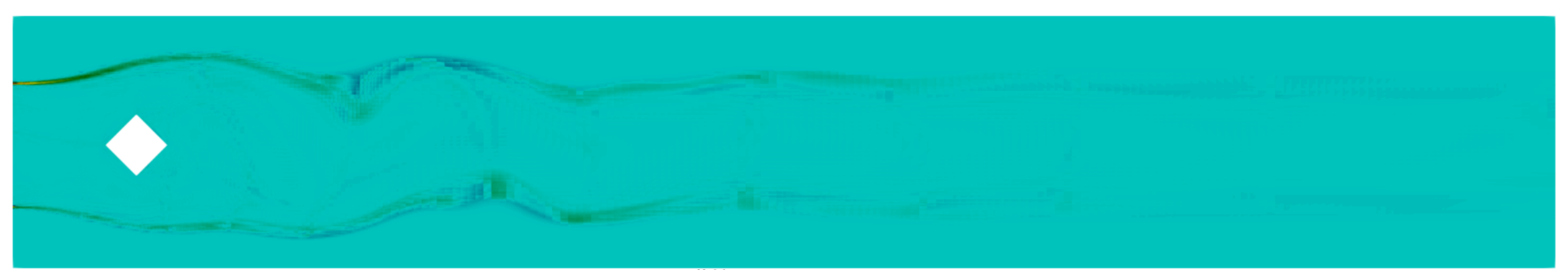}\\
		\includegraphics[width=\textwidth,trim={400 200 400 600}, clip]{figures/gnn/0005_diff_gnn_\testNumb_red.pdf}    	
	\end{minipage}
	\caption{\textbf{VV.} Scalar concentration advected by incompressible flow for $i=\testNumb$. Comparison of ROM approach at different viscosity levels $\kappa\in \{0.05, 0.01, 0.0005\}$ and GNN for $\kappa=0.0005$. FOMs on the left, reduced solution at the center and error on the right.}
	\label{fig:show_vv_50}
\end{figure}

\begin{figure}[!htpb]
	\newcommand{\testNumb}{99}
	\centering    
	\begin{minipage}{0.325\textwidth}\centering
		FOM $u$, $\kappa = 0.05$\\
		\includegraphics[width=\textwidth,trim={10 314 10 300}, clip]{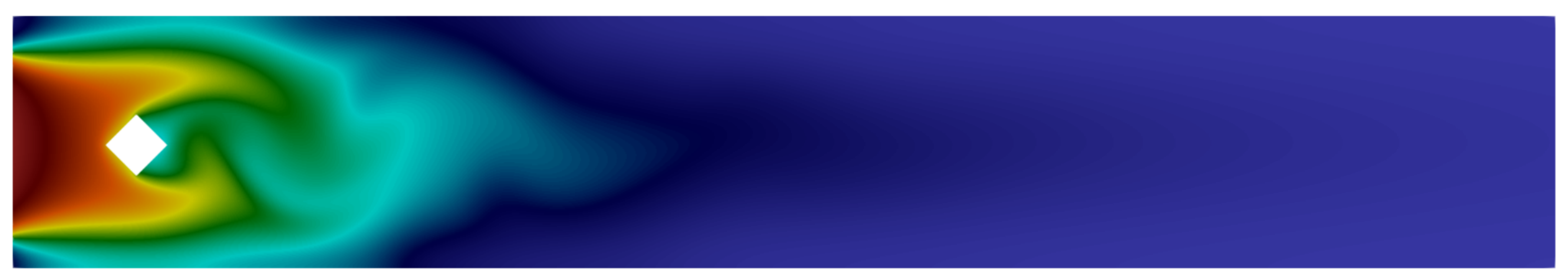}\\
		\includegraphics[width=\textwidth,trim={400 200 400 606}, clip]{figures/gnn/05_fom_\testNumb_red.pdf}    	
	\end{minipage}
	\begin{minipage}{0.325\textwidth}\centering
		ROM $u$, $\kappa = 0.05$\\
		\includegraphics[width=\textwidth,trim={10 314 10 300}, clip]{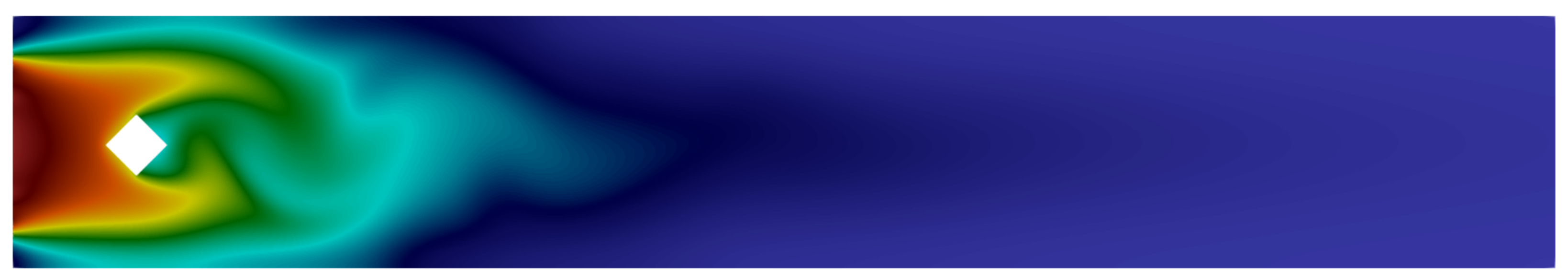}\\
		\includegraphics[width=\textwidth,trim={400 200 400 606}, clip]{figures/gnn/05_rom_\testNumb_red.pdf}    	
	\end{minipage}
	\begin{minipage}{0.325\textwidth}\centering
		Difference FOM--ROM $u$ with $\kappa = 0.05$\\
		\includegraphics[width=\textwidth,trim={10 314 10 300}, clip]{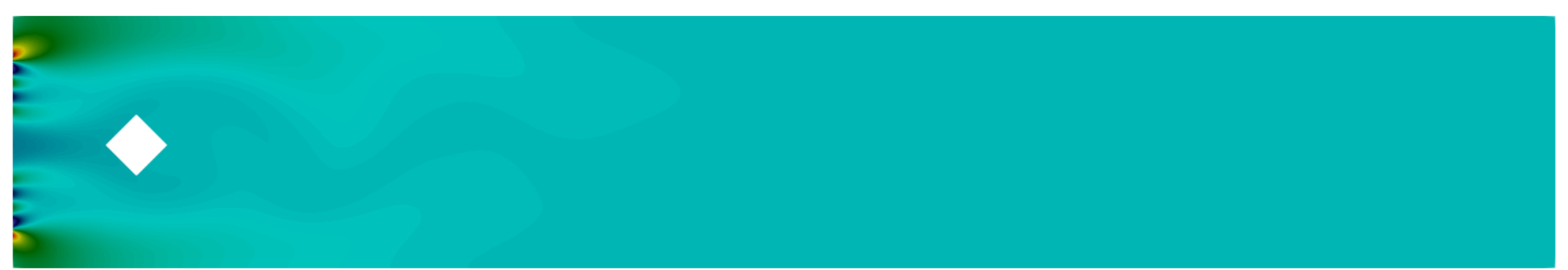}\\
		\includegraphics[width=\textwidth,trim={400 200 400 606}, clip]{figures/gnn/05_diff_\testNumb_red.pdf}    	
	\end{minipage}\\ \vspace*{2mm}
	\begin{minipage}{0.325\textwidth}\centering
		FOM $u$, $\kappa = 0.01$\\
		\includegraphics[width=\textwidth,trim={10 314 10 300}, clip]{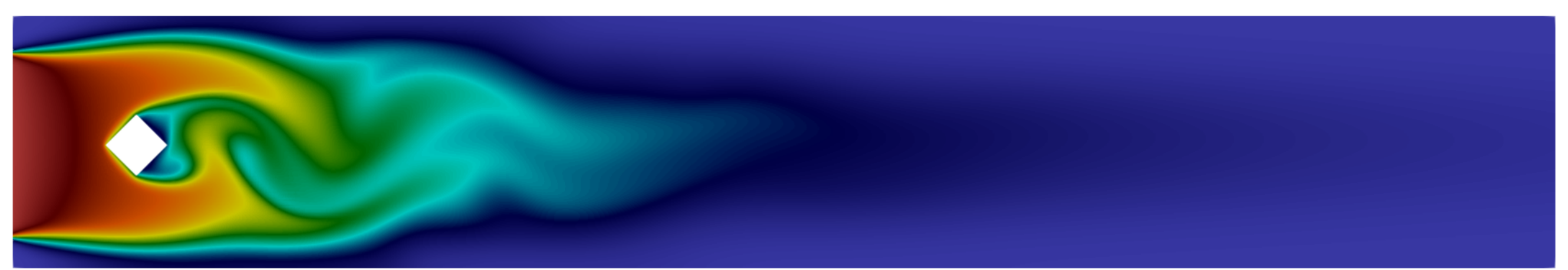}\\
		\includegraphics[width=\textwidth,trim={400 200 400 606}, clip]{figures/gnn/01_fom_\testNumb_red.pdf}    	
	\end{minipage}
	\begin{minipage}{0.325\textwidth}\centering
		ROM $u$,  $\kappa = 0.01$\\
		\includegraphics[width=\textwidth,trim={10 314 10 300}, clip]{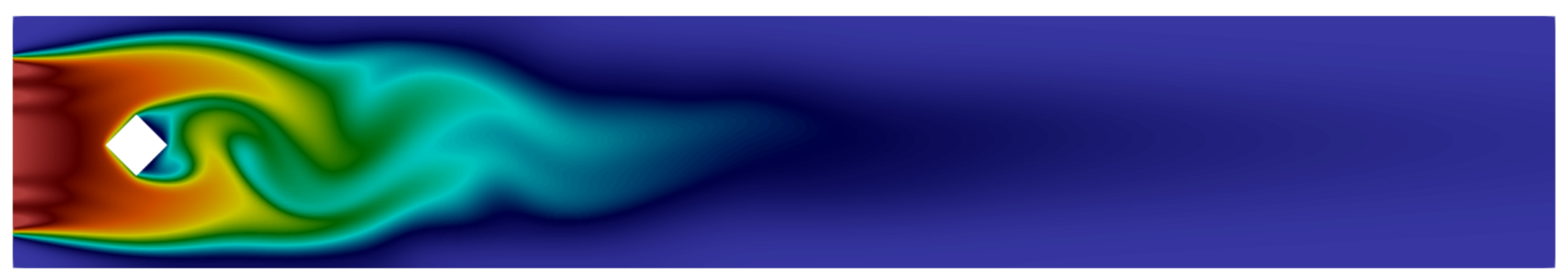}\\
		\includegraphics[width=\textwidth,trim={400 200 400 606}, clip]{figures/gnn/01_rom_\testNumb_red.pdf}    	
	\end{minipage}
	\begin{minipage}{0.325\textwidth}\centering
		Difference FOM--ROM $u$, $\kappa = 0.01$\\
		\includegraphics[width=\textwidth,trim={10 314 10 300}, clip]{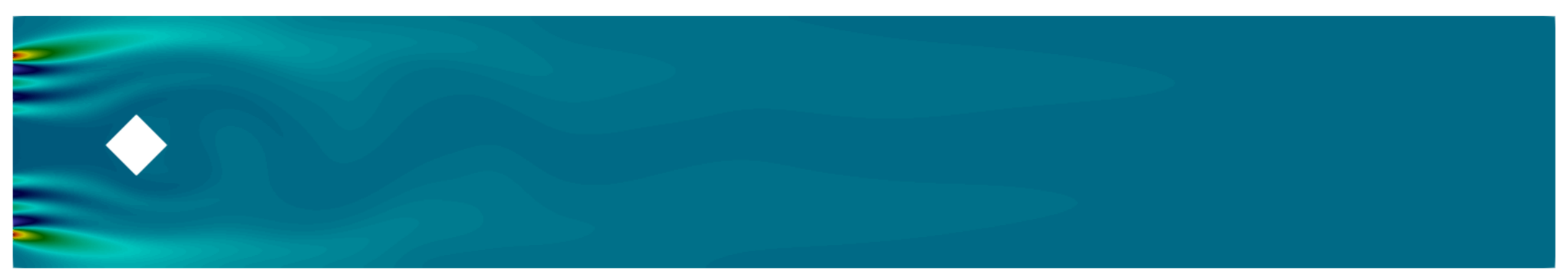}\\
		\includegraphics[width=\textwidth,trim={400 200 400 606}, clip]{figures/gnn/01_diff_\testNumb_red.pdf}    	
	\end{minipage}\\ \vspace*{2mm}    
	\begin{minipage}{0.325\textwidth}\centering
		FOM $u$, $\kappa = 0.0005$\\
		\includegraphics[width=\textwidth,trim={10 314 10 300}, clip]{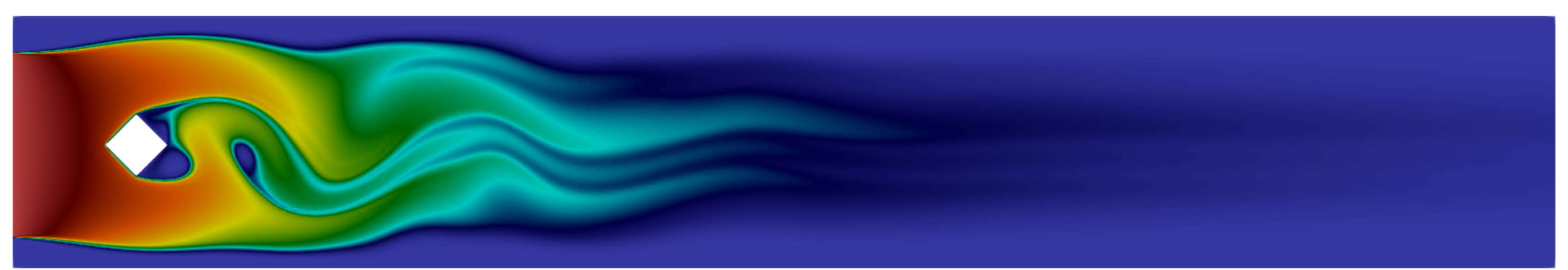}\\
		\includegraphics[width=\textwidth,trim={400 200 400 606}, clip]{figures/gnn/0005_fom_\testNumb_red.pdf}    	
	\end{minipage}
	\begin{minipage}{0.325\textwidth}\centering
		ROM $u$, $\kappa = 0.0005$\\
		\includegraphics[width=\textwidth,trim={10 314 10 300}, clip]{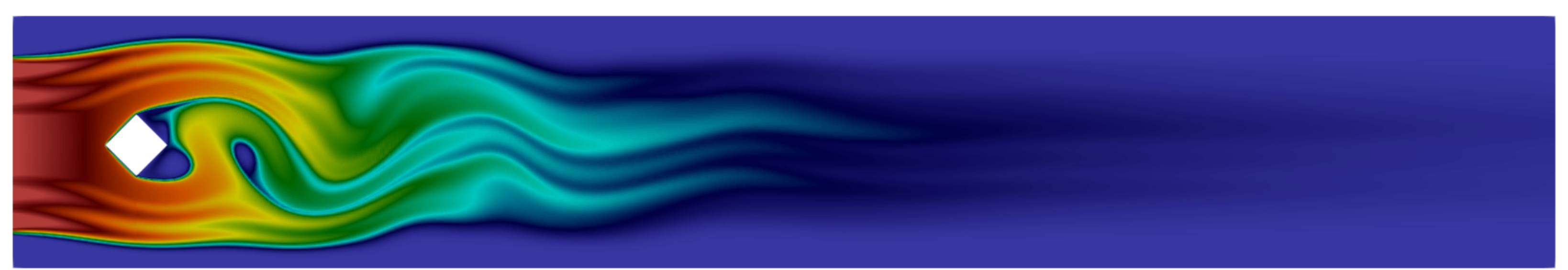}\\
		\includegraphics[width=\textwidth,trim={400 200 400 606}, clip]{figures/gnn/0005_rom_\testNumb_red.pdf}    	
	\end{minipage}
	\begin{minipage}{0.325\textwidth}\centering
		Difference FOM--ROM $u$, $\kappa = 0.0005$\\
		\includegraphics[width=\textwidth,trim={10 314 10 300}, clip]{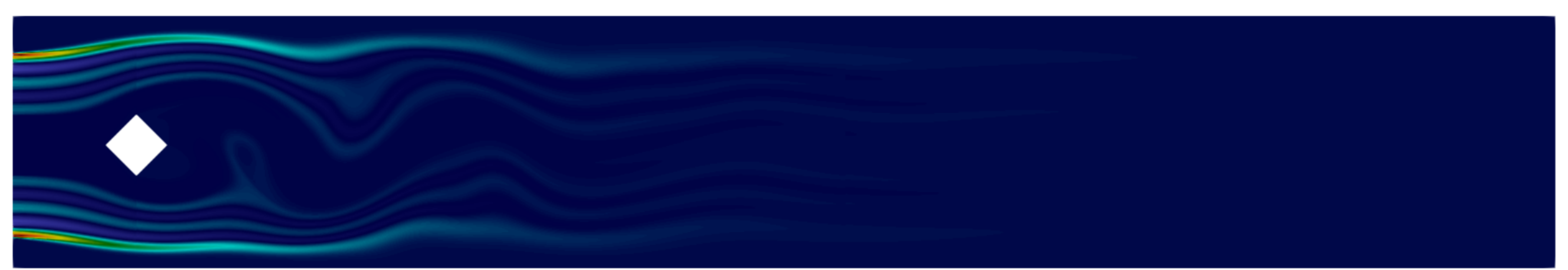}\\
		\includegraphics[width=\textwidth,trim={400 200 400 606}, clip]{figures/gnn/0005_diff_\testNumb_red.pdf}    	
	\end{minipage}\\ \vspace*{2mm}
	\begin{minipage}{0.325\textwidth}\centering
		FOM $u$, $\kappa = 0.0005$\\
		\includegraphics[width=\textwidth,trim={10 314 10 300}, clip]{figures/gnn/0005_fom_\testNumb_red.pdf}\\
		\includegraphics[width=\textwidth,trim={400 200 400 600}, clip]{figures/gnn/0005_fom_\testNumb_red.pdf}    	
	\end{minipage}
	\begin{minipage}{0.325\textwidth}\centering
		GNN $u$, $\kappa = 0.0005$\\
		\includegraphics[width=\textwidth,trim={10 314 10 300}, clip]{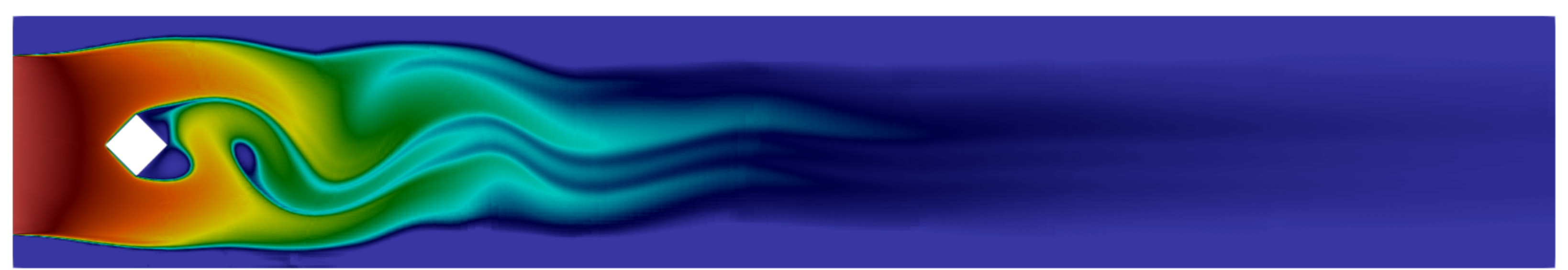}\\
		\includegraphics[width=\textwidth,trim={400 200 400 600}, clip]{figures/gnn/0005_gnn_\testNumb_red.pdf}    	
	\end{minipage}
	\begin{minipage}{0.325\textwidth}\centering
		Difference FOM--GNN $u$, $\kappa = 0.0005$\\
		\includegraphics[width=\textwidth,trim={10 314 10 300}, clip]{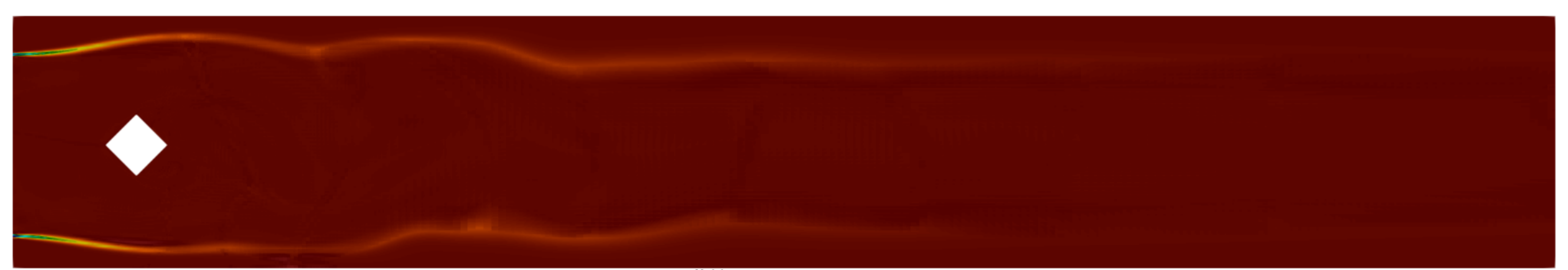}\\
		\includegraphics[width=\textwidth,trim={400 200 400 600}, clip]{figures/gnn/0005_diff_gnn_\testNumb_red.pdf}    	
	\end{minipage}
	\caption{\textbf{VV.} Scalar concentration advected by incompressible flow for ${i}=\testNumb$. Comparison of ROM approach at different viscosity levels $\kappa\in \{0.05, 0.01, 0.0005\}$ and GNN for $\kappa=0.0005$. FOMs on the left, reduced solution at the center and error on the right.}
	\label{fig:show_vv_99}
\end{figure}

Figures~\ref{fig:show_vv_0}, \ref{fig:show_vv_50} and \ref{fig:show_vv_99} show the results of the algorithm for parameters with index $i\in\lbrace 0,50,99 \rbrace$. In particular, we show on the left columns the FOM simulations, in the center column the ROM simulations and the error in the right column. 
Moreover, in the different rows, we have different viscosity levels. 
The first three rows use the classical DD-ROM approach. 
We can immediately see that the vanishing viscosity $\kappa = \kappa_3 =0.0005$ level shows strong numerical oscillations along the whole solution, which are not present in the FOM method. 
This phenomenon is observable also for higher viscosity levels but it is less pronounced and concentrated on the left of the domain, where the discontinuity are imposed as boundary conditions (see error plots).
Finally, in the last row, we show the results of the GNN approach, which uses the first two viscosity levels to predict the vanishing viscosity one.
Contrary the DD-ROM, we do not observe many numerical oscillations in the reduced solutions and they are much more physically meaningful.
Thinking about extending this approach for more complicated problems, as Euler's equations, one could guarantee the presence of the correct amount of shocks and the right location or maintaining the positivity of density and pressure close to discontinuities. 

\begin{figure}[h]
    \centering
    \includegraphics[width=1\textwidth]{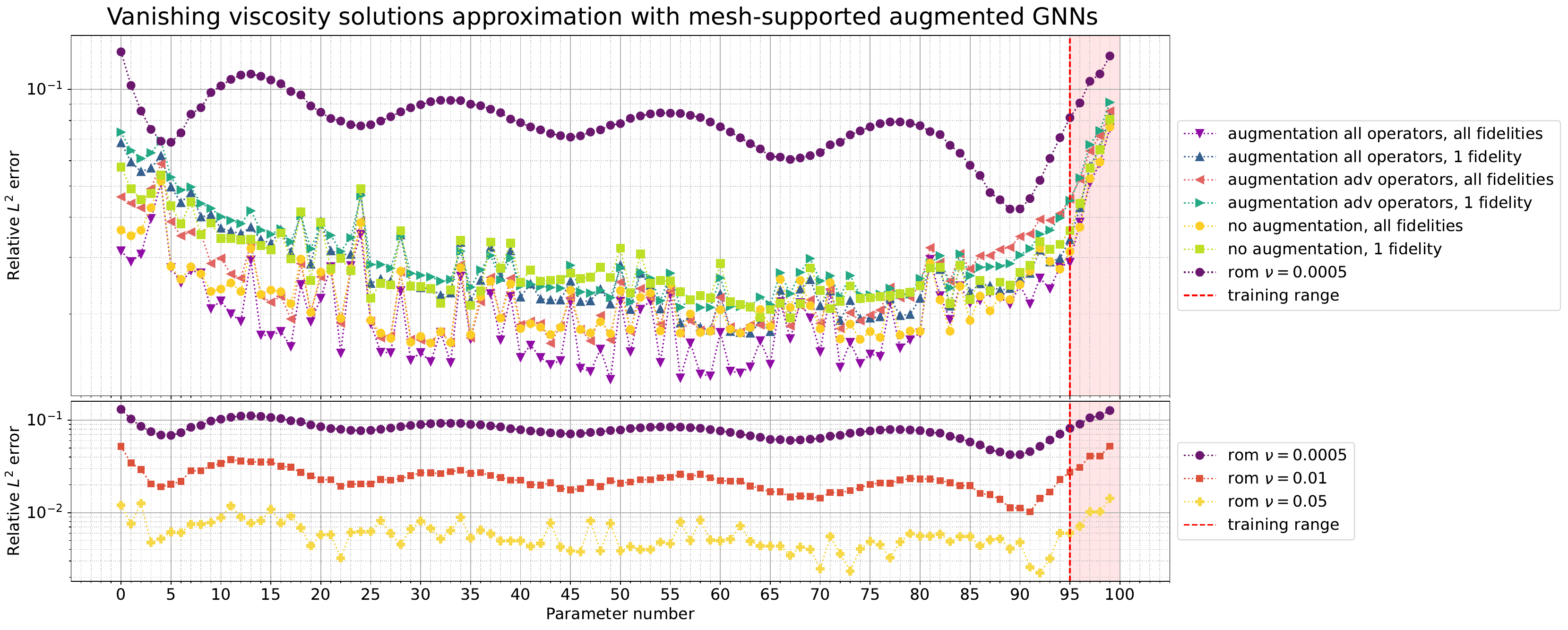}
    \caption{\textbf{VV.} Relative errors for the scalar conservation advected by incompressible flow problem. The parameters corresponding to the snapshots used for the GNNs and DD-ROMs training correspond to the abscissae $0,5,10,\dots,95$ the rest are test parameters. The dashed red background highlights the extrapolation range. \textbf{Top:} errors on train and test set with different GNN approaches given by the three augmentation $\mathcal{O}_1,\,\mathcal{O}_2$ and $\mathcal{O}_3$ and by using either 1 viscosity level (1 fidelity) or 2 (all fidelities) and errors for DD-ROM with the same viscosity level $\nu=0.0005$.  \textbf{Bottom:} errors for DD-ROM approaches at different viscosity levels. The reduced dimensions of the ROMs are $\{\NRB_{\Omega_i}\}_{i=1}^K = [5, 5, 5, 5]$ with $K=4$ partitions.}
		\label{fig:adr_error}
\end{figure}
In Figure~\ref{fig:adr_error}, we show a quantitative measure of the error of the reduced approaches presented in terms of relative $L^2$ error. 
Overall, we can immediately see that the new GNN approach can always reach errors of the order of $1-2\%$ for the vanishing viscosity solutions, with few peaks in the extrapolatory regime of $8\%$, while the classical DD-ROM on the vanishing viscosity solutions perform worse, with errors around 6-10\%.
On the other hand, the DD-ROM for higher viscosity levels have lower errors around 3\% for $\kappa_2$ and 0.5\% for $\kappa_1$, hence, they are still reliably representing those solutions.

On the different GNN approaches, in Figure~\ref{fig:adr_error} at the top we compare the different augmentations $\mathcal O_1,\,\mathcal O_2$ and $\mathcal O_3$ and how many levels of viscosity we keep into considerations to derive the vanishing viscosity solution.
The usage of multiple fidelity levels (two viscosity levels) is a great improvement for all the augmentations proposed and it can make gain a factor of 2 in terms of accuracy.
There are slight differences with the used augmentations and, in particular, we observe that the $\mathcal O_1$ augmentation, with all operators, guarantee better performance, while there are no appreciable differences between $\mathcal O_2$ and $\mathcal O_3$.
Clearly, one could come up with many other augmentation possibilities choosing more operators, but at a cost of increasing the dimensions of the GNN and the offline training costs.
We believe that all the presented options already perform much better with respect to classical approaches and can already be used without further changes.
\begin{table}
	\centering
	\caption{\textbf{VV.} Computational costs for scalar advected by a incompressible flow problem with GNNs approximating vanishing viscosity solutions (VV). The speedup is computed as the FOMs computational time over the ROM one. The speedup of the GNN is with respect to the FOM with viscosity $\nu=0.0005$. The FOM runs in parallel with $K=4$ cores as the DD-ROMs, so ``FOM time'' and ``DD-ROM time'' refers to wallclock time. Regarding the GNN results, ``Single forward GNN online time" refers to a single online evaluation while ``Total online time" refers to the evaluation of the $100$ training and test snapshots altogether with only two separate GNN forward evaluations with batches of $50$ inputs each. The speedup is evaluated as ``FOM time" over "Total online time" divided by $100$.}\label{tab:adr_GNN_times}
	\begin{tabular}{|c|c|c|c|c|c|c|}\hline
		&\multicolumn{2}{c|}{FOM}&\multicolumn{4}{c|}{DD-ROM}\\\hline 
	$\kappa$ &$N_h$ & time &  $\NRB_i$ & time &speedup & mean $L^2$ error\\ \hline\hline
	0.05 & 43776 & 3.243 [s] &  [5, 5, 5, 5] & 59.912 [$\mu$s] & ~54129& 0.00595\\ \hline
	0.01 & 43776 & 3.236 [s] &  [5, 5, 5, 5] & 79.798 [$\mu$s] & ~40552& 0.0235\\ \hline
	0.0005 & 175104 & 9.668 [s] &  [5, 5, 5, 5] & 95.844 [$\mu$s] & ~100872& 0.0796\\ \hline
	\end{tabular}
	
	\vspace{0.2cm}
	
	\begin{tabular}{|c|c|c|c|c|c|}\hline
	$\kappa$ & GNN training time & Single forward GNN online time & Total online time & GNN speedup & mean $L^2$ error \\ \hline\hline
	0.0005 & $\leq 60$ [min] & 2.661 [s] & 17.166 [s] & $\sim 56$ & 0.0217\\ \hline
\end{tabular}
\end{table}

In Table~\ref{tab:adr_GNN_times}, we compare the computational times necessary to compute the FOM solutions, the DD-ROM ones, the training time for the GNN and the online costs of the GNN. As mentioned before, we employ only one GPU NVIDIA Quadro RTX 4000 with 8GB of memory. Typical GNNs applications that involve autoencoders to perform nonlinear dimension reduction are much heavier. The training time of the GNNs for the different choice of augmentation operators vary between $48$ minutes and $60$ minutes approximately. We believe that in the near future more optimized implementations will reduce the training costs of GNNs. The computational time of the evaluation of a single forward of the GNN is on average $2.661$ seconds but vectorization ensures the evaluation of multiple online solutions altogether: with our limited memory budget we could predict all the $100$ training and test snapshots with just $2$ batches of $50$ stacked inputs each. The ``Total online time" computed as previously described is $17.166$ seconds that is $171.66$ milliseconds per online solution with a speedup of around $56$ with respect to the $9.668$ seconds for the FOM.

Although the speedup for the GNN simulations are not as remarkable as for the DD-ROM, we want to highlight that the accuracy of the GNN solutions are qualitatively much better than the DD-ROM for that viscosity level, and physically more meaningful. This aspect is a major advantage with respect to classical linear ROMs that is probably worth the loss of computational advantage. In perspective, when dealing with nonlinear and more expensive FOM for different equations, the GNN approach will not require any extra computational costs, while FOMs and ROMs model might need special treatments for the nonlinearity that would make their costs increase.

\section{Conclusions}
\label{sec:conclusions}
We argue that Friedrichs' systems represent a valuable framework to study and devise reduced order models of many parametric PDEs at the same time: among them the ones studied in this work and others, like mixed elliptic and hyperbolic problems, complex and time-dependent FS and also nonlinear PDEs whose linearization results in FS, e.g. the Euler equations. The advantages include the availability of \textit{a posteriori} error estimators and the easy to preserve mathematical properties of positivity and symmetry from the full-order formulations to the reduced-order ones. We underlined in section~\ref{subsec:ultraweak} how optimally stable reduced-order models can be obtained from the ultraweak formulation. A more efficient numerical solver for Friedrichs' systems is the hybridized discontinuous Galerkin method~\cite{chen2023unified}. These are possible future directions of research.

Working with discontinuous Galerkin discretizations is not only crucial from the possibly mixed elliptic and hyperbolic nature of Friedrichs' systems, but also to design domain decomposable reduced-order models with a minimum effort: in fact, penalties at the subdomains interfaces are inherited directly from the full-order models. We demonstrated with numerical experiments the limits and the ranges of application of domain decomposable ROMs: generally, with respect to single domain ROMs, there are benefits only when the model under study is truly decomposable, that is when the parameters affect independently different subdomains and the respective solutions are poorly correlated for unseen parametric instances. The results we showed in our academic benchmarks were obtained with the aim to tackle more complex multi-physics models like fluid-structure interaction systems. A typical application of DD-ROMs for FS is represented by parametric PDEs with a mixed elliptic and hyperbolic nature and possibly solution manifolds more and less linearly approximable respectively. The repartitioning strategies we developed are suited to adapt the reduced local dimension of the linear approximants, especially when the parameters influence only a limited region like in test case ADR~\ref{subsubsec:adr_test}. The implementation of \textit{ad hoc} physics inspired indicators can be a future direction of research.

The Friedrichs' systems formulation itself does not solve the problems caused by a slow decaying Kolmogorov n-width. DD-ROMs can help in this regard, isolating regions with a slow Kolmogorov n-width for which nonlinear approximants can be employed and regions with a fast decaying Kolmogorov n-width for which classical linear projection-based ROMs provide efficient and reliable predictions. Related to this subject and motivated also by the heavy computational resources that graph neural networks require when employed for model order reduction, we introduced a new paradigm for surrogate modeling: the inference with GNNs of vanishing viscosity solutions from a succession of higher viscosity projection-based ROMs. The approach is, of course, general and can be applied to PDEs that are not FS. The crucial hypotheses underneath this approach is the approximability with linear spaces of the solution manifolds corresponding to higher viscosity levels. We showed that the additional computational costs are not too large in our test case in~\cref{sec:vv}. Possible directions of research include more complex problems and different regularization or filtering choices, other than additional viscous terms.

\newpage
\appendix
\renewcommand\theequation{\arabic{equation}}
\section{Transformation into dissipative system}
\label{appendix:accretive}
In some cases, the term $\mathcal{A}^0=0$ or property~\eqref{eq:positivity_classical} is not satisfied, but there is a way to recover the previous framework. We want to recover a dissipative~\cite{mackenzie1995posteriori} or accretive system~\cite{jensen2004discontinuous}. For example the linearized Euler equations in entropy variables~\cite{sonar1998dual} have $\mathcal{A}^0=0$.

The condition of uniform positive definiteness
\begin{equation}
    \exists\mu_{0},\quad\mathcal{A}^{0}+(\mathcal{A}^0)^t -\mathcal{X} \geq 2\mu_0 \mathbb{I}_m\ \text{a.e. in }\Omega,
\end{equation}
is still valid if there exist $\boldsymbol{\xi}\in\mathbb{R}^d,\ \lVert\boldsymbol{\xi}\rVert = 1$ and $\beta\in\mathbb{R},\ \beta>0$ such that after the transformation
\begin{equation}
    v(x)=e^{-\beta \xi \cdot x} z ( x ),
\end{equation}
the resulting system
\begin{equation}
    \sum_i A^i \partial_i v(x) + \beta\sum_{i=1}^{d}\xi_i A^{i} v(x) = e^{-\beta(\xi \cdot x)} f,
\end{equation}
 satisfies, with the newly found $A^0 = \beta\sum_{i=1}^{d}\xi_i A^{i}$,
\begin{equation}
    \exists\mu_{0},\quad{A}^{0}+({A}^0)^t -\mathcal{X} =2\beta\sum_{i=1}^{d}\xi_i{A}^{i}-\mathcal{X}\geq 2\mu_0 \mathbb{I}_m\ \text{a.e. in }\Omega.
\end{equation}
In some cases, such ${\xi}$ and $\beta$ exist, for example if the symmetric matrix $\sum_{i=1}^{d}\xi_i{A}^{i}$ has at least one positive eigenvalue for some ${\xi}$ for almost every ${x}\in\Omega$, then taking $\beta$ sufficiently large is enough to satisfy the condition. It is also sufficient that $\sum_{i=1}^{d}\xi_i({x}){A}^{i}({x})$ has at least a positive eigenvalue for almost every ${x}\in\Omega$ where ${\xi}={\xi}({x})$, see \cite[Example 28]{jensen2004discontinuous}.

\begin{rmk}
    A more general transformation is
    \begin{equation}
        v(x)=w(x)z(x),
    \end{equation}
    so that the positive definiteness condition becomes
    \begin{equation}
        \exists\mu_{0},\quad\mathcal{A}^{0}+(\mathcal{A}^0)^t -\mathcal{X} =2\sum_{i=1}^{d}\partial_i (-\log w)\mathcal{A}^{i}-\mathcal{X}\geq 2\mu_0 \mathbb{I}_m\ \text{a.e. in }\Omega.
    \end{equation}
\end{rmk}

\section{Constructive method to define boundary operators}
\label{appendix:boundary_operator}

We report a procedure to define a boundary operator $M\in\mathcal{L}(V, V')$ starting from some specified boundary conditions. We exploit Theorem 4.3, Lemma 4.4 and Corallary 4.1 from~\cite{ern2007intrinsic}. It can be seen that the most common Dirichlet, Neumann and Robin boundary conditions can be found for some FS~\cite{ern2006discontinuous,Ern2006b,di2011mathematical}, following this procedure.

\begin{lemma}[Theorem 4.3, Lemma 4.4 and Corollary 4.1 from~\cite{ern2007intrinsic}]
    \label{lemma:constructive_M}
    Let us assume that $(V_0, V_0^*)$ satisfy~\eqref{eq: cone formalism} and that $V_0+V_0^*\subset V$ is closed. We denote with $P:V\rightarrow V_0$ and $Q:V:\rightarrow V_0^*$ the projectors onto the subspaces $V_0\subset V$ and $V_0^*\subset V$ of the Hilbert space $V$, respectively. Then, the boundary operator  $\mathcal M\in\mathcal{L}(V, V')$ defined as
    \begin{align}
    	\begin{split}
        \scpr{\mathcal Mu}{v}_{V', V} =& \scpr{\mathcal DPu}{Pv}_{V', V} - \scpr{\mathcal DQu}{Qv}_{V', V} +\\
        &\scpr{\mathcal D(P+Q-PQ)u}{ v}_{V', V} - \scpr{\mathcal Du}{(P+Q-PQ)v}_{V', V}
    	\end{split}
    \end{align}
    is admissible and satisfies $V_0=\ker(\mathcal D-\mathcal M)$ and $V_0^*=\ker(\mathcal D+\mathcal M^*)$. In particular, 
    \begin{enumerate}
        \item If $V = V_0 + V_0^*$, then $\mathcal M$ is self-adjoint and
        \begin{equation}
            \scpr{\mathcal Mu}{v}_{V', V} = \scpr{\mathcal DPu}{Pv}_{V', V} - \scpr{\mathcal DQu}{ Qv}_{V', V}.
        \end{equation}
        \item If $V_0=V_0^*$, then $\mathcal M$ is skew-symmetric and
        \begin{equation}
            \label{eq:constructive_M}
            \scpr{\mathcal Mu}{v}_{V', V} = \scpr{\mathcal DPu}{v}_{V', V} - \scpr{\mathcal DPv}{ u}_{V', V}.
        \end{equation}
    \end{enumerate}
\end{lemma}

We remark that, for fixed $(V_0, V_0^*)$, admissible boundary operators $\mathcal M\in\mathcal{L}(V, V')$ that satisfy $V_0=\ker(\mathcal D-\mathcal M)$ and $V_0^*=\ker(\mathcal D+\mathcal M^*)$ are not unique. The boundary operator defined in Lemma~\ref{lemma:constructive_M} is just a possible explicit definition, in general.

As an exercise, we show how to find the definition of the operator $\mathcal M$ for our linear compressible elasticity FS, from Section~\ref{subsubsec:comprLinElast_def}. We want to impose the boundary conditions $\mathbf{u}_{|\Gamma_D} = 0$ and $(\boldsymbol{\sigma}\cdot\mathbf{n})_{|\Gamma_N}=0$, so,
\begin{equation}
    V_0=V_0^* = \{(\mathbf{u}, \boldsymbol{\sigma})\in V\mid \mathbf{u}_{|\Gamma_D} = 0,\quad (\boldsymbol{\sigma}\cdot\mathbf{n})_{|\Gamma_N}=0\} = H_{\boldsymbol{\sigma},\Gamma_N}\times[H^1_{\Gamma_D}(\Omega)]^d,
\end{equation}
since we defined $V=H_{\boldsymbol{\sigma}}\times [H^1(\Omega)]^d$, with $H_{\boldsymbol{\sigma}}=\{\boldsymbol{\sigma}\in [L^2(\Omega)]^{d\times d}\mid \nabla\cdot(\boldsymbol{\sigma}+\boldsymbol{\sigma}^t)\in [L^2(\Omega)]^d\}$, the traces $\gamma_D:[H^1(\Omega)]^d\rightarrow [H^{\frac{1}{2}}(\Gamma_D)]^d$ and $\gamma_N:H_{\boldsymbol{\sigma}}^d\rightarrow [H^{-\frac{1}{2}}(\Gamma_N)]^d$ on $\Gamma_D$ and $\Gamma_N$ are well-defined. In particular,
\begin{equation}
    H_{\boldsymbol{\sigma},\Gamma_N}=\{\boldsymbol{\sigma}\in H_{\boldsymbol{\sigma}}\mid \gamma_N(\boldsymbol{\sigma})=(\boldsymbol{\sigma}\cdot\mathbf{n})_{|\Gamma_N}=0\},\qquad
    [H^1_{\Gamma_D}(\Omega)]^d = \{\mathbf{u}\in [H^1(\Omega)]^d\mid \gamma_D(\mathbf{u})=\mathbf{u}_{|\Gamma_D} = 0\}.
\end{equation}
Moreover, $(V_0, V_0^*)$ satisfy the properties of cone formalism~\eqref{eq: cone formalism}. Thus, we can use the definition~\eqref{eq:constructive_M} of Lemma~\ref{lemma:constructive_M}:
\begin{align}
\begin{split}
    \scpr{\mathcal M(\boldsymbol{\sigma}, \mathbf{u})}{(\boldsymbol{\tau}, \mathbf{v})}_{V', V} =& \scpr{\mathcal DP(\boldsymbol{\sigma}, \mathbf{u})}{ (\boldsymbol{\tau}, \mathbf{v})}_{V', V} - \scpr{\mathcal DP(\boldsymbol{\tau}, \mathbf{v})} {(\boldsymbol{\sigma}, \mathbf{u})}_{V', V}\\
    =&-\langle \tfrac{1}{2}(\boldsymbol{\sigma}+\boldsymbol{\sigma}^t)\cdot \n, \mathbf{v}\rangle_{-\frac{1}{2}, \frac{1}{2}, \Gamma_D}
	+ \langle  \tfrac{1}{2}(\boldsymbol{\tau}+\boldsymbol{\tau}^t)\cdot \n, \mathbf{u}\rangle_{-\frac{1}{2}, \frac{1}{2}, \Gamma_D} +\\
	&\langle \tfrac{1}{2}(\boldsymbol{\sigma}+\boldsymbol{\sigma}^t)\cdot \n, \mathbf{v}\rangle_{-\frac{1}{2}, \frac{1}{2}, \Gamma_N}
	- \langle  \tfrac{1}{2}(\boldsymbol{\tau}+\boldsymbol{\tau}^t)\cdot \n, \mathbf{u}\rangle_{-\frac{1}{2}, \frac{1}{2}, \Gamma_N},
\end{split}
\end{align}
where $P:V=H_{\boldsymbol{\sigma}}\times [H^1(\Omega)]^d\rightarrow V_0=H_{\boldsymbol{\sigma},\Gamma_N}\times[H^1_{\Gamma_D}(\Omega)]^d$ is the projector into the subspace $V_0$ of the Hilbert graph space $V$ with scalar product:
\begin{equation}
    ((\boldsymbol{\sigma}, \mathbf{u}), (\boldsymbol{\tau}, \mathbf{v}))_{V} = (\mathbf{u}, \mathbf{v})_{[L^d(\Omega)]^d} + (\boldsymbol{\sigma}, \boldsymbol{\tau})_{[L^2(\Omega)]^{d\times d}} + (A(\boldsymbol{\sigma}, \mathbf{u}), A(\boldsymbol{\tau}, \mathbf{v})).
\end{equation}

\section{ROM convergence studies}\label{app:convergence_ROM}
In this section, we validate the DD-ROM implementation, checking the convergence towards the FOM solutions with respect to the dimension of the reduced space. Uniform local reduced dimensions are employed $\{r_{\Omega_i}\}_{i=1}^{K}$ and $\{r_{\Omega_i}\}_{i=1}^{k}$. For each convergence study $20$ uniformly independent samples are used as training dataset and $50$ uniformly independent samples as test dataset.
\begin{figure}[!htpb]
    \centering    
    \includegraphics[width=1\textwidth]{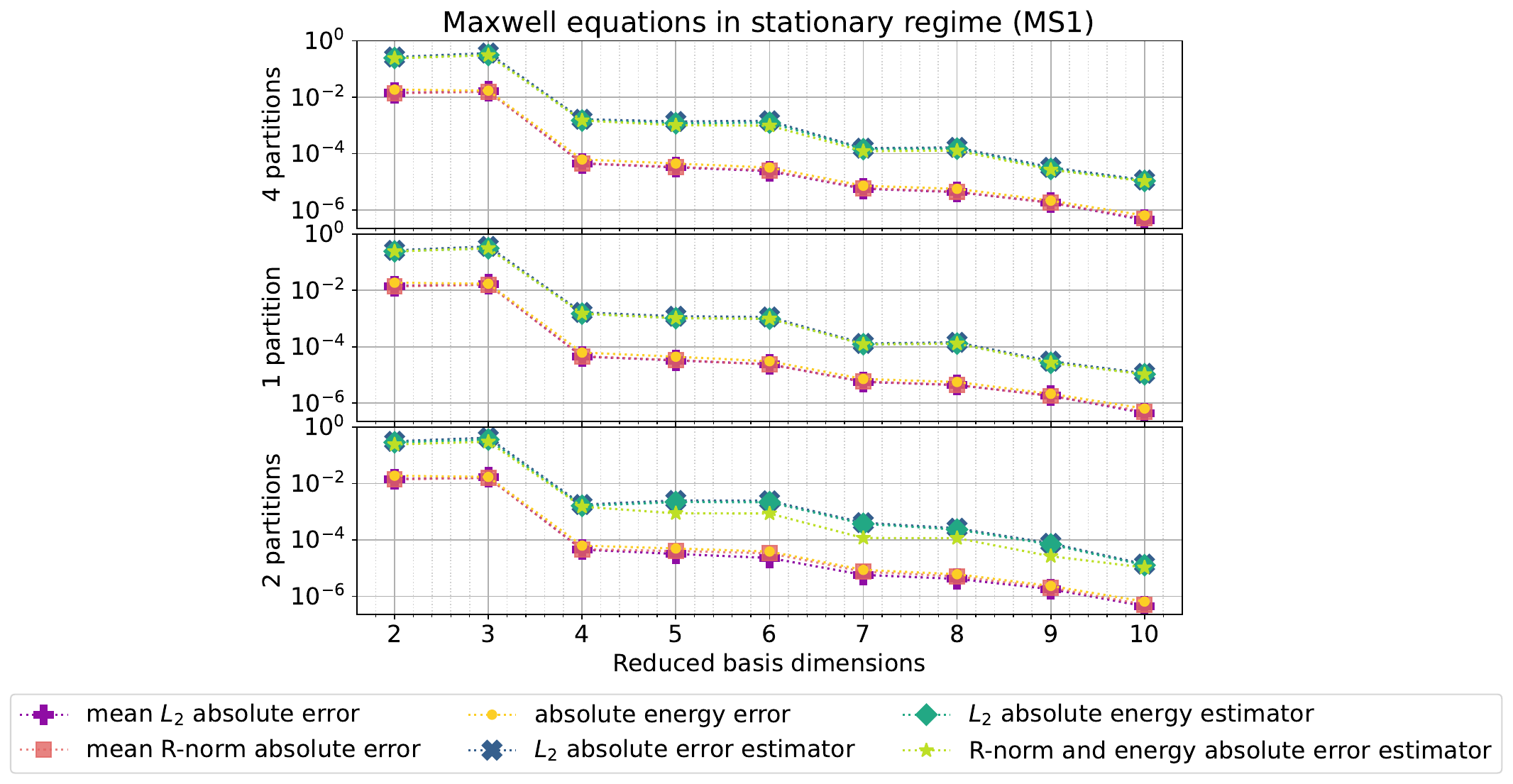}
    \caption{\textbf{MS1}. The convergence of DD-ROMS with uniform local reduced dimensions $\{r_{\Omega_i}\}_{i=1}^{K}$ and $\{r_{\Omega_i}\}_{i=1}^{k}$ is assessed. The uniform value of the local reduced dimensions is reported in the abscissae. For this test case an improvement of the accuracy with respect to the single domain reduced basis is not observed.}
    \label{fig:conv_torus}
\end{figure}

\begin{figure}[!htpb]
    \centering    
    \includegraphics[width=1\textwidth]{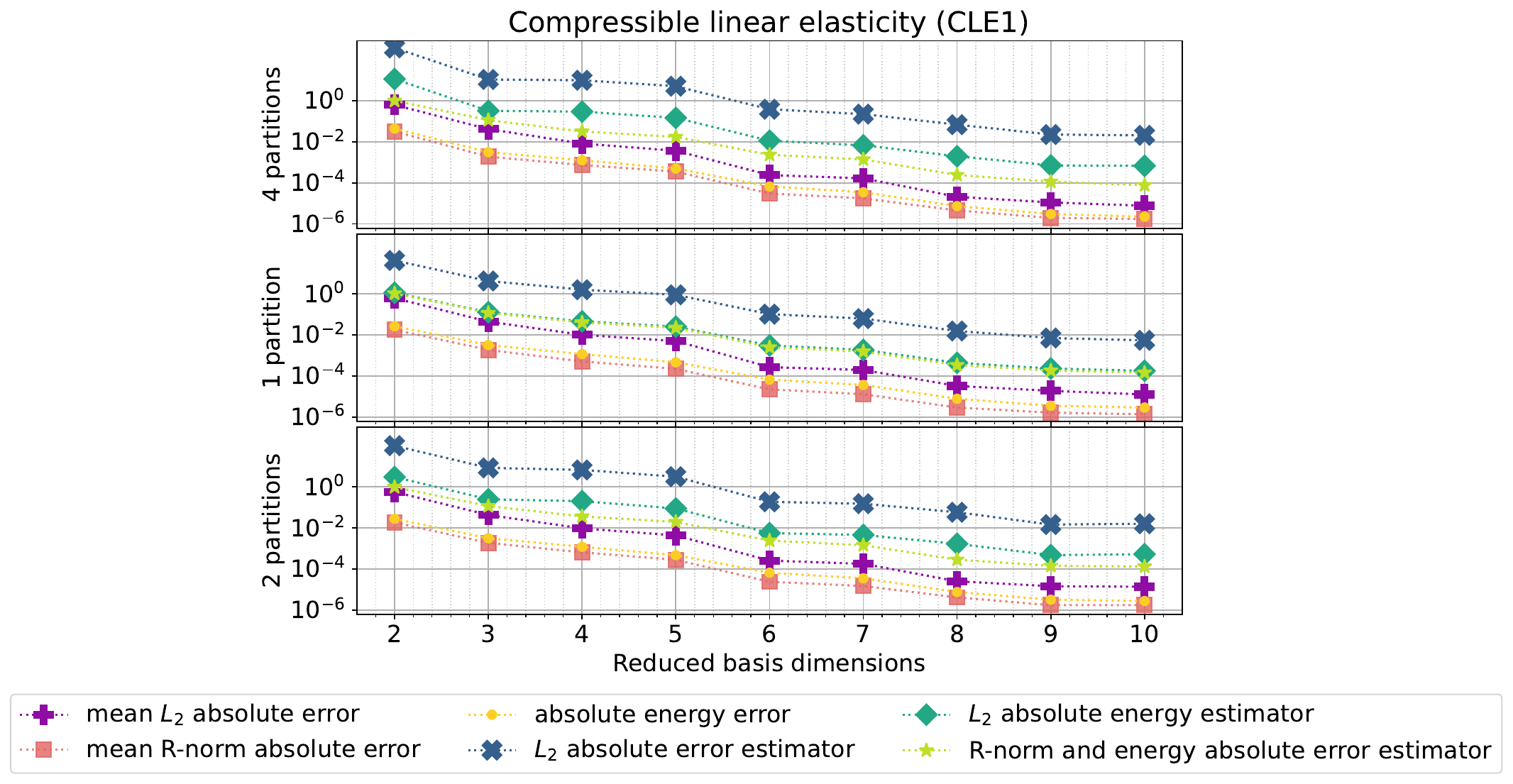}
    \caption{\textbf{CLE1}. The convergence of DD-ROMS with uniform local reduced dimensions $\{r_{\Omega_i}\}_{i=1}^{K}$ and $\{r_{\Omega_i}\}_{i=1}^{k}$ is assessed. The uniform value of the local reduced dimensions is reported in the abscissae. For this test case an improvement of the accuracy with respect to the single domain reduced basis is not observed.}
    \label{fig:conv_bar}
\end{figure}

In Figure~\ref{fig:conv_torus}, we show the $L^2$-error, the $R$-error and the energy error decay and their respective error estimators for the Maxwell equations test case \textbf{MS1}, ~\cref{subsubsec:max}, with constant parameters $\mu$ and $\sigma$ on the whole domain. We clearly see an exponential behavior in the error as we add basis functions. On the other hand, we do not observe strong differences between the ROM, DD-ROM with repartition and DD-ROM with \texttt{deal.II} subdomains, for this simple test case. Similar results can be observed in Figure~\ref{fig:conv_bar}, where the same analysis is applied for the compressible linear elasticity test \textbf{CLE1} from~\cref{subsubsection:cle_test}. 

From these results, it should be clear that the employment of local reduced basis is not always useful to increase the accuracy of the predictions. Nonetheless, it may be used to locally reduce the dimension of the linear approximants. Possible benefits include the adaptation of the computational resources (higher dimensional reduced basis are chosen only where it is necessary) and the possibility to speedup parametric studies and non-intrusive surrogate modelling thanks to the further reduced local dimensions~\cite{xiao2019domain,xiao2019domainb}.

Typical cases where DD-ROMs are effective to increase the accuracy of the predictions are truly decomposable systems where the parameters affect independently different regions of the computational domain, as in test case \textbf{MS2} in~\cref{subsubsec:max}.

\section*{Acknowledgements}
This work was partially funded by European Union Funding for Research and
Innovation --- Horizon 2020 Program --- in the framework of European Research
Council Executive Agency: H2020 ERC CoG 2015 AROMA-CFD project 681447 ``Advanced
Reduced Order Methods with Applications in Computational Fluid Dynamics'' P.I.
Professor Gianluigi Rozza. We also acknowledge the PRIN 2017 “Numerical Analysis for Full and Reduced Order Methods for the efficient and accurate solution of complex
systems governed by Partial Differential Equations” (NA-FROM-PDEs). Davide Torlo  has been funded by a SISSA Mathematical fellowship within Italian Excellence Departments initiative by Ministry of University and Research.

\bibliographystyle{abbrv}
\bibliography{biblio}

% \begin{thebibliography}{10} \end{thebibliography}

\end{document}